\documentclass[a4paper, 10pt]{article}
\usepackage[active]{srcltx}

\usepackage{multicol}
\usepackage{euscript}
\usepackage{amssymb}
\usepackage{amsmath}
\usepackage{amsthm}
\usepackage[dvips]{graphicx}
\usepackage{eufrak}
\usepackage{psfrag}
\newcommand{\integraltt}{\int_{\max \{ 2, \, 2 + \tau \}}^{\min \{ T, \, T + \tau \}}}
\newcommand{\integralt}{\int_{\max \{ 0, \, \tau \}}^{\min \{ T, \, T + \tau \}}}
\newcommand{\integralx}{\int_{ \max \{0, \, \xi \}}^{ \min \{ L, \, L + \xi \}}}
\newcommand{\semi}{p_t [\bar{u}_0, \, \bar{u}_{b \, 0},
                   \, \bar{u}_{b \, l}]}
\newcommand{\semie}{p^{\; \ee}_t \,  [\bar{u}_0, \, \bar{u}_{b \, 0},
                   \, \bar{u}_{b \, l}]}

\newcommand{\bigd}{\mathrm{D}}
\newcommand{\hr}{\hat{r}_{1}}

\newcommand{\uu}{\Upsilon}

\newcommand{\driftd}{\lambda_2^{\ast}}

\newcommand{\ee}{\varepsilon}
\newcommand{\ue}{u^{\varepsilon}}
\newcommand{\integral}{\int_{0}^{L} \int_{0}^{T}}
\newcommand{\energy}{\frac{\, v_1 ^2}{2}}
\newcommand{\argument}{\frac{w_1}{v_1} + \lambda_1^{\ast}}

\newcommand{\bv}{\mathrm{Tot \, Var}}

\newcommand{\unpo}{\mathcal{O}(1)}
\newcommand{\Jzero}{  J^{ \lambda_i^{\ast} \, 0}}
\newcommand{\Jelle}{  J^{ \lambda_i^{\ast} \, L}}
\newcommand{\deltabt}{\tilde{\Delta}^{\lambda_i^{\ast}}}
\newcommand{\deltab}{\Delta^{\lambda_i^{\ast}}}

\newcommand{\tr}{\tilde{r}_1}
\newcommand{\domain}{\mathcal{D}_0}
\newcommand{\domainb}{\mathcal{D}_{ \, b}}
\newcommand{\domainu}{\mathcal{U}_0}
\newcommand{\domainub}{\mathcal{U}_{\, b}}
\newtheorem{teo}{Theorem}[section]
\newtheorem{pro}{Proposition}[section]
\newtheorem{lem}{Lemma}[section]
\theoremstyle{definition}
\newtheorem{say}{Definition}[section]
\newtheorem{rem}{Remark}[section]

\begin{document}
\numberwithin{equation}{section}
\bibliographystyle{plain}

\title{Vanishing viscosity solutions of a
 $2 \times 2$ triangular hyperbolic system with Dirichlet conditions on two
boundaries}
\author{Laura V. Spinolo  \smallskip  \\
        { \it \small S.I.S.S.A., Via Beirut 4, Trieste 34014, Italy} \\
             {\it \small e-mail: spinolo@sissa.it }}
\date{}


\maketitle

\begin{abstract}
We consider the $2 \times 2$ parabolic systems
\begin{equation*}
      u^{\ee}_t + A(u^{\ee}) u^{\ee}_x = \ee u^{\ee}_{xx}
\end{equation*}
on a domain $(t, \, x) \in \, ]0, \, + \infty[ \times ]0, \, l[$
with Dirichlet boundary conditions imposed at $x=0$ and at $x=l$.
The matrix $A$ is assumed to be in triangular form and strictly
hyperbolic, and the boundary is not characteristic, i.e. the
eigenvalues of $A$ are different from $0$.

We show that, if the initial and boundary data have sufficiently
small total variation, then the solution $u^{\ee}$ exists for all
$t \geq 0$ and depends Lipschitz continuously in $L^1$ on the
initial and boundary data.

Moreover, as $\ee \to 0^+$, the solutions $u^{\ee}(t)$ converge in
$L^1$ to a unique limit $u(t)$, which can be seen as the {\it
vanishing viscosity solution} of the quasilinear hyperbolic system
\begin{equation*}
   u_t + A(u)u_x = 0,
   \quad x \in \, ]0, \, l[.
\end{equation*}
This solution $u(t)$ depends Lipschitz continuously in $L^1$ w.r.t
the initial and boundary data. We also characterize precisely in
which sense the boundary data are assumed by the solution of the
hyperbolic system. \vspace{0.3cm}

\vspace{0.3cm}


\end{abstract}

\section{Introduction}

This paper deals with the initial-two-boundaries value problem
\begin{equation}
\label{the_problem}
      \left\{
      \begin{array}{lllll}
            u_t + A (u) u_x =0 ,
            \quad
            x \in \, ]0, \, l[, \; \;
            t \in \, ]0, + \infty [  \\
            \\
            u(0, x) = \bar{u}_0 (x),\\
            \\
            u(t, 0) = \bar{u}_{b \, 0}(t), \qquad
            u(t, l) = \bar{u}_{b l}(t). \\
     \end{array}
     \right.
\end{equation}
The crucial hypotheses we assume are that the matrix $A$ is
strictly hyperbolic with eigenvalues different from $0$ and that
the initial and boundary data are small in $BV$ norm and close to
a constant state $u^\ast$.

An existence result for hyperbolic boundary value problems was
proved in \cite{Good, SabTou:mixte} using an adaptation of the
Glimm scheme introduced in \cite{Gli}. Improvements of the results
in \cite{Good, SabTou:mixte} have been obtained by a wave-front
tracking technique introduced in \cite{Breft}  and later used in a
series of papers (\cite{Bre:Gli, BreCol, BreCraPi,BreLY,BreLF,
BreG, BreLew}) to establish the well posedness of the Cauchy
problem. Such a wave-front tracking technique was adapted to the
initial-boundary value problem in \cite{Ama}, where a substantial
improvement of the results in \cite{Good, SabTou:mixte} was
achieved. The well posedness of the initial-boundary value problem
was then proved in \cite{DonMar} relying on the wave-front
tracking technique described in \cite{Ama}.

All the results quoted so far deal with conservative systems; a
comprehensive account of the stability and uniqueness results for
the Cauchy problem for a system of conservation laws can be found
in \cite{Bre}. We refer, instead, to \cite{Daf:book} and to
\cite{Serre:book} for a general introduction to the systems of
conservation laws.

In \cite{BiaBre:BV, BiaBre:case, BiaBre:center} and
\cite{BiaBrevv} a different problem was dealt with: let $u^{\ee}$
be a family of solutions to the parabolic systems
\begin{equation*}
      u^{\ee}_t + A(u^{\ee})u^{\ee}_x = \ee u^{\ee}_{xx}.
\end{equation*}
One expects that as $\ee \to 0^+$ the
solution $u^{\ee}$ converges in some sense to a solution of the
corresponding hyperbolic system
\begin{equation*}
      u_t + A(u)u_x=0.
\end{equation*}
The mathematical proof of this convergence was obtained via a
suitable decomposition of the gradient of the solution $u^{\ee}$
along travelling waves. We refer to \cite{BiaBrevv} for an account
of the proof of the convergence of the vanishing viscosity
approximation and of the uniqueness and the stability of the
vanishing viscosity limit: it is important to underline, however,
that in \cite{BiaBrevv} the systems considered are not necessarily
conservative.

The vanishing viscosity approximation of initial-boundary value
problems was studied in numerous works: in the following, we will
briefly refer to some of the principal results, without any sake
of completeness. Moreover, if not otherwise stated, the systems
considered are supposed to be in conservation form.

In particular, in \cite{Serre:solglob} it was considered the
vanishing viscosity approximation
\begin{equation*}
        \ue_t + f(\ue )_x = \ee \ue_{xx}
\end{equation*}
of an initial-boundary value problem and it was given a precise
description of the first term of the expansion of $\ue$ in the
neighborhood of a point where two shocks or a shock and a boundary
layer profile meet.

The works \cite{Gisclon:etudes, Gisclon-Serre:etudes} dealt with
the general parabolic approximation
\begin{equation}
\label{eq_approx}
        \ue_t + f(\ue)_x = \ee \big( B(\ue ) \ue_x \big)_x,
\end{equation}
where the viscosity $B(u)$ is invertible but in general different
from the identity. It was proved the existence of a $T> 0$ such
that $\ue$ converges in $L^{\infty}\big( (0, \, T);
L^2(\mathbb{R}^+) \big)$ to a solution of
\begin{equation*}
       u_t + f(u)_x =0
\end{equation*}
and it is given a precise characterization of the boundary
condition induced in the hyperbolic limit.

In \cite{Serre:couches} it was introduced an Evans function
machinery to study the stability of boundary layer profiles: the
parabolic approximation considered was in the form
\eqref{eq_approx},
in the case of an invertible viscosity matrix $B$ and of a non
characteristic boundary (i.e. all the eigenvalues of $Df(u)$ were
supposed to be different from zero). However, the analysis was
extended in a series of paper (\cite{Serre-Zum, Rousset:inviscid,
Rousset:residual, Rousset:char}) to the boundary characteristic
case and to very general parabolic approximations, with non
invertible viscosity matrices.

In \cite{AnBia} it was considered the family of
initial-one-boundary value problems
\begin{equation*}
      \left\{
      \begin{array}{lll}
            u^{\varepsilon}_t + A (u^{\varepsilon}) u^{\varepsilon}_x =
            \varepsilon u^{\varepsilon}_{xx},
            \quad
            x \in \, ]0, \, + \infty[,  \quad
            t \in \, ]0, + \infty [ \\
            \\
            u^{\varepsilon}(0, x) = \bar{u}_0 (x),
            \qquad
            u^{\varepsilon}(t, 0) = \bar{u}_b(t), \\
     \end{array}
     \right.
\end{equation*}
it is proved the (global in time) convergence of approximated
solutions and the stability and the uniqueness of the limit. In
\cite{AnBia} the boundary characteristic case was allowed (i.e.
one characteristic field was allowed to have speed close to that
of the boundary) and the crucial tool in the proof of the
convergence and the stability is the introduction of a suitable
decomposition of the gradient of the vanishing viscosity solution.
Moreover, we underline that, as in \cite{BiaBrevv}, the systems
considered were not necessarily in conservation form.

In the present paper we will consider the vanishing viscosity
approximation for the initial-two-boundaries value problem:
\begin{equation}
\label{vvapproximation}
      \left\{
      \begin{array}{lllll}
            u^{\varepsilon}_t + A (u^{\varepsilon}) u^{\varepsilon}_x =
            \varepsilon u^{\varepsilon}_{xx},
            \quad
            x \in \, ]0, \, l[,  \quad
            t \in \, ]0, + \infty [ \\
            \\
            u^{\varepsilon}(0, x) = \bar{u}_0 (x),\\
            \\
            u^{\varepsilon}(t, 0) = \bar{u}_{b \, 0}(t), \qquad
            u^{\varepsilon}(t, l) = \bar{u}_{b l}(t). \\
     \end{array}
     \right.
\end{equation}
We will assume that $A$ is in
triangular form, i.e.
\begin{equation}
\label{E:triang}
       A(u) =
       \left(
       \begin{array}{lll}
             \lambda_1(u_1) & 0 \\
             \\
             g(u_1, \, u_2) & \lambda_2(u_1, \, u_2) \\
       \end{array}
       \right),
\end{equation}
and sufficiently smooth in a compact neighborhood $K$ of a fixed
point $u^{\ast}$. Moreover, we assume $A$ to be uniformly strictly
hyperbolic, in particular we assume that there exists a constant
$c>0$ ($2c$ is then the "separation speed") such that
\begin{equation}
\label{eq_separation_speed}
       \lambda_1(u) < - c < 0 < c < \lambda_2 (u) \quad \forall \, u
       \in K.
\end{equation}
The above condition means that the speed of the boundary (in our
case $0$) is strictly different from the characteristic speeds of
the two families of waves.

We denote with $r_1(u)$ the first eigenvector of $A(u)$,
corresponding to the eigenvalues $\lambda_1(u)$, and with $r_2$
the second one. Due to the particular structure of $A$, we
normalize $r_1$ and $r_2$ as
\begin{equation}
\label{E:norm}
\langle (1,0), r_1(u) \rangle = 1, \quad r_2 = \left( \begin{array}{c}
0 \cr 1
\end{array} \right).
\end{equation}
The dual base of $(r_1(u), \, r_2)$ is denoted by $(\ell_1, \,
\ell_2(u))$.

We will assume that the initial data $\bar u_0$ and boundary
data $\bar u_{b0}$, $\bar u_{bl}$ have sufficiently small total variation, i.e.
\begin{equation}
\label{E:TVbd}
      \bv(\bar{u}_0), \; \bv(\bar{u}_{b \, 0}), \;
      \bv(\bar{u}_{b \, l}) \leq
      \delta_1
\end{equation}
for a suitable $\delta_1<<1$. Moreover, since we will study
boundary layers with small total variation, we assume that there
exists a value $u^{\ast}$ such that
\begin{equation} \label{E:noBD}
       \| \bar{u}_0 - u^{\ast} \|_{\infty} \leq \delta_1
       \quad
       \| \bar{u}_{b \, 0} - u^{\ast} \|_{\infty} \leq \delta_1
       \quad
       \| \bar{u}_{b \, l} - u^{\ast} \|_{\infty} \leq \delta_1.
\end{equation}
For technical reasons, we will also
assume some stronger regularity: the boundary and initial data
will be sufficiently smooth and will satisfy
\begin{equation}
\label{E:bounder} \|d^j \bar{u}_0 / dx^j\|_{L^1}, \; \|d^j
\bar{u}_{b \, 0} /
      dt^j\|_{L^1}, \; \|d^j \bar{u}_{b \, l} / dt^j\|_{L^1}
      \leq M < + \infty \quad j= 2, \, \dots n,
\end{equation}
for some $n \in \mathbb{N}$ and some large constant $M$. Some
observations about the extension of our results to the case of
boundary and initial data with weaker regularity will be made in
Remark \ref{rem_regularity}.

We will denote by $\mathcal{U}_{\, 0}$, $\mathcal{U}_{\, b}$ the
set of functions $u_0$, $u_b$ satisfying \eqref{E:TVbd},
 \eqref{E:noBD}, \eqref{E:bounder} in $]0,l[$ or $]0,+\infty[$,
respectively. We also define the sets $\domain \subseteq L^1(0, \,
l)$, $\domainb \subseteq L^1_{loc}(0, \, + \infty)$ of functions
such that
\begin{equation}
\label{eq_domain}
      \bv \big\{ \bar{u}_0 \big\} \leq \delta_1, \quad \bv
     \big\{ \bar{u}_b \big\} \leq \delta_1,
\end{equation}
respectively.
\begin{rem}
The fact that we will consider only $2 \times 2$ triangular
systems does not affect very deeply the structure of the problem,
but leads to some considerable simplification in the computations.
In particular, since the matrix $A$ is in triangular form, we will
see in Section \ref{gradient_decomposition} that the generalized
eigenvector of the travelling wave profile of the second family is
constant, and so it is the generalized eigenvector of the boundary
layer profile of the second family: such a feature simplifies the
computation of source terms, which is performed in the Appendix
\ref{explicit_source_t}. Since also the expression itself of the
source terms is simpler, the consequent estimates, carried on in
Section \ref{BV_estimates}, are easier in the case of a triangular
system than in the general one.

We refer, instead, to Remark \ref{rem_regularity} for some
considerations about the hypotheses of regularity we have assumed.
\end{rem}

The first theorem concerns the existence of a solution to the
parabolic problem \eqref{vvapproximation}; moreover, it ensures
that such a solution satisfies stability estimates independent on
$\ee$.

\begin{teo}
\label{main_result}
      Suppose $\bar{u}_0 \in \domainu$, $\bar{u}_{b \, 0}, \; \bar{u}_{b \, l} \in
      \domainub$ and $A$ is of the form \eqref{E:triang}
      and satisfies \eqref{eq_separation_speed}. Then, for any $\ee> 0$, the system
      \eqref{vvapproximation} has a unique solution
      $u^{\ee}(t)$ defined for all $t \ge 0$.

      This solution depends Lipschitz continuously
      in $L^1$ on the initial and boundary data: indeed, let
      $\bar{v}_0 \in \domainu$, $\bar{v}_{b \, 0}$, $\bar{v}_{b \, l} \in \domainub$
      be the initial and boundary data of a solution
      $v^{\ee}(t)$ of \eqref{vvapproximation}. Then for some
      constants $L_1$ and $L_2$,
      depending only on the matrix $A$ and the bound on the initial and boundary
      data $\delta_1$,
      the following holds:
      \begin{equation}
      \label{stability_wrt_l1_tpar}
      \begin{split}
              \|v^\ee(t) - u^\ee(t)\|_{L^1}
      &       \leq
              L_1 \Big(  \|\bar{v}_{ 0} - \bar{u}_{0}\|_{L^1(0, \, l)}+
              \|\bar{v}_{b 0} - \bar{u}_{b 0}\|_{L^1(0, \, + \infty)}
              +  \|\bar{v}_{b l} - \bar{u}_{b l}\|_{L^1
              (0, \, + \infty)} \Big)  \\
      &       + L_2 \Big( |t -s| + |\sqrt{t} - \sqrt{s} | \Big).
      \end{split}
      \end{equation}
      \end{teo}

The second theorem concerns the limit as $\ee \to 0^+$. Since we
have a uniform bound on the total variation, by Helly's theorem
there is a subsequence of $u^\ee$ converging in $L^1$ to a limit
function $u(t)$ on a countable dense set of times $t_n$. By the
stability estimate \eqref{stability_wrt_l1_tpar}, the convergence
is on the whole $\mathbb{R}^+$.

However, different subsequences could a priori converge to
different limits: we will actually prove that the limit is unique
and that moreover the semigroup property holds.

      \begin{teo}
      \label{T:2}
      As ${\varepsilon \to 0^+}$, the sequence $u^{\varepsilon}(t)$ of solutions
      of \eqref{vvapproximation} converges to
      a unique function $u(t)$ for all $t \geq 0$: we denote such a limit by
      \begin{equation*}
            u(t) = \semi.
      \end{equation*}
      This convergence defines a unique semigroup
      \begin{equation}
      \label{eq_thm_semigroup}
      \begin{split}
      &      S: [0, \, + \infty] \times \mathcal{U}_{\, 0}
            \times \mathcal{U}_{\, b} \times \mathcal{U}_{\, b}
            \to
            \mathcal{D}_{\, 0}
            \times \mathcal{U}_{\, b} \times \mathcal{U}_{\, b} \\
     &      \quad \; \; (t, \, u_0, \, u_{b \, 0}, \, u_{b \, l})
            \qquad \quad \; \; \mapsto
           \bigg( p_t [u_0, \, u_{b \, 0}, \, u_{b \, l}], \,
           u_{b \, 0}( \, \cdot \, + t), \,
           u_{b \, l}(\, \cdot \, + t)\bigg) \\
      \end{split}
      \end{equation}
      which satisfies the
      following stability estimates in $L^1(0,l)$:
      \begin{equation}
      \label{stability_wrt_l1_t}
      \begin{split}
             \Bigl\| \semi  - p_s [ \bar{v}_{0}, \, \bar{v}_{b \,  0}, \,
             \bar{v}_{b \, l} ] \Bigr\|_{L^1} \leq
      &      L_1 \bigg(  \|\bar{v}_{ 0} - \bar{u}_{0}\|_{L^1(0, \, l)}+
            \|\bar{v}_{b 0} - \bar{u}_{b 0}\|_{L^1(0, \, + \infty)}
            \\
      &     +  \|\bar{v}_{b l} - \bar{u}_{b l}\|_{L^1(0, \, + \infty)} \bigg)+
             L_2 |t -s|,        \\
      \end{split}
      \end{equation}
      for some constant $L_1$, $L_2$ depending only on $A$ and on $\delta_1$.
\end{teo}
\begin{rem}
\label{rem_regularity} By the stability estimate
\eqref{stability_wrt_l1_t} the semigroup $S$ defined by
\eqref{eq_thm_semigroup} can be extended to initial and boundary
data that satisfy much weaker regularity assumptions, i.e.
$\bar{u}_0 \in \domain$ and $\bar{u}_{b \, 0}, \; \bar{u}_{b \, l}
\in \domainb$. Indeed, let $\{ \rho_k \}$ be a sequence of
regularizing kernels and let $\bar{u}_0 \in \domain$. Then $\rho_k
\ast \bar{u}_0$, $\rho_k \ast \bar{u}_{b \, 0}$ and $\rho_k \ast
\bar{u}_{b \, l}$ are initial and boundary data that satisfy the
hypothesis \eqref{E:bounder}: they are smooth and
\begin{equation*}
\begin{split}
&    \| d (\rho_k \ast \bar{u}_0) /   d x \|_{L^1} \leq
    \bv \big\{ \bar{u}_0 \big\} \leq
    \delta_1
     \qquad
     \| d (\rho_k \ast \bar{u}_{b \, 0}) /   d x \|_{L^1}
     \leq \delta_1
     \qquad
     \| d (\rho_k \ast \bar{u}_{b \, l}) /
     d x \|_{L^1} \leq \delta_1 \phantom{\Big( } \\
&    \Big\| d^j (\rho_k \ast \bar{u}_0) /   d x^j \big\|_{L^1} =
     \Big\| d \Big( (d^{j-1} \rho_k / d x^{j-1})
     \ast \bar{u}_0 \Big) / dx \Big\|_{L^1}
     \leq M(k, \, j) \delta_1 \quad \qquad  \qquad j=1, \dots n \\
&    \Big\| d^j \Big(\rho_k \ast \bar{u}_{b \, 0}\Big) /   d x^j
     \Big\|_{L^1}
     \leq M(k, \, j) \delta_1 \qquad
      \Big\| d^j \Big(\rho_k \ast \bar{u}_0) \Big) / d x^j \Big\|_{L^1}
     \leq M(k, \, j) \delta_1 \qquad j=1, \dots n.
\end{split}
\end{equation*}
The last estimates ensures that, for any fixed $k$, the $L^1$ norm
of the derivatives is finite: the bound is not uniform with
respect to $k$ but, since the constant $L_1$ in
\eqref{stability_wrt_l1_t} does not depend on the bound $M$ in
\eqref{E:bounder}, it is enough to prove the extendibility of the
semigroup to the whole domain $\domain$. Indeed, let $u^{\ee}_k$
the sequence of solutions to the systems
\begin{equation*}
\left\{
\begin{array}{lllll}
       \Big( u^{\ee}_k \Big)_t + A(u^{\ee}_k ) \Big( u^{\ee}_k
       \Big)_x = \ee \Big( u^{\ee}_k \Big)_{xx} \\
       \\
       u^{\ee}_k (0, \, x) = \rho_k \ast \bar{u}_0 \phantom{\Big(} \\
       \\
       u^{\ee}_k(t, \, 0) = \rho_k \ast \bar{u}_{b \, 0}
       \quad
       u^{\ee}_k (t, \, l) = \rho_k \ast \bar{u}_{b \, l}
       \phantom{\Big(}
\end{array}
\right.
\end{equation*}
Theorem \ref{T:2} ensures that, for any $k \in \mathbb{N}$ and for
any $t \ge 0$, the sequence $u^{\ee}_k (t)$ converges as $\ee \to
0^+$ to some limit function we will call $u_k (t)$. Then $u_k (t)$
is a Cauchy sequence since by \eqref{stability_wrt_l1_t}
\begin{equation*}
      \| u_k (t) - u_h (t) \|_{L^1(0, \, l)} \leq
      L_1 \Big( \| ( \rho_k - \rho_h ) \ast \bar{u}_0 \|_{L^1(0, \, l)} +
      \| ( \rho_k - \rho_h ) \ast \bar{u}_{b \, 0} \|_{L^1(0, \, + \infty)} +
      \| ( \rho_k - \rho_h ) \ast \bar{u}_{b \, l} \|_{L^1(0, \, + \infty )}
      \Big).
\end{equation*}
The same estimate \eqref{stability_wrt_l1_t} ensures that the
limit $\lim_{k \to + \infty} u_k (t)$ does not depend on the
choice of the sequence $\rho_k$ and therefore the extension
\begin{equation*}
       \semi = \lim_{k \to + \infty} u^k(t)
\end{equation*}
is well defined.

For simplicity, in the following we won't prove that, if
$(\bar{u}_0, \, \bar{u}_{b \, 0}, \, \bar{u}_{b \, l})$ belongs to
$\mathcal{D}_0 \times \mathcal{D}_b \times \mathcal{D}_b$ but not
to $\domainu \times \domainub \times \domainub$, then the solution
of the system \eqref{vvapproximation} converges as $\ee_n \to 0^+$
to $\semi$. However, we will exploit the extendibility property
described before, in particular in Section \ref{par_riemann} we
will consider the {\it vanishing viscosity solution} of the
Riemann and of the boundary Riemann problem, actually meaning the
{\it extension of the semigroup of the vanishing viscosity
solution} to piecewise constant initial and boundary data.
\end{rem}
The function $u(t) = \semi$ is the {\it vanishing viscosity
solution} to
\begin{equation}
\label{E:hyp6}
       u_t + A(u) u_x = 0.
\end{equation}
Note that it is not a weak solution, unless the system is
conservative, but one can prove that it is a {\it viscosity
solution}, in the sense of \cite{AmaCol}. In particular, we obtain
that, for a.e. $t$, the limits
\begin{equation}
\label{E:limB}
       \lim_{x \to 0^+} u(t, \, x) = u(t, \, 0^+),
       \quad \lim_{x \to l^-} u(t, \, x) = u(t, \, l^-)
\end{equation}
and the boundary data $\bar{u}_{b \, 0}(t)$, $\bar{u}_{b \, l}(t)$
can be connected by boundary profiles, i.e. there exists a
solution of the boundary value problem
\begin{equation*}
\left\{
\begin{array}{ll}
       A(v)v_x = v_{xx}, \quad x \in \, ]0, \, + \infty[ \\
       v(0)= \bar{u}_{b \, 0}(t),
       \quad
       \lim_{x \to + \infty } v(x)= u(t, \, 0^+)
\end{array}
\right. \quad \mathrm{and} \quad
\left\{
\begin{array}{ll}
       A(v)v_x = v_{xx}, \quad  x \in \, ]- \infty, \, 0[ \\
       v(0)= \bar{u}_{b \, l}(t),
       \quad
       \lim_{x \to - \infty } v(x)= u(t, \, l^-)
\end{array}
\right.
\end{equation*}
respectively. This means that the boundary datum $\bar{u}_{b  \,
0}$ lies on the stable manifold of $u(t, \, 0+)$, and the boundary
datum $\bar{u}_{b \, l}$ lies on the unstable manifold of $u(t,\,
l^-)$.

The paper is organized as follows.

First of all we make a change of variables in
\eqref{vvapproximation}: let $u(x, \, t) : = u^{\varepsilon}(x /
\varepsilon, \,  t / \varepsilon$). Then \eqref{vvapproximation}
is equivalent to the system
\begin{equation}
\label{rescaled}
      \left\{
      \begin{array}{lllll}
            u_t + A (u) u_x =
            u_{xx},
            \quad
            x \in \, ]0, \, L[, \;
            t \in \, ]0, + \infty [ \\
            \\
            u(0, x) = {u}_0 (x), \\
            \\
            u(t, 0) = {u}_{b \, 0}(t),  \qquad
            u(t, L) = {u}_{b L}(t) \\
     \end{array}
     \right.
\end{equation}
where $L = l / \varepsilon, \; u_{b \, 0}(t)=\bar{u}_{b \, 0}(t /
\varepsilon) , \; u_{b \, L}(t)= \bar{u}_{b \, l}(t / \varepsilon)
, \; {u}_{ 0}(x) = \bar{u}_ 0(x / \varepsilon)$. One can easily
check that
\begin{equation*}
\begin{split}
&     \bv \big\{ \bar{u}_{b \, 0} \big\} = \bv \big\{ {u}_{b \,
      0}\big\}
      \leq \delta_1
      \qquad
      \bv \big\{\bar{u}_{b \, l}\big\} = \bv({u}_{b \,L})
       \leq \delta_1 \\
&     \qquad \qquad \qquad
      \bv \big\{\bar{u}_{ 0}\big\} = \bv \big\{{u}_{0} \big\} \leq \delta_1.
\end{split}
\end{equation*}
Moreover, the derivatives of the boundary and initial data
satisfy
\begin{equation}
\label{estimate_upper}
      \|d^j {u}_0 / dx^j\|_{L^1}, \; \|d^j {u}_{b \, 0} /
      dt^j\|_{L^1}, \; \|d^j {u}_{b \, L} / dt^j\|_{L^1}
      \leq M \varepsilon^{j-1}
      < \delta_1 \quad j= 2, \, \dots n
\end{equation}
for $\ee$ small enough.

The crucial tool in the proof of the convergence of the solution
of \eqref{rescaled} as the scaling parameter $\varepsilon \to 0^+$
is Helly's theorem. One needs therefore to prove a uniform bound
on the total variation, independent on the length of the interval
$L$ and on the $L^1$ norm of the boundary and initial data.

In Section \ref{parabolic_estimates} we prove a priori bounds on
the solution of \eqref{rescaled} that ensure the local existence
and smoothness of solution. Moreover, we will show that, as long
as the total variation of the solution remains small, the $L^1$
norm of $u_{xx}$ is small too and the solution itself can be
prolonged in time. The proof is based on the following
observation: \eqref{rescaled} can be seen as a perturbed heat
equation and therefore one is led to introduce suitable
convolution kernels. Since the technique used in this section does
not depend on the dimension of the solution $u$, we perform the
computations for the $n \times n$ system.

In Section \ref{gradient_decomposition} we introduce the crucial
tool in the proof of the $BV$ estimates: a suitable decomposition
of the gradient of the solution. In the boundary free case
\cite{BiaBrevv}, the gradient $u_x$ is decomposed along a suitable
set of unit vectors $\tilde r_i$, $i=1,\dots,n$, which correspond
to the tangent vectors of the travelling wave profiles of
\[
u_t + A(u) u_x = u_{xx}.
\]
In the single boundary case \cite{AnBia}, instead, the gradient
$u_x$ is decomposed along $n$ travelling wave profiles (the same
as in the boundary free case) and along a boundary profile, i.e. a
solution to the stationary system
\[
   u_{xx} = A(u) u_x.
\]
Such a boundary profile lays on a manifold whose dimension is
related to the number of negative eigenvalues of $A(u)$, i.e. to
the number of characteristic fields that leave the domain $x >0$.

In our case, the basic idea is to split the part of the gradient
due to the presence of the initial datum from the part due to the
boundary data: the first part will be decomposed along the same
tangent vectors $\tilde r_1$, $\tilde r_2$ to travelling wave
profiles introduced in \cite{BiaBrevv}. Moreover, following the
same ideas as in \cite{AnBia}, in order to decompose the part of
the gradient due to the boundary data we use double boundary
profiles, i.e. suitable solutions of the stationary system
\begin{equation}
\label{eq_system_bp}
\left\{
\begin{array}{ll}
     u_x = p, \\
     p_x = A(u)p.
\end{array}
\right.
\end{equation}
In the linear case the two components of the system
\eqref{eq_system_bp} are decoupled and one can show that there is
a solution of the boundary value problem
\begin{equation}
\label{eq_system_bvp}
\left\{
\begin{array}{lll}
     u_x = p, \\
     p_x = A(u)p, \\
     u(0) = U_{b \, 0}, \quad
     u(L) = U_{b \, L}
\end{array}
\right.
\end{equation}
with total variation uniformly bounded with respect to $L$.

In the general case, the idea is to emulate the linear case, using
the center-stable manifold theorem coupled with a contraction
mapping argument: one finds that, provided the difference ${|U_{b
\, 0}- U_{b \, L}|}$ is small, there is a solution of
\eqref{eq_system_bvp} with uniformly bounded total variation. Such
a solution can be seen as the sum of two components, one
exponentially decreasing as $x \to + \infty$, the other as $x \to
- \infty$: we will denote by $\hat{r}_1$ and $\hat{r}_2$ the
tangent vectors to the first and the second part respectively. It
is important to underline, however, that in the non linear case
the two components are coupled: indeed, one finds that
$\hat{\lambda}_1$, the speed of exponential decay of the first
component, depends also on the second component, and viceversa
$\hat{\lambda}_2$ depends on the first component. The introduction
of the generalized eigenvalues $\hat{\lambda}_1$ and
$\hat{\lambda}_2$ allows the equations satisfied by the components
of the decomposition to be exactly in conservation form.

The decomposition of the gradient along travelling waves profiles
and double boundary layers takes the form
\begin{equation}
\label{eq_the_decomposition}
      u_x = v_1 \tilde{r}_1 +
          v_2 \tilde{r}_2 +
          p_1 \hat{r}_1 +
          p_2 \hat{r}_2.
\end{equation}
In Section \ref{par_double_bp} we will show that, because of the
triangular structure of the matrix $A$, the vector $\tilde{r}_2$
and $\hat{r}_2$ can be chosen to be identically equal to $r_2 =
(0, \, 1)$ and $\hat{\lambda}_1$ is identically equal to
$\lambda_1$.

Note that \eqref{eq_the_decomposition} is a system of 2 equations
in 4 unknowns: this allows some freedom in choosing in the most
suitable way the boundary and initial conditions. The precise
expression of all the boundary and Cauchy data we will impose on
$v_1, \; v_2, \; p_1$ and $p_2$ can be found in Section
\ref{boundary_conditions}, in the following however we will sketch
the crucial ideas involved in the choice of those conditions.

We need a preliminary observation: besides that in the choice of
the boundary conditions, some freedom is also allowed in the
attribution of the source terms. Indeed, if one inserts
\eqref{eq_the_decomposition} in the system
\begin{equation*}
       u_t + A(u) u_x - u_{xx} =0
\end{equation*}
obtains the equations
\begin{equation*}
\begin{split}
&      v_{1 t } +(\lambda_1 v_1)_x - v_{1 xx}+
       p_{1 t } +(\lambda_1 p_1)_x - p_{1 xx} = 0 \\
&      v_{2 t } +(\lambda_2 v_2)_x - v_{2 xx}+
       p_{2 t } +(\hat{\lambda}_2 p_2)_x - p_{2 xx} = \tilde{s}_1(t, \, x) \\
\end{split}
\end{equation*}
for some function $\tilde{s}_1$ whose exact expression can be
found in the Appendix \ref{explicit_source_t} and is not important
at the moment: however, it is crucial to observe that it is
identically zero when the solution is exactly a travelling wave or
a double boundary profile. Moreover, in general such a source term
is spread on the whole interval $]0, \, L[$: since $p_2$, the part
of the double boundary layer exponentially decaying as $x \to -
\infty$, should be affected only by the datum in $x=L$, it seems
reasonable to impose
\begin{equation}
\label{cons_laws}
\begin{split}
&      v_{1 t } +(\lambda_1 v_1)_x - v_{1 xx}=0
       \qquad \qquad \; \; \;
       p_{1 t } +(\lambda_1 p_1)_x - p_{1 xx} = 0 \\
&      v_{2 t } +(\lambda_2 v_2)_x - v_{2 xx} = \tilde{s}_1(t,
       \,x)
       \qquad
      p_{2 t } +(\hat{\lambda}_2 p_2)_x - p_{2 xx} =
       0. \\
\end{split}
\end{equation}

As regards the boundary and initial data we impose on the
components $p_1$, $p_2$, $v_1$ and $v_2$, we first observe that,
since $p_1$ and $p_2$ are the components of $u_x$ along double
boundary profiles, we don't want them to be influenced by the
initial datum. Hence we impose
\begin{equation*}
      p_1 (0, \, x) \equiv 0
      \qquad
      p_2 (0, \, x) \equiv 0.
\end{equation*}
Moreover, $p_1$ is the exponential decreasing component of the
boundary profile and hence it should not be affected too much by
the datum on the boundary $x=L$: more precisely, since the goal is
to establish a uniform bound on the $L^1$ norm of $p_1$, it seems
reasonable to look for some boundary condition that minimizes the
increment of $\| p_1\|_{L^1(0, \, L)}$ due to the datum on the
boundary $x=L$. An integration by parts ensures that
\begin{equation*}
      \frac{d }{ dt} \int_0^L |p_1 (t, \, x)| \leq
      |p_{1 x} - \lambda_1 p_1 |(t, \, L) +
       |p_{1 x} - \lambda_1 p_1 |(t, \, 0)
\end{equation*}
and therefore we will impose
\begin{equation*}
       |p_{1 x} - \lambda_1 p_1 |(t, \, L) \equiv 0
\end{equation*}
and, by analogous considerations,
\begin{equation*}
       |p_{2 x} - \hat{\lambda}_2 p_2 |(t, \, 0) \equiv 0.
\end{equation*}
On the other hand, $v_1$ and $v_2$ are the components of $u_x$
along travelling profiles and therefore we don't want them to be
strongly influenced by the presence of the boundary data. We
observe that, in the hyperbolic limit
\begin{equation*}
      u_t + A(u) u_x =0,
\end{equation*}
the waves of the first family go out from the domain through the
boundary $x=0$: we would like to emulate such a behavior in the
parabolic approximation. More precisely, since the aim is to show
a uniform bound on the $L^1$ norm of $v_1$, we look for some
boundary condition that ensures that the derivative of the wave in
the parabolic approximation crosses the boundary, as in the
hyperbolic limit. To make the situation clearer, it is useful to
consider the simple examples that follow: consider the linear
scalar equation
\begin{equation}
\label{eq_linear}
      z_t + \lambda_1^{\ast}z_x - z_{xx} =0
\end{equation}
with some Dirichlet condition imposed on the boundaries $x=0$ and
$x=L$, for example
\begin{equation}
\label{eq_dirichlet}
      z(t, \, 0) \equiv 0, \qquad
      z(t, \, L ) \equiv 1.
\end{equation}
Moreover, let $z^{D}(t, \, x)$ be a solution of \eqref{eq_linear}
and \eqref{eq_dirichlet}: the initial condition is not important
at the moment, but suppose for simplicity that $\bv \big\{
z^{D}(0, \, x) \big\} =1$. For sure $\bv \big\{ z^{D}(t)\big\} \ge
1$ and hence the derivative of $z^{D}$ cannot cross the boundary
$x=0$, or at least the loss of total variation that occurs at
$x=0$ has to be compensated by an increase at $x=L$.

On the other hand, let $z^{N}(t, \, x)$ be a solution of
\eqref{eq_linear} that satisfies a homogeneous Neumann condition
at $x=0$, for example
\begin{equation*}
      z_x^{N}(t, \, 0) \equiv 0,
      \qquad z^{N}(t, \, L) \equiv 1,
\end{equation*}
then an integration by parts ensures that
\begin{equation*}
       \frac{d}{dt} \int_0^L |z^N_x(t, \, x)| dx \leq
       - |z^N_{xx}(t, \, 0)|,
\end{equation*}
and hence the total variation of $z^{N}$ is flowing out from the
domain through the boundary $x=0$.

Hence we are are led by the previous considerations to impose on
the boundary $x=0$ a homogeneous Dirichlet condition on the
function $v_1$, which corresponds to the derivative of a
travelling wave of the first family:
\begin{equation*}
      v_1 (t, \, 0) \equiv 0.
\end{equation*}
The considerations that motivate the choice
\begin{equation*}
       v_2 (t, \, L) \equiv 0
\end{equation*}
are completely analogous.

In Section \ref{BV_estimates} we exploit the decomposition
\eqref{eq_the_decomposition} to prove that the total variation is
uniformly bounded by $\unpo \delta_1$. As we will see, the crucial
point is to prove that, if $\bv \big\{ u_x(\sigma) \big\} \leq
\unpo \delta_1$ for all $\sigma \leq t$, then it holds an estimate
of order two on the integrals of the source term:
\begin{equation}
\label{eq_bound_order2}
      \int_0^{t} \int_0^L |\tilde{s}_1(\sigma, \, x)| dx d \sigma \leq \unpo
      \delta_1^2.
\end{equation}
To show \eqref{eq_bound_order2} we will basically deal with each
of the term that appear in the expression of $\tilde{s}_1$
separately. Some of the estimates are based on the same techniques
 described in \cite{BiaBrevv}: in particular we will use the interaction,
area and length functional introduced in the boundary free case.
Some estimates, on the other hand, require quite long computations
and can be found in the appendix.

In Section \ref{stability_estimates} we will prove the stability
of the vanishing viscosity approximation with respect to $L^1$
perturbations. More precisely, let $u_0, \; u_{b \, 0}, \; u_{b \,
L}$ and $v_0, \; v_{b \, 0}, \; v_{b \, L}$ be the initial and
boundary data of two solutions $u$ and $v$ of problem
\eqref{rescaled}: we will show that there exists a constant $L_1$
such that
\begin{equation*}
              \| u( t) - v(t) \|_{L^1(0,L)}
              \leq L_1 \Big( \|u_0 - v_0 \|_{L^1(0, \, L)} +
              \| u_{b 0} - v_{b 0}\|_{L^1(0, \, t)} +
              \| u_{b L} - v_{b L}\|_{L^1(0, \, t)} \Big).
\end{equation*}
Moreover, one has also stability with respect to time: if $u$ is a
solution to \eqref{rescaled} then
\begin{equation*}
       \|u(t) - u(s)\|_{L^1} \leq
       L_2 \big( |t - s| + | \sqrt{t} - \sqrt{s}| \big)
\end{equation*}
for a suitable constant $L_2$. We will see that the constants
$L_1$ and $L_2$ depend uniquely on the matrix $A$ and on the bound
$\delta_1$ on the total variation of the initial and boundary
data. We will actually give just a sketch of the proof of the
stability, since we will show that one can employ the same tools
used to prove the $BV$ estimates and repeat with minor changes the
computations of Section \ref{BV_estimates}.

One can then get back to the solution $\ue$ of the original
problem \eqref{vvapproximation} and obtain that for all $\ee >0$
it satisfies
\begin{equation}
\label{estimate_stability}
\begin{split}
&      \bv \big\{ u^{\varepsilon}(t) \big\} \leq \unpo \delta_1
       \quad \forall \, t>0
       \qquad
       \| \ue(t) - u^* \|_{\infty} \leq
       \unpo \delta_1 \quad \forall \, t > 0 \\
&       \| u^{\ee }( t) - v^{\ee}(t) \|_{L^1(0, \, L) }
        \leq L_1 \big( \|\bar{u}_0 -
        \bar{v}_0 \|_{L^1(0, \, L)} + \| \bar{u}_{b 0} -
        \bar{v}_{b 0}\|_{L^1(0, \, t)} + \| \bar{u}_{b 0} -
        \bar{v}_{b L}\|_{L^1(0, \, t)} \big) \\
&       \|\ue(t) - \ue(s)\|_{L^1} \leq
       {L}_2 \big( |t - s| + \sqrt{\ee} \, | \sqrt{t} - \sqrt{s}|
       \big). \\
\end{split}
\end{equation}
In the last estimate, $\bar{u}_0, \; \bar{u}_{b 0} \; \bar{u}_{b
L}$ and $\bar{v}_0, \; \bar{v}_{b 0} \; \bar{v}_{b L}$ are the
initial and boundary data for two solutions $u^{\ee }$ and
$v^{\ee}$ of \eqref{vvapproximation}.
%

The uniform bound on the total variation of the solutions
$u^{\ee}$ of \eqref{vvapproximation} ensures that for any
$(\bar{u}_0, \, \bar{u}_{b \, 0}, \, \bar{u}_{b \, l }) \in \,
\mathcal{U}_{\, 0} \times \mathcal{U}_{\, b} \times \domainub$,
for any $t> 0$ and $\ee_n \to 0^+$ there is a subsequence
$\ee_{n_k}$ such that $u^{\ee_{n_k}}(t)$ converges in $L^1(0, \,
l)$ to some limit function we will denote by $\semi$. Letting $\ee
\to 0^+$ in \eqref{estimate_stability} one finds that the limit
satisfies the stability estimate
\begin{equation}
\label{stability_wrt_l1_t_II}
      \begin{split}
             \Bigl\| \semi  - p_s [ \bar{v}_{0}, \, \bar{v}_{b \,  0}, \,
             \bar{v}_{b \, l} ] \Bigr\|_{L^1} \leq
      &      L_1 \bigg(  \|\bar{v}_{ 0} - \bar{u}_{0}\|_{L^1(0, \, l)}+
            \|\bar{v}_{b 0} - \bar{u}_{b 0}\|_{L^1(0, \, + \infty)}\\
      &     +  \|\bar{v}_{b l} - \bar{u}_{b l}\|_{L^1(0, \, + \infty)} \bigg)+
             L_2 |t -s|.       \\
      \end{split}
\end{equation}
By a standard diagonalization procedure one can show that there is
a subsequence that converges for any rational time $t$ and for any
$(\bar{u}_0, \, \bar{u}_{b \, 0}, \, \bar{u}_{b \, l})$ in a
countable dense set of $\mathcal{U}_0 \times \mathcal{U}_{\, b}
\times \mathcal{U}_{\, b}$; the density is here intended in the
$L^1$ norm. Then by the estimate \eqref{stability_wrt_l1_t_II}
$\semi$ must be defined on close sets of times and boundary and
initial data. Hence $\semi$ is defined for any $t \ge 0$ and for
all $ (\bar{u}_0, \, \bar{u}_{b \, 0}, \, \bar{u}_{b \, l}) \in
\mathcal{U}_{\, 0} \times \mathcal{U}_{\, b} \times
\mathcal{U}_{\, b}$.

One can actually check that the operator
\begin{equation*}
      \begin{split}
      &      S: [0, \, + \infty] \times \mathcal{U}_{\, 0}
            \times \mathcal{U}_{\, b} \times \mathcal{U}_{\, b}
            \to
            \mathcal{D}_{\, 0}
            \times \mathcal{U}_{\, b} \times \mathcal{U}_{\, b} \\
     &      \quad \; \; (t, \, \bar{u}_0, \, \bar{u}_{b \, 0}, \,
            \bar{u}_{b \, l})
            \qquad \quad \; \; \mapsto
           \bigg( \semi, \,
           \bar{u}_{b \, 0}( \, \cdot \, + t), \,
           \bar{u}_{b \, l}(\, \cdot \, + t) \bigg) \\
      \end{split}
      \end{equation*}
satisfies the semigroup property

To complete the proof of Theorem \ref{T:2} one is therefore left
to show the uniqueness of the semigroup of vanishing viscosity
solutions: indeed, different sequences $u^{\ee_n}(t)$,
$u^{\nu_n}(t)$ could a priori converge to different limits.

The proof of the uniqueness of the vanishing viscosity limit can
be found in Section \ref{par_viscosity_solutions} and, following
the same ideas as in \cite{BiaBrevv}, the crucial step will be to
show that the semigroup defined via vanishing viscosity
approximation is actually a {\it viscosity solution} in the sense
of \cite{AmaCol}.

We refer to Section \ref{par_viscosity_solutions} for the precise
statement, here however we underline that the definition of
viscosity solution is based on local estimates that ensure,
roughly speaking, a "good behavior" in comparison with the
solutions of a suitable Riemann problem and of a suitable linear
problem.

 The
notion of viscosity solution was first described in the
conservative boundary free case in \cite{Bre:Gli} and was strictly
connected to the definition of Standard Riemann Semigroup (SRS)
that was introduced in the same paper. For completeness, we recall
here that a SRS is Lipschitz continuous with respect to the $L^1$
norm and in the case of piecewise constant initial data locally
coincides with the standard Riemann solver defined by Lax in
\cite{Lax}. In \cite{Bre:Gli} it is proved that if a SRS semigroup
exists, then it necessarily coincides with the wave-front tracking
limit and with the viscosity solution. One of the main advantages
one gains introducing the notion of viscosity solution is
therefore the characterization of global behaviors through local
ones.

The definition of SRS semigroup and of viscosity solution was
extended to conservative boundary value problems in \cite{AmaCol}.
Moreover, in the same paper it was proved that, also in the case
of an initial-boundary value problem, if a SRS exist then it
necessarily coincides with the wave-front tracking limit and with
the viscosity solution. Hence the uniqueness of the SRS semigroup
comes from the uniqueness of the wave-front tracking limit, proved
in \cite{DonMar}.

From the previous works it is clear that a crucial step in the
definition of viscosity solution is the description of the
Riemann solver and of the boundary Riemann solver.  \\
As mentioned before, a solution of the Riemann problem in the
boundary free case was introduced by Lax (\cite{Lax}) for
conservative systems in the case of linearly degenerate or
genuinely non linear fields. Such a definition was then extended
by Liu (\cite{Liu:riemann}) to very general conservative systems.
The characterization of the Riemann solver for non conservative
systems was introduced in \cite{BiaBrevv}, where it was also
proved the effective convergence of the vanishing viscosity
solutions and it was extended in the natural way the notion of SRS
and of viscosity solution.

As concerns boundary Riemann solvers, a solution of the initial
boundary value problem
\begin{equation}
\label{eq_boundary_Riemann}
\left\{
\begin{array}{ll}
       u_t + A(u) u_x =0 \\
       \\
       u(t, \, 0) \equiv \bar{u}_b
       \quad
       u(0, \, x) \equiv \bar{u}_0,
\end{array}
\right.
\end{equation}
was proposed in \cite{DubLeFl} in the case of systems in
conservation form with only linearly degenerate or genuinely non
linear fields: such a boundary Riemann solver is in general
different from the one defined by the vanishing viscosity limit
(some more precise considerations can be found in Remark
\ref{rem_comparison}). On the other side, in \cite{Good,
SabTou:mixte, Ama} and \cite{AmaCol} it was considered a quite
general boundary condition, which turns out to be compatible with
the one defined by the limit of vanishing viscosity
approximations: we refer again to Remark \ref{rem_comparison} for
a more precise statement. We underline, moreover, that a study of
the boundary conditions defined by the limit of the general
parabolic approximation
\begin{equation*}
       \ue_t + f(\ue)_x = \ee \big(  B(\ue) \ue_x \big)_x
\end{equation*}
can be found in \cite{Gisclon-Serre:etudes, Gisclon:etudes,
Rousset:inviscid, Serre-Zum, Rousset:residual, Rousset:char} in
the case of systems in conservation form. Finally, the Riemann
solver for boundary value problems non necessarily in conservation
form was first described in \cite{AnBia}; in this paper it was
also extended in the natural way the notion of SRS and of
viscosity solution.

In Section \ref{par_riemann} we will describe the Riemann solver
and the boundary Riemann solver defined by the vanishing viscosity
limit, which however have an interest in their own. The problem
dealt with is actually a particular case of the one solved in
\cite{AnBia}, where also the characteristic case was considered,
but since the reduction to our case is not completely trivial, we
will describe it explicitly. In particular, we will consider the
vanishing viscosity solution of the boundary Riemann problem
\eqref{eq_boundary_Riemann}. Let $u(0^+)=\lim_{x \to 0^+}u(t, \,
x)$ be the trace of the solution on the axis $x=0$, which does not
depend on time since the solution $u$ is self-similar. We will
show that there exists a solution of the ODE
\begin{equation}
\label{eq_parabolic_one}
       A(U) U_x = U_{xx}
\end{equation}
such that
\begin{equation*}
       U( 0) = \bar{u}_b,
       \quad
       \lim_{x \to + \infty} U(x) = u(0^+).
\end{equation*}
In other words, the boundary datum $\bar{u}_b$ does not
necessarily coincide with the trace $u(0^+)$, but it certainly
lays on the stable manifold of $u(0^+)$ with respect to the ODE
\eqref{eq_parabolic_one}.

 \vspace{0.5cm}

\begin{rem} The fact that the bounds on the total variation are uniform
with respect to the length $L$ of the interval implies that, for
any fixed $\ee>0$, one can let $L \to + \infty$ in
\eqref{rescaled}. Hence, coming back to the original system
\eqref{vvapproximation} one finds that also the solutions of
\begin{equation*}
      \left\{
      \begin{array}{lllll}
            u^{\varepsilon}_t + A (u^{\varepsilon}) u^{\varepsilon}_x =
            \varepsilon u^{\varepsilon}_{xx},
            \quad
            x \in \, ]0, \, + \infty[, \quad
            t \in \, ]0, + \infty [ \\
            \\
            u^{\varepsilon}(0, x) = \bar{u}_0 (x)\\
            \\
            u^{\varepsilon}(t, 0) = \bar{u}_{b \, 0}(t) \\
     \end{array}
     \right.
\end{equation*}
have total variation uniformly bounded with respect to $\ee$.

Hence the analysis of the vanishing viscosity approximations of
the initial-one-boundary value problem can be deduced as a limit
case from the study of the two boundaries case.
\end{rem}

\section{Parabolic estimates}
\label{parabolic_estimates} In this section we will find a
representation formula for the solution to \eqref{rescaled}
\begin{equation}
\label{rescaled4}
      \left\{
      \begin{array}{lllll}
            u_t + A (u) u_x =
            u_{xx},
            \quad
            x \in \, ]0, \, L[, \;
            t \in \, ]0, \, + \infty [ \\
            \\
            u(0, \, x) = {u}_0 (x), \\
            \\
            u(t, \, 0) = {u}_{b \, 0}(t), \qquad
            u(t, \, L) = {u}_{b L}(t) \\
     \end{array}
     \right.
\end{equation}
with initial and boundary data satisfying \eqref{E:TVbd},
\eqref{E:noBD} and \eqref{E:bounder}. The aim is to prove that the
solution of \eqref{rescaled4} is regular and that the $L^1$ norm
of the second derivative $\| u_{xx} (t) \|_{L^1(0, L)}$ is
bounded, as soon as the total variation of $u(t)$ remains small.
We will regard \eqref{rescaled4} as a perturbation of the linear
parabolic system with constant coefficients
\begin{equation}
\label{E:vctlin}
     u_t + A^{\ast}  u_x - u_{xx} = 0.
\end{equation}
Here and in the following we will assume $ A^{\ast}= A( u^{\ast})$
and $\lambda^{\ast}_i = \lambda_i (u^{\ast})$.

\subsection{The convolution kernels}
\label{par_conv_kernels }

The fundamental step is to study the equation \eqref{rescaled1}
in the scalar case, because the Green kernel for the general
vector case \eqref{E:vctlin} follows by using the base of eigenvectors of $A^\ast$.
Thanks to the linearity, we split the Green kernel of the equation
\begin{equation}
\label{linear_equation}
      z_t + \lambda_i^{\ast} z_x - z_{xx} =
      0
\end{equation}
into 3 parts:
\begin{enumerate}
\item $ \Delta^{ \lambda_i^{\ast}}(t,x,y)$ is
          the solution of \eqref{linear_equation}
          with zero boundary conditions and
          initial condition
          \begin{equation}
          \label{Delta}
          \nonumber
              \Delta^{\lambda_i^{\ast}} ( 0, \, x, \,  y ) =
              \delta_y
              \,
              \quad
              y \in \, ]0, L[ \, .
          \end{equation}
           This function is given by
%
%
\begin{equation}
\label{Delta_product} 
    \Delta^{ \lambda_i^{\ast}}(t, x, y) =
    \bigg( \sum_{m \,  = \,- \infty}^{m \,  = \, + \infty}
      G (t, x + 2mL -y) -
      G (t, x+ 2mL +y ) \bigg)
    \phi^{\, \lambda_i^{\ast}}(t, \, x, \, y),
\end{equation}
where $G(t, \, x) =  \big( e^{- x^2/ 4t} \big) / 2 \sqrt{\pi t}$
 is the standard heat kernel and
\begin{equation*}
   \phi^{\, \lambda_i^{\ast}}(t, x, y) =
   \exp
      \bigg(
           \frac{
                 \, \, \lambda_i^{\ast}}{2}
           \, (x -y) -
           \frac{
                \, \, (  \lambda_i^{\ast})^2 }{4} \, t \,
      \bigg).
\end{equation*}
\item $J^{\lambda_i^{\ast} \, 0}(t, x)$ is
          the solution of \eqref{linear_equation})
          with zero initial datum and boundary
          conditions
          \begin{equation}
          \label{eq_jzero}
               J^{ \lambda_i^{\ast} 0}( t, \, 0 ) = 1
               \quad
               J^{ \lambda_i^{\ast} 0}( t, \, L) = 0.
          \end{equation}
It follows that
\begin{equation}
\label{J_0}
     J^{ \lambda_i^{\ast} \, 0}(t, x) =
      A \exp \big(
                 \lambda_i^{\ast} x
           \big)
    + B -
    \int_{0}^{L}
          \Delta^{ \lambda_i^{\ast}} (t, x, y)
          \bigg(
                 A \exp \big(
                              \lambda_i^{\ast} y
                        \big)
                + B
         \bigg)
         dy,
\end{equation}
with
\begin{equation*}
       A = - \frac{1}{
                     e^{
                         \lambda_i^{
                                  \ast
                                 } L
                       } - 1
                   } \qquad
     B = \frac{ e^{
                    \lambda_i^{
                              \ast
                            }   L
                  }
              }
       {
                     e^{
                         \lambda^{
                                  \ast
                                 }  L
                       } - 1
                   }.
\end{equation*}
\item $J^{ \lambda_i^{\ast} \, L}(t,x)$
          is the solution of \eqref{linear_equation}
          with zero Cauchy datum and boundary
          conditions
          \begin{equation}
          \label{eq_jelle}
                J^{\lambda_i^{\ast} L}( t, 0 ) = 0
               \quad
               J^{\lambda_i^{\ast} L}( t, L) = 1
          \end{equation}
and it is given by
\begin{equation}
\label{eq_J_L}
       J^{ \lambda_i^{\ast} \, L}(t,x)=
       C \exp \big(
                 \lambda_i^{\ast} x
           \big)
       + D -
       \int_{0}^{L}
          \Delta^{ \lambda{\ast}} (t, x, y)
          \bigg(
                 C \exp \big(
                              \lambda_i^{\ast} y
                        \big)
                + D
         \bigg)
         dy ,
\end{equation}
where
\begin{equation*}
       C =  \frac{1}{
                     e^{
                         \lambda^{
                                  \ast
                                 } L
                       } - 1
                   }
     \qquad
     D  = A = -  \frac{ 1
              }
       {
                     e^{
                         \lambda^{
                                  \ast
                                 }  L
                       } - 1
                   }.
\end{equation*}
\end{enumerate}
Note that all the coefficients $A$, $B$, $C$, $D$ remain bounded
as $ L \to + \infty$. Moreover, one can apply the maximum
principle and, via a comparison with the constant solutions, finds
that $
  0 \leq J^{ \lambda_i^{\ast} \, 0}(t, x), \;
  J^{ \lambda_i^{\ast} \, L}(t, x)
  \leq
  1
$. Hence the integrals
\begin{equation*}
        \int_{0}^{T} J^{
        \lambda_i^{\ast}\, 0}(t, x) v' (t) dt
        \quad 
        \int_{0}^{T} J^{
        \lambda_i^{\ast}\, 0}(t, x) v' (t) dt
\end{equation*}
are well defined for every function $v(t) \in BV(0, \, + \infty)$
and for every $T$.

In the following, we will also need a further convolution kernel
$\tilde{\Delta}^{\lambda_i^{\ast}}(t, \, x, \, y)$ such that
\begin{equation*}
      \tilde{\Delta}^{\lambda_i^{\ast}}_y(t,\, x, \, y) +
      \Delta^{\lambda_i^{\ast}}_x(t, \, x, \, y) = 0,
\end{equation*}
i.e.
\begin{equation}
\label{eq_def_tilde_delta}
      \deltabt(t, \, x, \, y) = \int_y^L
      \Delta_x^{\lambda_i^{\ast}}(t, \, x, \, z) dz.
\end{equation}
To get the previous formula we have arbitrarily imposed $
\deltabt(t, \, x, \, L) = 0$.

Note that $\tilde{\Delta}^{\lambda_i^{\ast}}(t, \, x, \, 0)$ is
the derivative with respect to $x$ of a function $z(t, \, x)$
which satisfies
$$
  z(t, \, x) + J^{\lambda_i^{\ast} \, 0}(t, \, x) +
  J^{\lambda_i^{\ast} \, L}(t, \, x) = 1.
$$
Hence,
\begin{equation}
\label{delta_and_Jx}
      \tilde{\Delta}^{\lambda_i^{\ast}}(t, \, x, \, 0) +
      J_x^{\lambda_i^{\ast} \, 0}(t, \, x) +
      J_x^{\lambda_i^{\ast} \, L}(t, \, x)= 0.
\end{equation}
The following proposition provides some basic estimates on the
convolution kernels we will need later.

\begin{pro}
\label{estimate_kernels_pro}
      The convolution kernel $\Delta^{\lambda_i^{\ast}}$
      satisfies
      \begin{equation}
      \label{estimate_kernels}
               \|\deltab(t, \, y)\|_{L^1} \leq \unpo
               \quad
               \| \deltab_x(t,  \, y)\|_{L^1} \leq
               \unpo / \sqrt{t}  \qquad
               \forall \, t < 1, \; y \in ]0, \, L[ \, .
      \end{equation}
      The following estimates hold for the boundary
      kernels $J^{\lambda_i^\ast \, 0}$, $J^{\lambda_i^\ast \, L}$:
      \begin{equation}
      \label{estimate_kernels_I}
      \begin{array}{ccccc}
             0 \leq J^{ \lambda_i^{\ast} \, 0}(t, x), \;
                      J^{ \lambda_i^{\ast} \, L}(t, x)
                      \leq  1 & \forall \, t \geq 0, \; x \in ]0,L[ \\
             \\
             \|\Jzero_x (t)\|_{L^1},
                      \| \Jelle_x(t)\|_{L^1} \leq \unpo
                      \quad \forall \; \; 0 < t < 1, \\
             \\
             \|\Jzero_{xx} (t)\|_{L^1},
                      \| \Jelle_{xx} (t)\|_{L^1} \leq \unpo/\sqrt{t}
                      \quad \forall \; \; 0 < t <1.
      \end{array}
      \end{equation}
      The auxiliary convolution kernel $\deltabt$ satisfies
      estimates analogous to those of $\deltab:$
      \begin{equation}
      \label{estimate_kernelsII}
             \|\deltabt(t, \,  y)\|_{L^1}
             \leq \unpo \quad
             \|\deltabt_x(t, \, y)\|_{L^1} \leq
            \unpo / \sqrt{t}
            \quad  \forall \; 0 < t < 1, \;  y \in \, ]0, \, L[ \, .
      \end{equation}
\end{pro}
The proof of the proposition can be found in the Appendix \ref{proof_pro_kernels}.

Now we are ready to deal with the vector case. Let $
 \; r^{\ast}_i, \;  l^{\ast}_i   \; i=1, 2
$ be respectively the left and the right eigenvectors of
$A^{\ast}= A(u^{\ast})$. We define the matrix kernels
\begin{equation}
\begin{array}{lll}
           \Delta^{ \ast} : =
           \sum_{ i= 1}^{2}
           \Delta^{ \lambda_i^{\ast} }
           r^{\ast}_i \otimes l^{\ast}_i,
           &
           \tilde{ \Delta}^{ \ast   }: =
           \sum_{ i= 1}^{2}
           \tilde{ \Delta}^{ \lambda_i^{\ast} }
           r^{\ast}_i \otimes l^{\ast}_i,
           \\
           \\
           J^{\,  \ast \, 0} : =
           \sum_{ i= 1}^{2}
           J^{ \lambda_i^{\ast} \, 0 }
           r^{\ast}_i \otimes l^{\ast}_i,
           &
           J^{ \, \ast \, L} : =
           \sum_{ i= 1}^{2}
           J^{ \lambda_i^{\ast} \, L }
           r^{\ast}_i \otimes l^{\ast}_i.
      \end{array}
\end{equation}
By construction these are the matrix kernels for the initial data
corresponding to the cases 1, 2 and 3 considered above (equations
\eqref{Delta}, \eqref{eq_jzero} and \eqref{eq_jelle}
respectively).

\subsection{Parabolic estimates}
\label{sub_paraboli_estimates} The solution of equation
\eqref{rescaled4} can be written as
\begin{equation}
\label{solution}
\begin{split}
         u (t, \, x) =&~
         \int_{0}^{L} \Delta^{\ast}(t, x, y) u_0 (y) dy +
         u_0 (0) J^{\, \ast \,  0}(t, \, x) +
         \int_{0}^{t}
         J^{\, \ast\,  0}( t - s, x ) u'_{b 0} (s) ds
         + u_0 (L) J^{\,   \ast \,  L }(t, \, x)
         \\ &~ +
         \int_{0}^{t}
         J^{\, \ast \,  L}( t - s, x ) u'_{b L} (s) ds
         + \int_{0}^{t} \int_{0}^{L}
         \Delta^{ \ast}(t-s, x, y)
         \big(
              A^{\ast} - A(u)
         \big) u_y (s, y)dy ds, \\
   \end{split}
\end{equation}
and therefore, recalling \eqref{delta_and_Jx} and integrating by
parts,
\begin{equation}
\label{solution_x}
   \begin{split}
          u_x (t, \, x) =&~
         \int_{0}^{L} \tilde{\Delta}^{\ast}
         (t, x, y) u'_{0}(y) dy +
         \int_{0}^{t}
         J_x^{ \, \ast \,  0}( t - s, \, x ) u'_{b 0} (s) ds +
         \int_{0}^{t}
         J_x^{ \, \ast \,  L}( t - s, \, x ) u'_{b L} (s) ds \\
   &~  +
         \int_{0}^{t} \int_{0}^{L}
         \tilde{\Delta}^{\ast}(t-s, \,  x, \, y)
         \bigg(
              \Big(  A^{\ast} - A(u) \Big)
              u_{yy} - DA(u) \Big( u_y \otimes u_y \Big)
         \bigg)     (s, \, y)dy ds \\
   &~ +
         J_x^{\, \ast \, L}(t, \, x)
         \big( u_0(L) - u_0 ( 0) \big) -
         \int_0^t
         \big( J_x^{\, \ast \, 0}+
                J_x^{\, \ast \, L}
         \big) (t-s , \, x)
         (    A^{\ast} - A(u) ) u_x (s, \, 0) ds. \\
   \end{split}
\end{equation}
From the previous expression we immediately have that, as long as
it can be prolonged, the solution is regular. Moreover, the local
existence of a solution of equation \eqref{rescaled4} follows from
the representation formulas \eqref{solution} and
\eqref{solution_x} via the contraction map theorem.

We can now use the representation \eqref{solution_x} to prove the
following proposition.

\begin{pro}
\label{pro_u_xx}
       If
       $
         \| u_x (t )\|_{L^1} \leq \mathcal{O}(1) \delta_1
       $ for all $t  \in [ 0, \, 1  ]$,
       then
        $$
         \| u_{xx} (t )\|_{L^1} \leq \frac{\mathcal{O}(1) \delta_1}{\sqrt{t}}
         \quad \forall \; t \in [  0, \, 1  ].
       $$
\end{pro}
\begin{proof}
From \eqref{solution_x} we get
\begin{equation}
\label{eq_u_xx}
      \begin{split}
             u_{xx}(t, \, x)
            =&~
             \int_{0}^{L}
             \tilde{\Delta}_x^{\ast}(t, \, x, \, y)
              u'_{0 } (y) dy +
               \int_{0}^{t}
               J_{xx}^{\,  \ast \,  0}( t - s, x ) u'_{b 0} (s) ds +
              \int_{0}^{t}
              J_{xx}^{ \, \ast \,  L}( t - s, x ) u'_{b L} (s) ds \\
          &~+
              \int_{0}^{t} \int_{0}^{L}
              \tilde{\Delta}_{x }^{\ast}(t-s, \,  x, \, y)
              \bigg(
              \Big(
                    A^{\ast} - A(u)
              \Big) u_{yy}
              - DA (u)\Big(  u_y \otimes u_y \Big)
              \bigg)(s, y) dy ds \\
       &~ +
              J_{xx}^{\, \ast \, L}(t, \, x)
              \big( u_0(L) -u_0 (0) \big) -
              \int_{0}^{t}
              \big( J_{xx}^{ \, \ast \, 0} +
                    J_{xx}^{ \, \ast \, L}
              \big) ( t-s, \, x)
               \Big(
                   A^{\ast} - A(u)
               \Big) u_x (s, 0) ds. \\
        \end{split}
\end{equation}
The previous representation formula shows that the function $ t
\mapsto \|u_{xx} (t)\|_{L^1}$ is continuous.

We claim that there is a constant $C$ independent from $L$ such
that
$$
   \|u_{xx} (t) \|_{L^1} \leq \frac{C\delta_1}{ \sqrt{t}}
   \qquad \forall \, t <1.
$$
Indeed, for a fixed large constant $C$, define
$$
   \tau = \inf \biggl\{ t : \;
   \|u_{xx} (t) \|_{L^1} \geq \frac{C}{ \sqrt{t}} \delta_1 \biggr\}.
$$
The time $\tau$ is strictly bigger than $0$ if $C$ is sufficiently large,
since by hypothesis $\|u_0 ''\|_{L^1}$ is finite.
Moreover, one has  $\|u_{xx} (\tau)\|_{L^1} = C \delta_1 /
\sqrt{\tau}$ thanks to the continuity of the map $ t \mapsto
\|u_{xx} (t)\|_{L^1}$.

From \eqref{eq_u_xx} it follows that
\begin{equation*}
\begin{split}
     \|u_{xx}(\tau)\|_{L^1} = \frac{C}{ \sqrt{\tau}} \delta_1
&     \leq
      \|\Delta_x^{\ast} (\tau)\|_{L^1} \, \|u'_0\|_{L^1} +
      \unpo \delta_1 \int_0^{\tau} \|u_{yy} (s)\|_{L^1} \,
      \|\Delta_x^{\ast}(\tau - s)\|_{L^1} ds +
      2 \delta_1 \int_0^{\tau}
      \frac{\unpo}{\sqrt{\tau -s}} \,  ds \\
&     \quad +  \frac{ \unpo}{\sqrt{\tau}} \delta_1  +
      2 \delta_1^2 \int_0^t \frac{\unpo \, C}{\sqrt{s ( \tau - s})} \, ds \\
&     \leq \frac{ 2 \unpo \delta_1}{\sqrt{\tau}}+
      2 \unpo \, C \delta_1^2 +
      2 \unpo \, \sqrt{\tau} \delta_1, \\
\end{split}
\end{equation*}
which is a contradiction if $C$ is large enough and $\delta_1$ sufficiently small.
In the previous
estimate we have used the bounds
$$
  \qquad  \qquad  \qquad
  \| u_{b \, 0}'\|_{L^{\infty} } \leq
  \| u_{ b \, 0} ''\|_{ L^1 } \leq
  \delta_1
  \qquad \qquad 
  \int_0^{\tau} \frac{1}{\sqrt{s (\tau - s)}} \, ds = \pi.
$$
\end{proof}
If $ t > 1$ and $\|u_x(s)\|_{L^1} \leq \unpo \delta_1$ for any $s
\in [0, \, t]$ , we can apply the previous proposition to the
interval $ [ t - 1, \, t ]$ and obtain
$$
  \| u_{ xx} ( t) \|_{L^1 } \leq \mathrm \mathrm{O}(1) \delta_1
  \quad t \geq 1.
$$
Since the derivative $u_x$ is regular, this implies in particular
that, if $\|u_{x} (s)\|_{L^1} \leq \mathcal{O}(1) \delta_1$ for
any $s \leq t$, then $ \|u_x (t) \|_{L^{\infty}} \leq
\mathcal{O}(1) \delta_1$ if $ t \geq 1$: in other words, as long
as $u_x$ remains small in the $L^1$ norm, it remains small in the
$L^{\infty}$ norm too.
\section{Gradient decomposition}
\label{gradient_decomposition}
\subsection{Double boundary layers and travelling waves}
\label{par_double_bp} In this section we will introduce a suitable
decomposition of the gradient of the solution to
 \eqref{rescaled},
\begin{equation*}
\left\{
      \begin{array}{lllll}
            u_t + A (u) u_x =
            u_{xx},
            \quad
            x \in \, ]0, \, L[, \quad
            t \in \, ]0, + \infty [ \\
            \\
            u(0, x) = {u}_0 (x), \\
            \\
            u(t, 0) = {u}_{b \, 0}(t), \qquad
            u(t, L) = {u}_{b L}(t). \\
     \end{array}
     \right.
\end{equation*}
We will employ a decomposition in the form
\begin{equation}
\label{eq_decomposition}
    u_x = v_1 \tilde{r}_1 +
          v_2 \tilde{r}_2 +
          p_1 \hat{r}_1 +
          p_2 \hat{r}_2,
\end{equation}
where the first two terms correspond to
derivatives of travelling waves and the
last two correspond to the derivative of a double boundary profile.
More precisely, $p_1$ is the part of the double boundary profile exponentially
decaying as $x \to +\infty$, $p_2$ is the part exponentially
decaying as $x \to -\infty$.

The principal results of this section are the construction of the
vectors $\hat r_1$, $\hat r_2$, the description of a decomposition
of $u_x$ in the form \eqref{eq_decomposition}, the computations of
the equations for the $4$ components $v_1$, $v_2$, $p_1$, $p_2$
and finally the choice of the boundary conditions for the same
components. In the description of the decomposition we will focus
mainly on the construction of the double boundary profiles,
because the construction of the travelling wave profiles follows
the same steps as in \cite{BiaBrevv}.

The construction of the double boundary profile is based on the
following idea: in the linear case, one finds that there is a
solution of the boundary value problem
\begin{equation}
\label{eq_bvp}
\left\{
\begin{array}{lll}
     u_x = p, \\
     p_x = A(u)p, \\
     u(0) = U_{b \, 0}, \quad
     u(L) = U_{b \, L}
\end{array}
\right.
\end{equation}
and such a solution is the sum two components: one exponentially
decaying as $x \to + \infty$, the other as $x \to - \infty$.
Moreover, when the length $L$ is very large the solution has the
behavior illustrated in figure \ref{fig_contractions} (on the
left): it is very steep near the boundary $x=0$ because of the
presence of the exponentially decreasing component, then it is
almost horizontal in a large interval and then it is steep again
near the boundary $x=L$ because of the presence of the exponential
decreasing part.

The idea is to try to simulate such a spatial behavior also in the
non linear case: in this way, when $L$ is large enough the
derivative of the double boundary profile is concentrated near the
boundaries $x=0$ and $x=L$ and therefore there is essentially no
interaction with the travelling wave profiles inside the domain.
This behavior is the same one observes in the hyperbolic limit,
where in $]0,L[$ the solution is generated only by travelling wave
profiles. We will find out that, if $|U_{b \, 0} - U_{b \, L}|$ is
small enough, then there exists indeed a solution of the boundary
value problem \eqref{eq_bvp} with the behavior illustrated in
figure \ref{fig_contractions}.

In this way, we construct the functions $p_1 \hat{r}_1 (u, \,
p_1)$ and $p_2 \hat{r}_2 (u, \, p_2)$: however, since the
decomposition \eqref{eq_decomposition} is a 2-dimensional vector
equation in 4 scalar unknowns, we have some freedom in assigning
the initial and boundary data for $v_1$, $v_2$, $p_1$ and $p_2$.
The detailed description of the boundary conditions can be found
in Section \ref{boundary_conditions}, but the crucial idea is to
impose some conditions that allow the component $p_1$ and $p_2$ to
behave like the derivative of a double boundary layer, and thus to
be independent from the choice of the initial datum and to be
concentrated near the boundary $x=0$ or $x=L$, respectively. On
the other hand, we want to impose some conditions on the
components $v_1$ and $v_2$ that forces them to behave like the
derivative of waves in the hyperbolic limit, thus flowing out from
the domain through the boundary $x=0$ (waves of the first family)
or through the boundary $x=L$ (waves of the second family).

Moreover, we have also some freedom in assigning the source terms,
as it will be clear in Section \ref{par_equations}: again the
basic idea we will follow is that $p_2$, which corresponds to the
component of the double boundary profile exponential decaying as
$x \to + \infty$, should be affected only by the datum in $x=L$.
Since in general the source term are spread on the whole interval
$]0, \, L [$, we will impose that the equation for $p_2$ has no
source term.

\subsubsection{Double boundary profiles}
As a first step, we characterize the solutions of the system
\begin{equation}
\label{boundary_profiles}
      \left\{
      \begin{array}{ll}
           u_x=p \\
           p_x= A(u) p
      \end{array}
      \right.
\end{equation}
that converge with exponential decay to some value $(\bar{u}, \,
0)$ with $\bar{u}$ in a small enough neighborhood of the value
$u^{\ast}$ defined by the relation \eqref{E:noBD}. Since
$(u^{\ast}, 0)$ is an equilibrium point, we can consider the
linearized system, whose center and stable subspaces are given by
\begin{equation*}
      V^c = \{ p = 0 \},
      \qquad
      V^s = \mathrm{span} \langle r_1 (u^{\ast})  \rangle,
      \qquad
      V^u = \mathrm{span} \langle r_2(u^{\ast}) \rangle.
\end{equation*}

\begin{figure}
\caption{the center-stable manifold $\mathcal{M}^{cs}$ and the
center-unstable manifold $\mathcal{M}^{cu}$ with orbits
exponentially decaying to an equilibrium point as $x \to + \infty$
or $x \to - \infty$, respectively} \label{fig_manifolds}
\begin{center}
\psfrag{Z}{$u$} 
\psfrag{M}{$\mathcal{M}^{cu}$} \psfrag{U}{$\mathcal{M}^{cs}$}
\psfrag{V}{$p_1$} \psfrag{W}{$p_2$}
\includegraphics[scale=0.5]{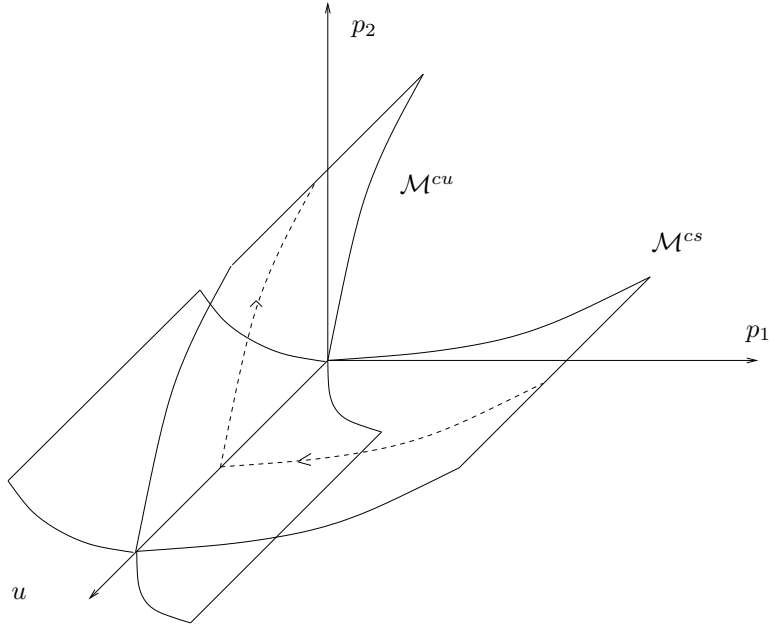}
\end{center}
\end{figure}

Let $(p_1, \, p_2)$ be the coordinates of $p$ with respect to the
base defined by the eigenvalues $r_1(u^{\ast})$ and
$r_2(u^{\ast})$ of $A(u^{\ast})$: thanks to the center-stable
manifold theorem, there exists a regular function
\begin{eqnarray*}
          \phi:
            \{ (u, p_1) : |u- u^{\ast}|
                   , |p_1| \leq \varepsilon
                   \} \subseteq V^c \oplus V^s
            & \to & \mathbb{R} \\
          ( u, p_1) & \mapsto &
                      p_2 = \phi ( u, p_1),
\end{eqnarray*}
which parameterizes the solutions of \eqref{boundary_profiles}
that do not blow up exponentially for $x \to + \infty$. In our
case, one can see that this manifold is made by the orbits which
converge for $x \to + \infty$ to an equilibrium $(\bar{u}, \, 0)$,
with $\bar{u}$ close to $u^{\ast}$ (figure \ref{fig_manifolds}).
In particular this manifold is unique.

The dimension of this manifold is $dim \, V^c + dim \, V^s$, i.e.
3 in our case. Since $ p_1 = 0$ implies $p_2 = \phi(u, p_1) = 0$,
we can set $\phi(u, p_1) = p_1 h( u, p_1)$ and $\mathcal{M}^{cs}$
can be described by the following condition:
\begin{equation*}
      p = p_1 r_1(u^{\ast}) +
          p_1 h(u, p_1) r_2(u^{\ast}) 
=      p_1
      \left(
      \begin{array}{ll}
            \quad 1 \\
            f(u, p_1)
      \end{array}
      \right):  =
      p_1 \hat{r}_1 (u, p_1). 
\end{equation*}
Inserting the previous expression in the system
\eqref{boundary_profiles}, one obtains
\begin{equation*}
      A(u) p_1 \hr= \big( p_1 \hr \big)_x = p_{1 x} \hr +
      (p_1)^2 \bigd \hr \hr + p_1 p_{1 x} \hat{r}_{1 \, p}.
\end{equation*}
Let $\ell_1 = (1, \, 0)$: if we multiply the previous expression
by $\ell_1$ we obtain, since $A$ is triangular,
\begin{equation*}
       \lambda_1 p_1  = p_{1 x } ,
\end{equation*}
and hence
\begin{equation}
\label{E:commu}
       A(u) p_1 \hr = \lambda_1 p_1 \hr + (p_1)^2  \bigd \hr \hr +
      \lambda_1
      p_1^2 \hat{r}_{1 \, p}.
\end{equation}
It follows that
\begin{equation*}
      \hat{r}_1 ( u, \, 0) = r_1(u ) \quad  \forall \; u,
\end{equation*}
and therefore
\begin{equation*}
      | \hat{r}_1 (u, \, p_1) - r_1(u) |\leq
      \mathcal{O}(1)|p_1|.
\end{equation*}

In a similar way one can also define a regular, 3-dimensional
center-unstable manifold $\mathcal{M}^{cu}$ containing all the
orbits that as $x \to - \infty$ converge with exponential decay to
some point $(\bar{u}, \, 0)$ with $\bar{u}$ close to $u^{\ast}.$
The manifold is parameterized by $V^c \oplus V^u$; moreover, since
the matrix $A$ is triangular, one can choose
\begin{equation*}
     \hat{r}_2 \equiv r_2 (u) \equiv
     \left(
     \begin{array}{ll}
       0 \\
       1
     \end{array}
     \right).
\end{equation*}
The manifold $\mathcal{M}^{cu}$ is thus described by the
relation $p= p_2 r_2$.

As a second step, we show that the functions $p_1 \hr$ and $p_2
r_2$ indeed allow us to construct a solution of the two-boundaries
value problem
\begin{equation}
\label{eq_bvpII} \left\{
\begin{array}{lll}
     z_{xx} = A(z) z_x, \\
     z(0) = U_{b \, 0} \quad
     z(L) = U_{b \, L}
\end{array}
\right.
\end{equation}
Decomposing $z_x$ as
\[
z_x  = p_1 \hr (z, \, p_1) + p_2 r_2
\]
and using the relation \eqref{E:commu}, we obtain the system
\begin{equation}
\label{E:systedec}
\left\{
\begin{array}{lll}
     z_x  = p_1 \hr (z, \, p_1) + p_2 r_2, \\
     p_{1x} = \lambda_1(z) p_1, \\
     p_{2x} = \hat \lambda_2(z, \, p_1) p_2
\end{array}
\right.
\end{equation}
where we have defined
\begin{equation}
\label{eq_hat_lambda}
      \hat{\lambda}_2(u, \, p_1) :=
      \lambda_2(u) - p_1
      \langle \hat{\ell}_2, \,
        \mathrm{D} \hat{r}_1 r_2  \rangle ,
\end{equation}
where the vector $\hat{\ell}_2$ satisfies $\langle \hat{\ell}_2,
\, \hr \rangle=0$ and $\langle \hat{\ell}_2, \, r_2 \rangle=1$.
Hence, while in the linear case the two components of the solution
of the system \eqref{eq_bvpII} are decoupled, in the general case
there is a coupling in the equation of $z$, and in the choice of
$\hat \lambda_2$, which is in some sense the effective eigenvalue
for $p_2$. Note that
\begin{equation}
\label{eq_gen_eigen}
      \big|
      \hat \lambda_2(u,p_1) - \lambda_2(u)
      \big|\leq
      \mathcal{O}(1) p_1.
\end{equation}
An application of contraction principle ensures that, if $|U_{b \,
0} - U_{b \, L}| \leq \delta_1$ for a small enough $\delta_1$,
then the above system with boundary data $z(0) = U_{b \, 0}$,
$z(L) = U_{b \, L}$
has a unique solution. 
Moreover, one also finds that $\big| \hat \lambda_2(u, \, p_1) -
\lambda(u)\big| \leq \mathcal{O}(1) \delta_1$.

\begin{figure}
\begin{center}
\caption{the graphic and the orbit of a double boundary layer when
the length $L$ of the interval is large} \label{fig_contractions}
\psfrag{a}{$U_{b \, 0}$} \psfrag{b}{$U_{b \, L}$} \psfrag{L}{$L$}
\vfill
\includegraphics[scale=0.3]{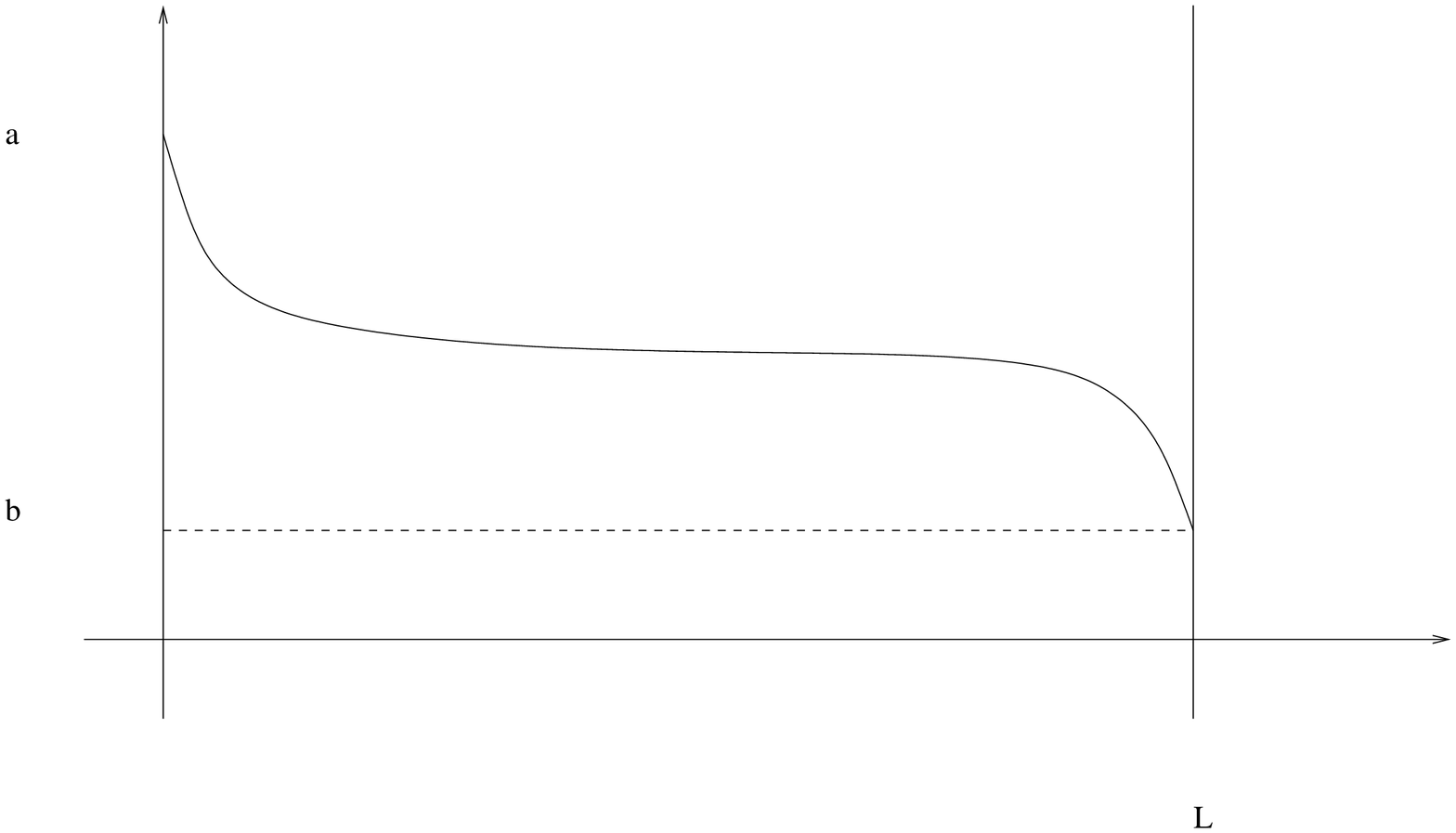}
\hfill \psfrag{Z}{$u$} \psfrag{M}{$\mathcal{M}^{cu}$}
\psfrag{U}{$\mathcal{M}^{cs}$} \psfrag{V}{$p_1$} \psfrag{W}{$p_2$}
\psfrag{a}{$U_{b \, 0}$} \psfrag{b}{$U_{b \, L}$}
\includegraphics[scale=0.3]{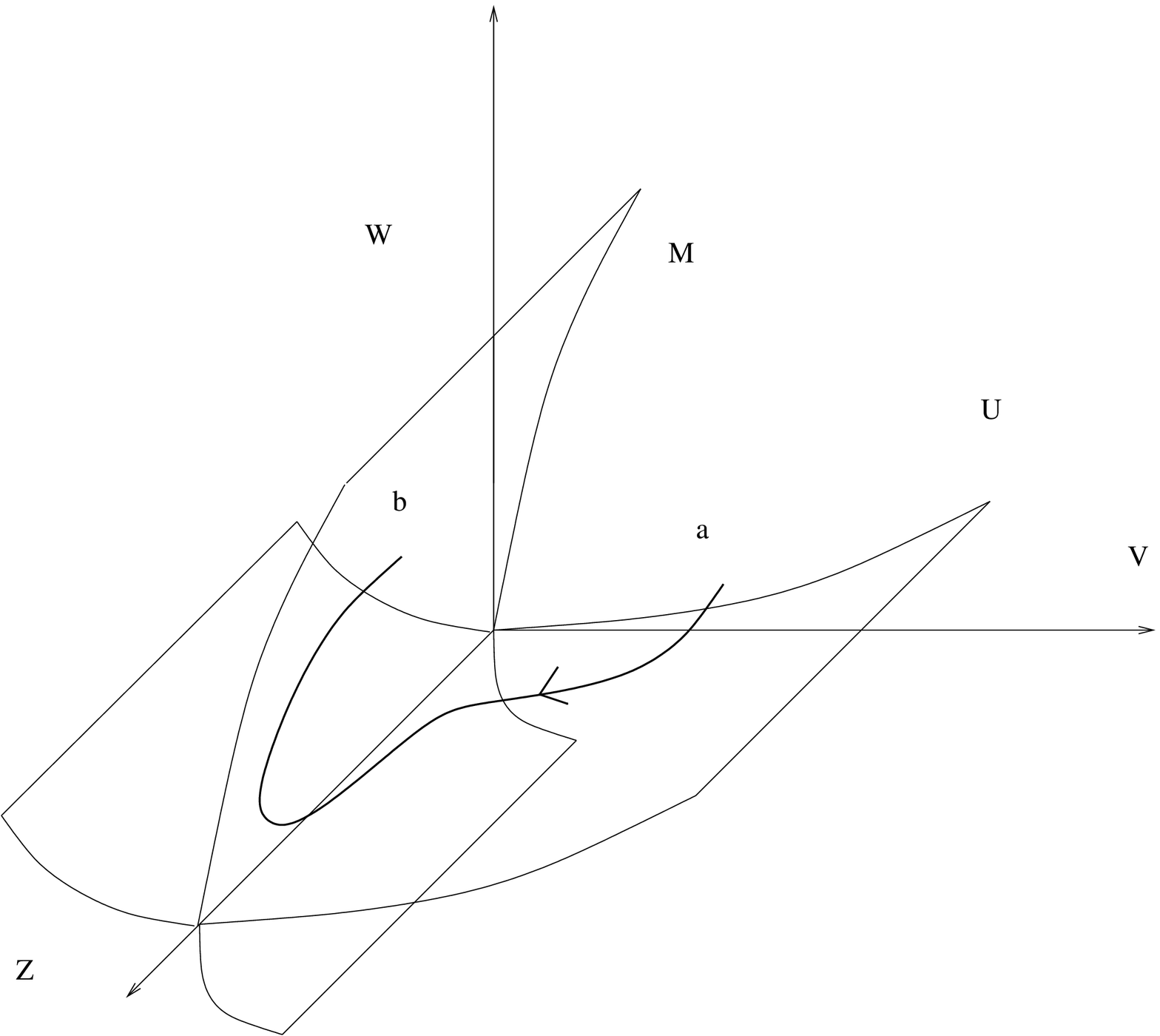}
\end{center}
\end{figure}

Since $\lambda_1 < 0$, $\hat \lambda_2 > 0$ for $\delta_1 \ll 1$,
we obtain that $p_1$ is exponentially decaying, while $p_2$ is
exponentially increasing. We can thus figure the double boundary
profile as follows (figure \ref{fig_contractions}): when the
length $L$ of the interval is very large, the solution will be
steep near zero, because in that region $p_1$ varies exponentially
fast. Then it will be almost horizontal for a long interval and
becomes again very steep in a left neighborhood of $x=L$, because
$p_2$ increases exponentially.

\subsubsection{Travelling waves}
\label{par_travelling_waves}

We refer to \cite{BiaBrevv} for an
exhaustive account of the analysis that allow the definition of
the decomposition along travelling waves: here we will only recall
for
completeness the crucial steps.

Consider the system
\begin{equation}
\label{eq_travelling_waves}
\left\{
\begin{array}{lll}
     u_x = p \\
     p_x = \big( A( u ) - \sigma I )p \\
     \sigma_x=0
\end{array}
\right.
\end{equation}
and an equilibrium point $(u^{\ast}, \, 0, \, \lambda_i
(u^{\ast}))$. The center manifold theorem ensures that the center
space $V^c = \big\{ p=0 \big\}$ parameterizes a center manifold
$\mathcal{M}^c$. This manifold contains all the solutions of
\eqref{eq_travelling_waves} that do not diverge exponentially
neither as $x \to - \infty$ nor as $x \to + \infty$.

It can be shown that the center manifold $\mathcal{M}^c$ around
the equilibrium $(u^{\ast}, \, 0, \, \lambda_i (u^{\ast}))$ is
described by a function $v_i \tilde{r}_i(u, \, v_i, \, \sigma_i).$
Since $A$ is triangular, one can take
\begin{equation*}
      \tilde{r}_1(u, \, v_1 , \, \sigma_1) =
      \left(
      \begin{array}{cc}
           1 \\
           m(u, \, v_1, \, \sigma_1)
      \end{array}
      \right),
      \qquad
      \tilde{r}_2 (u, \, v_2, \, \sigma_2)
      \equiv \left(
             \begin{array}{ll}
                    0 \\
                    1
             \end{array}
             \right),
\end{equation*}
for some suitable function $m$ (in general different from the
function $f$ in the vector $\hat r_1$). One can moreover show that
the following equations hold:
\begin{equation*}
      \begin{split}
      &    \qquad \qquad
           A(u) \tilde{r}_1 = \lambda_1 \tilde{r}_1 +
           v_1 \mathrm{D} \tilde{r}_1 \tilde{r}_1 +
           v_1 ( \lambda_1 - \sigma_1 ) \tilde{r}_{1 v}, \\
      &    \tilde{r}_1 (u, \, 0, \, \sigma_1) = r_1 (u)
           \quad \forall \, u, \; \sigma_1,
           \quad
           | \tilde{r} (u, \, v_1, \, \sigma_1) - r_1 (u)| =
           \mathcal{O}(1) v_1,
           \quad
           \tilde{r}_{1 \, \sigma} = \mathcal{O}(1) v_1. \\
      \end{split}
\end{equation*}
Here and in the following we will denote by $( \ell_1, \; \:
\tilde{\ell}_2)$ the dual base of $(\tr, \; \: r_2)$.

\subsubsection{Gradient decomposition}
\label{par_gradient_decomposition} We set
\begin{equation}
\label{decomposition}
      \left\{
       \begin{array}{ll}
       u_x = v_1 \tilde{r}_1 ({u,\, v_1, \, \sigma_1}) +
                   v_2 r_2 +
                   p_1 \hat{r}_1 (u, p_1 ) +
                   p_2 r_2 \\
       &   \\
        u_t = w_1 \tilde{r}_1 ({u,\, v_1, \, \sigma_1}) +
                   w_2 r_2
       \end{array}
       \right.
       \quad   \sigma_1 = \lambda_1(u^{\ast})
                              - \theta
                              \bigg(
                                    \frac{ w_1}{v_1} + \lambda_1(u^{\ast})
                              \bigg).
\end{equation}
The function $\theta$ is here and in the following an odd cutoff
such that
\begin{equation}
\label{eq_theta}
      \theta(s) =
      \left\{
      \begin{array}{lll}
             s \quad \quad \textrm{if} \; |s|\leq \hat{\delta}     \\
             0 \quad \quad \textrm{if} \; |s|\geq 3 \hat{\delta}    \\
             \textrm{smooth connection if}
                     \quad \hat{\delta} \leq s \leq 3 \hat{\delta}
      \end{array}
      \right.
      \delta_1 < < \hat{\delta} \leq \frac{1}{3}.
\end{equation}
The choice of the speed $\sigma$ follows from the analysis of the
boundary free case, \cite{BiaBrevv}.

Note that \eqref{decomposition} is a system of 4 equations in 6
unknowns: as we underlined in the introduction, this will allow
some freedom in choosing the boundary conditions for $v_i, \;
i=1,2$ and $p_i, \; i=1,2$. More precisely, we will proceed as
follows.
\begin{enumerate}
\item We will insert \eqref{decomposition} in the parabolic
equation \eqref{rescaled}. This will generate a system of $4$
equations in $6$ unknown. \item We will obtain the equations for
$v_i, \; w_i, \; p_i, \; i=1, \; 2$ by assigning in a suitable way
the terms obtained. \item We will impose boundary and initial
conditions on each of the $6$ equations obtained. This procedure
selects one and only one solution for each of those equations.
\end{enumerate}
The decomposition \eqref{decomposition} is thus complete. We
observe that the idea is to let the equations to choose the
components in the decomposition, by only imposing reasonable
initial-boundary conditions and by assigning carefully the terms
obtained by inserting \eqref{decomposition} in the system
\eqref{rescaled}.

\subsection{The equations satisfied by
$ \boldsymbol{ v_i, \; p_i, \; w_i \; \; i=1, \, 2}$}
\label{par_equations}

These equations are obtained via the computations in Appendix
\ref{explicit_source_t}: inserting the components $v_i, \; p_i \;
w_i,$ $i=1, \, 2$ in the equation
\begin{equation*}
      u_t + A(u) u_x - u_{xx} =0,
\end{equation*}
we find
\begin{equation*}
\begin{split}
&      v_{1 t } +(\lambda_1 v_1)_x - v_{1 xx}+
       p_{1 t } +(\lambda_1 p_1)_x - p_{1 xx} = 0 \\
&      v_{2 t } +(\lambda_2 v_2)_x - v_{2 xx}+
       p_{2 t } +(\hat{\lambda}_2 p_2)_x - p_{2 xx} = \tilde{s}_1(t, \, x)\\
&      w_{1 t } +(\lambda_1 w_1)_x - w_{1 xx}=0 \\
&      w_{2 t } +(\lambda_2 w_2)_x - w_{2 xx} = \tilde{s}_2(t, \, x) \\
\end{split}
\end{equation*}
for some function $\tilde{s}_i(t, \, x) \; i=1, \, 2$ whose
explicit expression can be found in the appendix. Moreover, as it
is shown in the Appendix \ref{explicit_source_t}, from the
equation
\begin{equation*}
      u_t = u_{xx} - A(u) u_x
\end{equation*}
one gets the relations
\begin{equation}
\label{eq_vt}
\begin{split}
&     w_1 = v_{1 x} - \lambda_1 v_1 + p_{1 x} - \lambda_1 p_1  \\
&      w_2 = v_{2 x} - \lambda_2 v_2 + p_{2 x} - \hat{\lambda}_2
      p_2 + e(t, \, x) \\
\end{split}
\end{equation}
for a suitable error term $e(t,  \, x)$. The following Proposition
(whose proof can be found in the Appendix 
\ref{source_pro}) gives the form of the source terms:
\begin{pro}
\label{reasons_of_source_term} The following estimate holds:
\begin{equation}
\label{reasons_of_source_term_eq}
      \begin{split}
      | \tilde{s}_1&(t, \, x) |, \; | \tilde{s}_2 (t, \, x)| , \;
      |e(t, \, x)| \leq
                \mathcal{O}(1) \Biggl\{ \sum_{i \neq j}
              \Big[ |v_i| \Big( |v_j |+|v_{j x}| + |w_j| +
              |w_{j x} | \Big) +
              |w_i| \Big( |w_j|+ |v_{j x}| \Big) \Big] \\
      & + \sum_{i, \, j} \Big( |p_i| + |p_{ix }| \Big) \Big(|v_j| +
      |v_{jx }| + |w_j| + |w_{j x}| \Big)
      + | p_{1x} - \lambda_1 p_1 |
                        \Big( | p_{1x} | + |p_2| \Big)  \\
&      + \Big| w_1 v_{1 x} - v_1 w_{1 x} \Big|
      + v_1 ^2 \bigg|
      \bigg(  \frac{w_1}{v_1}
      \bigg)_x  \bigg|^2 \chi_{ \{ |w_1| \leq \delta_1 |v_1| \} }
      +  | w_1 + \sigma_1 v_1 |
       \Big( |v_1|+|v_{ 1 x }| + | w_1| + |w_{1 x}| \Big) \Biggr\}.
     \end{split}
\end{equation}
\end{pro}

Following the denomination of \cite{AnBia},
we will denote the above terms as follows:
\begin{enumerate}
\item interaction between waves of family 1 and family 2
\[
\sum_{i \neq j}
              |v_i|\Big( |v_j |+|v_{j x}| +
              |w_j| + |w_{j x} |\Big) +
              |w_i| \Big(|w_j|+ |v_{j x}| \Big);
\]
\item interaction of travelling waves with boundary profiles
\[
\sum_{i, \, j}\Big( |p_i| + |p_{ix }|\Big)\Big(|v_j| +
      |v_{jx }| + |w_j| + |w_{j x}|\Big);
\]
\item interaction among boundary profiles
\[
| p_{1x} - \lambda_1 p_1 |
                        \Big( | p_{1x} | + |p_2| \Big) ;
\]
\item $\sigma_1$ is not constant
\[
| w_1 v_{1 x} - v_1 w_{1 x}| + v_1 ^2 \bigg|
\bigg(  \frac{w_1}{v_1}
   \bigg)_x  \bigg|^2 \chi_{ \{ |w_1| \leq \delta_1 |v_1| \} };
\]
\item the cutoff function $\theta$ is active
\[
  | w_1 + \sigma_1 v_1 |
                         \Big(|v_1|+|v_{ 1 x }| + | w_1| + |w_{1 x}|\Big).
\]
\end{enumerate}
Since the component $p_2$ of the boundary profile should remain
close to the boundary $x=L$, and the source $\tilde s_1$ is in
general spread in the whole interval $[0,L]$, we split the
previous expression as follows:
\begin{equation*}
\begin{split}
&      v_{1 t } +(\lambda_1 v_1)_x - v_{1 xx} = 0
       \qquad \qquad \qquad \quad
       p_{1 t } +(\lambda_1 p_1)_x - p_{1 xx} = 0 \\
&      v_{2 t } +(\lambda_2 v_2)_x - v_{2 xx} =
       \tilde{s}_1 (t, \, x) \qquad \qquad \, \,
       p_{2 t } +(\hat{\lambda}_2 p_2)_x - p_{2 xx} = 0 \\
\end{split}
\end{equation*}

\subsection{Boundary conditions}
\label{boundary_conditions}

To conclude the characterization of the equations satisfied by
$v_i$, $p_i$, $w_i$, we have to assign the boundary conditions.
The basic idea is that each component $v_i$, $p_i$, $i=1,$ should
behave like a travelling wave or a boundary profile, respectively.
More precisely, we can make the following observations:
 \vskip .2cm \noindent 1) In order to behave like a double
 boundary profile, $p_1$ and $p_2$ should be independent from the
 initial datum, hence we are led to impose
\begin{equation*}
      p_1 (0, \, x) \equiv 0, \; \; p_2 (0, \, x) \equiv 0.
\end{equation*}
It follows that the initial data for $v_1$ and $v_2$ are given by 
\begin{equation*}
      v_1 (0, \, x) = \langle \ell_1, \; u'_0(x) \rangle
      \qquad
      v_2 (0, \, x) = \langle \tilde{\ell}_2, \, u'_0(x) \rangle .
\end{equation*}
\vskip .2cm \noindent 2) To emulate the behavior observed in the
hyperbolic limit, the waves of the first family should disappear
when hitting the boundary $x=0$, and the waves of the second
family should disappear at $x=L$. To understand what kind of
boundary condition it is convenient to impose, one can observe
that an integration by parts leads to
\begin{equation*}
\begin{split}
&     \frac{d}{dt} \int_0^L |v_1 (t, \, x)| dx =
      \int_0^L \text{sign} v_1
      \Big( v_{1 x} - \lambda_1 v_1 \Big)_x dx \\
&     \qquad \qquad \qquad \qquad =
      \int_0^L \delta_{v = 0} (v_{1 x} - \lambda_1 v_1 ) dx +
      \bigg[ \text{sign} v_1 (v_{1x} - \lambda_1 v_1 )\bigg]^L_0 \leq
      \bigg[ \text{sign} v_1 (v_{1x} - \lambda_1 v_1 )\bigg]^L_0, \\
&     \frac{d}{dt} \int_0^L |v_2 (t, \, x)| dx \leq
      \int_0^t \int_0^L |\tilde{s}_1 (s, \, x)| ds dx +
      \bigg[ \text{sign} v_2 (v_{2x} - \lambda_2 v_2 )\bigg]^L_0. \\
\end{split}
\end{equation*}
(we have used the inequality $\delta_{v =0} v_x \leq 0$). To
minimize the increment of $\|v_i(t)\|_{L^1}$ due to the
interactions with the boundary we impose
\begin{equation*}
      v_1 (t, \, 0) \equiv 0,
      \qquad
      v_2 (t, \, L) \equiv 0,
\end{equation*}
and integrating with respect to $t$ the previous equations we get
\begin{equation}
\label{estimate_v}
\begin{split}
&      \int_0^L |v_1 (t, \, x)| dx \leq \int_0^L |v_1(0, \, x)| dx
       + \int_0^t |v_{1 x} - \lambda_1 v_1|(s, \, L) ds, \\
&      \int_0^L |v_2 (t, \, x)| dx \leq \int_0^L |v_2(0, \, x)| dx
       + \int_0^t \int_0^L |\tilde{s}_1 (s, \, x)| ds dx
       + \int_0^t |v_{2 x} - \lambda_2 v_2|(s, \, 0) ds.\\
\end{split}
\end{equation}
We have used the following observations:
\begin{equation}
\label{observation_sign}
\begin{split}
&     v_1 (0)=0 \implies \lim_{x \to 0^+} \text{sign}( v_1) v_{1x}(x)
      \ge 0 \\
&     v_2 (L)=0 \implies \lim_{x \to L^-} \text{sign}(v_2) v_{2x}(x)
      \leq 0 . \\
\end{split}
\end{equation}
If one inserts the previous Dirichlet condition on $v_i \; \; i=1,
\, 2$ in the decomposition \eqref{decomposition}, obtains the
followings boundary conditions for $p_i$:
\begin{equation}
\label{bc}
      p_1 (t, \, 0) =  \langle \ell_1, \, u_x (t, \, 0) \rangle ,
      \qquad
      p_2 (t, \, L) = \langle \tilde{\ell}_2, \, u_x (t, \, L)  \rangle  -
                      p_1 \langle \tilde{\ell}_2, \, \hat{r}_1  \rangle .
\end{equation}
\vskip .2cm \noindent 3) Since $p_1$ should be located near $x=0$,
and $p_2$ near $x=L$, we would like to impose that the increment
of $\| p_1 \|_{L^1}$ due to the datum at $x=L$ is minimal, and
similarly that the increment of $\| p_2 \|_{L^1}$ caused by the
boundary datum in $x=0$ is as low as possible. Since the values
$p_1(t, \, 0)$ and $p_2(t, \, L)$ are already determined, we will
impose on
$p_1$ some condition at $x=L$ and on $p_2$ at $x=0$. \\
We observe that an integration by parts like the ones performed
before leads to
\begin{equation*}
      \int_0^L |p_1 (t, \, x)| \leq \int_0^t |p_{1 x} - \lambda_1 p_1
      |(s, \, 0) ds + \int_0^t |p_{1 x} - \lambda_1 p_1|(s, \, L)
      ds.
\end{equation*}
Hence we are led by the previous considerations to impose
\begin{equation}
\label{eq_condition}
      (p_{1 x} - \lambda_1 p_1 ) (t, \, L) \equiv 0.
\end{equation}
Similarly, we impose
\begin{equation}
\label{eq_condition2}
      (p_{2 x} - \hat{\lambda}_2 p_2 )(t, \, 0) \equiv 0.
\end{equation}
From these two equations we obtain the boundary conditions for
$v_1$, $v_2$: indeed, we have
\[
  \Big( v_{1,x} - \lambda_1 v_1 \Big)(t, \, L) = \langle \ell_1, \,
  u_{t}(t, \, L)
  \rangle
\]
and
\[
   \Big( v_{2,x} - \lambda_2 v_2 \Big)(t, \, 0) =
   \langle  \tilde{\ell}_2, \, u_t (t, \, 0) \rangle  - e(t, \, 0).
\]
\vskip .2cm

At this point, the initial-boundary data
are perfectly determined for all the components
$v_i$, $p_i$, $i=1,2$, and thus the decomposition is complete.


\section{BV estimates}
\label{BV_estimates}

Aim of this section is to prove the following theorem, which
constitutes the first part of Theorem \ref{main_result}.

\begin{teo}
\label{BV}
      Let $u(t, \, x)$ be the local in time solution of the $2 \times 2$ system
\begin{equation}
\label{rescaled1}
      \left\{
      \begin{array}{lllll}
            u_t + A (u) u_x =
            u_{xx} \\
            \\
            u(0, x) = {u}_0 (x)\\
            \\
            u(t, 0) = {u}_{b \, 0}(t) \qquad
            u(t, L) = {u}_{b L}(t) \\
     \end{array}
     \right.
\end{equation}
      and suppose that the boundary and initial conditions are regular
      and
      satisfy
      \begin{equation*}
            \bigg\| \frac{d^k u_0}{dx^k} \bigg\|_{L^1(0, \, L)},
            \bigg\| \frac{d^k u_{b \, 0}}{dt^k}
            \bigg\|_{L^1(0, \, + \infty)},
            \bigg\| \frac{d^k u_{b \, L}}{dt^k}
            \bigg\|_{L^1(0, \, + \infty)}
            \leq \delta_1 \quad k = 1, \dots n,
      \end{equation*}
      for some $\delta_1$ sufficiently small.

      Then $u(t, \, x)$ is defined $\forall \; t > 0$
      and its total variation is uniformly bounded:
      \begin{equation}
      \|u_x(t) \|_{L^1(0, \, L)} \leq C  \delta_1
      \end{equation}
      for some constant $C$ which does not depend on $L$.
\end{teo}

It is enough to prove that there is a constant
$\delta_0$ such that $k \delta_1 \leq \delta_0 <<1 $ with $k$
small enough and such that the following holds: if $\delta_1$ is
small enough and $ \|u_x ( s)\|_{L^1} \leq C \delta_1 \; \forall s
\in \, [0, \, t]$ then
\begin{equation}
\label{observation_BV}
      \begin{array}{ll}
      {\displaystyle \int_0^t \int_0^L |\tilde{s}_1 (\sigma, \, x) | dx
      d \sigma
             \leq \unpo \delta_0^2},
      &
      {\displaystyle        \int_0^t \int_0^L |\tilde{s}_2 (\sigma, \, x) | dx
      d\sigma
             \leq \unpo \delta_0^2, }\\
      {\displaystyle       \int_0^{t}|v_{2 x} - \lambda_2 v_2| (\sigma, \, 0)
             d\sigma \leq m \delta_1,}
      &
      {\displaystyle        \int_0^{t}|v_{1x} - \lambda_1 v_1 | (\sigma, \, L)
             d\sigma \leq m \delta_1,} \\
      {\displaystyle       \int_0^{t} |p_{1 x} - \lambda_1 p_1|
             (\sigma, \, 0) d\sigma
             \leq m \delta_1,}
      &
      {\displaystyle        \int_0^{t} |p_{2 x} - \hat{\lambda}_2 p_2|
             (\sigma, \, 0) d\sigma \leq m \delta_1,} \\
      \end{array}
\end{equation}
for some constant $m$ that does not depend on $C$.

Indeed, suppose the previous implication holds. From the
representation formula \eqref{solution_x} it immediately follows
that the function $t \mapsto \|u_x( t) \|_{L^1}$ is continuous:
hence, it will satisfy $\|u_x (t)\|_{L^1} < C \delta_1$ if $t$ is
small enough, since the total variation of the initial datum is
bounded by $\delta_1$.

Suppose by contradiction that $\tau$ is the first time such that
$\|u_x (\tau)\|_{L^1} = C \delta_1$. Then we use the equations
\begin{equation*}
\begin{split}
&      v_{1 t } +(\lambda_1 v_1)_x - v_{1 xx}=0
       \qquad \qquad \; \; \;
       p_{1 t } +(\hat{\lambda}_1 p_1)_x - p_{1 xx} = 0 \\
&      v_{2 t } +(\lambda_2 v_2)_x - v_{2 xx} = \tilde{s}_1(t, \,
x)
       \qquad
      p_{2 t } +(\hat{\lambda}_2 p_2)_x - p_{2 xx} =
       0 \\
\end{split}
\end{equation*}
and the boundary conditions described in Section
\ref{boundary_conditions} and, integrating by parts, we get
\begin{equation*}
\begin{split}
     \int_0^L |u_x (\tau, \, x)| dx \leq
&    \sum_{i = 1}^2 \int_0^L |v_i (\tau, \, x) | +
     \int_0^L |p_i (\tau, \, x) | dx \leq
     \sum_{i= 1}^2 \int_0^L |v_i (0, \, x) | +
     \int_0^{\tau} \int_0^L |\tilde{s}_1 (\sigma, \, x)| dx d\sigma \\
&    \quad +
     \int_0^{\tau} |v_{2 x } -
     \lambda_2 v_2|(\sigma, \, 0) d\sigma +
     \int_0^{\tau} |v_{1x } -
     \lambda_1 v_1|(\sigma, \, L) d\sigma +
     \int_0^{\tau}  |p_{1x } - \lambda_1 p_1|(\sigma, \, 0)
     d\sigma  \\
&    \quad +
     \int_0^{\tau} |p_{2x } -
     \hat{\lambda}_2 p_2|(\sigma, \, L) d\sigma
     \leq
     (4 m + 2 ) \delta_1+ \unpo \delta_0^2
     < C \delta_1,
     \phantom{\int}
\end{split}
\end{equation*}
if $C$ is large enough: this contradicts the assumption
$\|u_x (\tau)\|_{L^1} = C \delta_1$.
%

Note that since all the functions in the right hand side of
\eqref{reasons_of_source_term_eq} are continuous (and hence
bounded on $[0, \, L]$), we have that
\begin{equation}
\label{eq_bound_source}
      \int_0^s \int_0^L |\tilde{s}_i (\sigma, \, x)| dx d \sigma
      \leq \unpo \delta_1 \quad i=1, \, 2,
\end{equation}
for $s$ small enough. Hence to prove \eqref{observation_BV} we can
suppose that \eqref{eq_bound_source} holds for any $s \in [0, \,
t]$: since we will show that actually
\begin{equation*}
      \int_0^t \int_0^L |\tilde{s}_i (\sigma, \, x)| dx d \sigma
      \leq \unpo \delta_0^2, \quad i=1, \, 2,
\end{equation*}
the assumption will be a posteriori justified since
$k \delta_1 \leq \delta_0 <<1$.

We will proceed as follows: in Section
\ref{elementary_estimates} we will show some elementary estimates,
while in Section \ref{interaction_functionals} we will introduce
suitable functionals that allow the estimates
\begin{equation*}
\begin{split}
&      \int_0^t \int_0^L \sum_{i \neq j}
       \big( |v_i|( |v_j |+|v_{j x}| + |w_j| + |w_{j x} |) +
           |w_i| (|w_j|+ |v_{j x}| ) \big) (\sigma \, x)d\sigma dx \leq
       \unpo \delta_1^2, \\
&      \int_0^t \int_0^L
       | w_1 v_{1 x} - v_1 w_{1 x}| (\sigma, \, x) d\sigma \leq \unpo
       \delta_1^2, \\
&      \int_0^t \int_0^L \bigg| v_1 ^2
       \bigg( \frac{w_1}{v_1}\bigg)_x
       \bigg|^2 \chi_{ \{ |w_1| \leq \delta_1 |v_1| \} }(\sigma, \, x)
       d\sigma dx \leq \unpo \delta_1^2 .\\
\end{split}
\end{equation*}
In Section \ref{energy_estimates} we will consider the term
\begin{equation*}
      \int_0^t \int_0^L
      | w_1 + \sigma_1 v_1 |
      (|v_1|+|v_{ 1 x }| + | w_1|)(\sigma, \, x) d\sigma dx,
\end{equation*}
and prove a bound of order $\delta_1^2$.
\subsection{Elementary estimates}
\label{elementary_estimates}

This section is devoted to the estimates which can be obtained by
elementary techniques, like the maximum principle. We will in
particular show that the components $p_i$, $i=1,2$ are
exponentially decaying as one moves far away from the boundary,
and that their decay exponent does not depend on the interval
length $L$. Moreover, by introducing various functional, we
estimate the boundary data assigned to the components $v_1$, $v_2$
and prove that the functions $v_i$ are integrable along all
vertical lines $\{x = \text{const}\}$. This means that, as in the
boundary free case, the profiles of travelling waves just cross
the vertical lines.

\subsubsection{Estimates via maximum principle}
\label{par_maximum}


We will first deal with $p_1$. The results in Section
\ref{sub_paraboli_estimates} ensures that
\begin{equation*}
      \|u_x (t)\|_{L^{\infty}} \leq \| u_{xx}(t)\|_{L^1} \leq
      \unpo \delta_1.
\end{equation*}
Hence it follows that
\begin{equation*}
      |p_1 (t, \, 0)| = |\langle l_1 , \, u_x (t, \, 0 )  \rangle | \leq k
      \delta_1,
\end{equation*}
for some $k$ large enough.

The equation satisfied by $p_1$ is
\begin{equation*}
     p_{1 t} + \lambda_1(u) p_{1 x} + \lambda_{1 x} (u) p_1 -
     p_{1 xx}=0.
\end{equation*}
This is a linear equation, with coefficients depending on the
solution $u(t,x)$. Let $2c$ be the separation speed defined in
\eqref{eq_separation_speed} and
\begin{equation*}
      q(x) = k \delta_1 \exp \big( - c x /2 \big).
\end{equation*}
Since $|\lambda_{1 x}| \leq \unpo \delta_1$ and $\delta_1 <<1$,
$q$ satisfies
\begin{equation*}
      q_t + \lambda_1 q_x + \lambda_{1 x} q - q_{xx} > 0.
\end{equation*}
Hence the difference $(q - p_1 )$ satisfies
\begin{equation*}
\left\{
\begin{array}{lll}
      (q- p_1)_t + \lambda_1 ( q - p_1)_x +
       \lambda_{1 x} (q - p_1) -
        (q_{xx} - p_{1 xx}) > 0  \phantom{\bigg(} \\
      (q - p_1 )(t, \, 0) \ge 0 \phantom{\bigg(} \\
      \bigg( (q - p_1 ) - \lambda_1 (q -p_1 )_x \bigg)(t, \, L) >
      0.
\end{array}
\right.
\end{equation*}
By standard techniques it follows that $(q - p_1)(t, \, x) \ge 0$
for any $t, \; x$ and hence
\begin{equation}
\label{estimate_exp_decay1}
      |p_1(t, \, x)| \leq k \delta_1 \exp(- c \, x /2 ).
\end{equation}
The boundary condition on $p_2$ satisfies the following bound:
\begin{equation*}
       |p_2(t, \, L)| = |\langle  \hat{l}_2, \, u_x(t, \, L) \rangle
       - p_1 \langle  \tilde{l}_2, \, \hat{r}_1  \rangle  | \leq \unpo \delta_1,
       \qquad \forall \, t, \; x.
\end{equation*}
Since $|p_1(t, \, x)| \leq k \delta_1$, then from
\eqref{eq_gen_eigen} it follows that $|\lambda_2 -
\hat{\lambda}_2| \leq \unpo \delta_1$ and hence in the same way as
before one can prove
\begin{equation}
\label{estimate_exp_decay2}
      |p_2 (t, \, x)| \leq \unpo \delta_1 \exp (c(x - L)/2 ),
      \qquad \forall \, t, \; x.
\end{equation}
From \eqref{estimate_exp_decay1} it follows
\begin{equation*}
      \|p_1 (t)\|_{L^1} \leq \unpo \delta_1,
      \qquad
      \|v_1(t)\|_{L^1} \leq \unpo \delta_1
\end{equation*}
and, since $\|u_x\|_{L^{\infty}} \leq \unpo \delta_1$,
\begin{equation*}
      \|v_1\|_{L^{\infty}} \leq \unpo \delta_1.
\end{equation*}
Analogously, from \eqref{estimate_exp_decay2} it follow
\begin{equation*}
      \|p_2 (t)\|_{L^1} \leq \unpo \delta_1,
      \qquad
      \|v_2(t)\|_{L^1} \leq \unpo \delta_1,
      \qquad
      \|v_2\|_{L^{\infty}} \leq \unpo \delta_1.
\end{equation*}
The following proposition summarizes the results obtained in this
paragraph:
\begin{pro}
\label{exp_decay}
      Let $p_i, \; v_i$ be the solutions of \eqref{cons_laws} with
      the boundary conditions described in Section
      \ref{boundary_conditions}. Then
      \begin{equation*}
            |p_1(t, \, x) | \leq \unpo \delta_1 \exp (- cx / 2),
            \qquad
            |p_2(t, \, x) | \leq \unpo \delta_1 \exp (c ( x -L)/2),
      \end{equation*}
      where $2c$ is the separation speed defined by
      \eqref{eq_separation_speed}.

      The previous estimates imply
      \begin{equation*}
            \|p_i (t)\|_{L^1} \leq \unpo \delta_1,
            \qquad
            \|v_i (t)\|_{L^1} \leq \unpo \delta_1,
            \qquad
            \|v_i(t)\|_{\infty} \leq \unpo \delta_1,
            \quad i=1, \, 2.
      \end{equation*}
\end{pro}
\begin{rem}
\label{rem_other_way} The estimate of $\|v_i(t)\|_{L^1}$ can also
be obtained directly from \eqref{estimate_v}: indeed, since
\begin{equation*}
       {(p_{1x} - \lambda_1 p_1)(t, \, L) \equiv 0}
\end{equation*}
and the total variation of $u_{b \, L}$ is bounded by $\delta_1$,
from \eqref{eq_vt} one gets
\begin{equation*}
      \int_0^t |v_{1 x} - \lambda_1 v_1 |(s, \, L) ds \leq
      \delta_1,
\end{equation*}
and hence $\|v_1 (t)\|_{L^1} \leq 2 \delta_1$.

To obtain the estimate on $v_2$ from \eqref{estimate_v} one has to
start supposing
\begin{equation}
\label{estimate_assumption_e}
      \int_0^t |e(s, \, 0)| ds \leq  \delta_1.
\end{equation}
With the same computations as before one gets $\|v_2 (t)\|_{L^1}
\leq \unpo \delta_1$. As it will be clear from the next sections,
the assumption \eqref{estimate_assumption_e} actually leads to the
estimate
\begin{equation*}
      \int_0^t |e(s, \, 0)| ds \leq \unpo \delta_1^2,
\end{equation*}
and therefore it is a posteriori well justified.
\end{rem}

\subsubsection{Integrability with respect to time}

The following lemma, which can be proved by a simple integration
by parts, introduces a useful estimate we will widely use in the
following.
\begin{lem}
\label{functional_pro}
      Let $P(x)$ be a non negative $\mathcal{C}^2$
      function defined on $\mathbb{R}$
      and let $q$ be a solution of
      \begin{equation*}
            q_t + (\lambda q)_x - q_{xx}= s(t, \, x).
      \end{equation*}
      Then the following estimate holds:
      \begin{equation*}
      \begin{split}
             \frac{d }{dt } \int_0^L |q(t, \, x)| P(x) dx \leq
      &      \int_0^L |s(t, \, x) | P(x) dx +
             \int_0^L |q(t, \, x)| ( \lambda P' + P'')( x)
             dx \\
      &      -
             \bigg[ P'
             |q(t)|\bigg]^{x =L }_{x =0} +
             \bigg[ P \,
             \text{\rm sign} (q) \, (q_x - \lambda q)(t)
             \bigg]^{x = L}_{x =0}. \\
      \end{split}
      \end{equation*}
\end{lem}

Before applying the previous lemma, we recall that the boundary
data of the scaled problem \eqref{rescaled} belongs to $BV(0, \, +
\infty)$ and that the $L^1$ norms of $u'_{b \, 0}$ and $u'_{b \,
L}$ are bounded by $\delta_1$. From the decomposition $u_t = w_1
\tr + w_2 r_2$, we immediately have
\begin{equation*}
       \|w_i ( x = 0)\|_{L^1 (0, \, + \infty)} \leq
       \delta_1
       \qquad
       \|w_i ( x = L)\|_{L^1 (0, \, + \infty)} \leq
       \delta_1
       \quad i=1, \, 2.
\end{equation*}
Moreover, in Section \ref{par_equations} we found that $w_i \;
i=1, \, 2$ can be decomposed as follows:
\begin{equation}
\label{eq_decomposition_w}
\begin{split}
&     w_1 = p_{1x} - \lambda_1 p_1 + v_{1 x} - \lambda_1 v_1 \\
&     w_2 = p_{2x} - \hat{ \lambda}_2 p_2 + v_{2 x} - \lambda_2 v_ 2 +
      e(t, \, x),\\
\end{split}
\end{equation}
where the error term $e(t, \, x)$ satisfies the estimate
\eqref{reasons_of_source_term_eq}. As we anticipated in Remark
\ref{rem_other_way}, we will suppose
\begin{equation*}
       \int_0^t |e(s, \, x)|ds \leq  \delta_1
       \quad \forall \, x \, \in \, [0, \, L].
\end{equation*}
Since we will obtain an estimate of order $\delta_1^2 \leq
\delta_1$, this
assumption is a posteriori well justified.

From the boundary condition \eqref{eq_condition2} $(p_{2 x} -
\hat{\lambda}_2 p_2)(t, \, 0) \equiv 0$ and from the decomposition
\eqref{eq_decomposition_w} we get
\begin{equation*}
      \int_0^t |v_{2 x} - \lambda_2 v_2 |(s, \, 0) ds \leq 2
      \delta_1.
\end{equation*}
Similarly, one obtains that
\begin{equation*}
      \int_0^t |v_{1 x} - \lambda_1 v_1|(s, \, L) ds \leq
      \delta_1.
\end{equation*}

An application of Lemma \ref{functional_pro} with $P \equiv 1$ and
$q= v_2$ leads by observation \eqref{observation_sign} to
\begin{equation*}
\begin{split}
     \int_0^t |v_{2 x}(s, \, L)| ds \leq&~
      \int_0^t \int_0^L |\tilde{s}_2 (s, \, x)| dx ds +
      \int_0^t |v_{2x } - \lambda_2 v_2| (s, \, 0) ds +
      \int_0^L |v_2 (0, \, x)| dx \\
     \leq&~ \unpo \delta_1 + 2 \delta_1 +
      \unpo \delta_1 \leq \unpo \delta_1, \\
\end{split}
\end{equation*}
and similarly
\begin{equation*}
      \int_0^t |v_{1 x}(s, \, 0)| ds
      \leq \unpo \delta_1.
\end{equation*}
Let $2c$ be the separation speed defined by
\eqref{eq_separation_speed}: the application of Lemma
\ref{functional_pro} with $q(t, \, x)= v_2(t, \, x)$ and
\begin{equation*}
      P(x) = P_y(x)=
      \left\{
      \begin{array}{lll}
            1 / c & x \leq y \\
            \\
            \exp \Big(  c (y- x) \Big) / c
            & x > y \\
      \end{array}
      \right.
      \quad y \, \in \, [0, \, L[
\end{equation*}
leads to the estimate
\begin{equation*}
\begin{split}
      \int_0^t |v_2 (s, \, y)| ds
&     \leq
      \int_0^L |v_2 (0, \, x) | dx +
      \frac{1}{c} \int_0^t \int_0^L |\tilde{s}_1(s, \, x)| ds dx ~   \\
&     \quad
      + P_y (0)\int_0^L  | v_{2 x} - \lambda_2 v_2|(s, \, 0) ds +
      P_y(L) \int_0^t |v_{2 x}(s, \, L)| ds \\
&     \leq \unpo \delta_1 + \unpo \delta_1
      \leq \unpo \delta_1
      \qquad \forall \, y \in \, [ 0, \, L[.
      \phantom{\int}
\end{split}
\end{equation*}
Analogously, we get
\begin{equation*}
      \int_0^t |v_1(s, \, y) | ds \leq \unpo \delta_1
      \quad \forall \, y \, \in \, ]0, \, L].
\end{equation*}

The following proposition summarizes what we have proved so far:
\begin{pro}
\label{functional_estimates_pro}
       Let $v_i, \; p_i \; \; i=1, \, 2$ be the solutions to the
       equations \eqref{cons_laws} with the boundary conditions
       described in Section \ref{boundary_conditions}. Then it
       holds
       \begin{equation*}
       \begin{split}
            \int_0^t |v_{2 x } - \lambda_2 v_2|(s, \, 0) ds
             \leq 2 \delta_1,
             & \qquad
             \int_0^t |v_{ 1 x } - \lambda_1 v_1|(s, \, L) ds
             \leq \delta_1, \\
            \int_0^t |v_{1 x}(s, \, 0)| ds \leq \unpo \delta_1,
       &      \qquad
             \int_0^t |v_{2 x} (s, \, L)|ds \leq \unpo \delta_1, \\
       \end{split}
       \end{equation*}
       and
       \begin{equation*}
              \int_0^t |v_i (s, \, y)| ds \leq \unpo \delta_1,
              \quad \forall \, y \, \in \, [0, \, L ]
              \quad i=1, \, 2.
       \end{equation*}
\end{pro}
Further computations (Appendix \ref{exp_decay_px_proof}) ensure
that
\begin{equation}
\label{exp_decay_px}
      |p_{1 x}(t, \, x)| \leq \unpo \delta_1 \exp (- cx / 2),
      \qquad
      |p_{2 x}(t, \, x)| \leq \unpo \delta_1 \exp \big( c( x -L)/ 2 \big).
\end{equation}

The following proposition deals with other estimates of
integrals with respect to time: the proof is quite long and requires
the introduction of new convolution kernels. It can be found
in the Appendix \ref{other_wrt_time_par}.
\begin{pro}
\label{other_wrt_time_pro}
      In the same hypothesis of Proposition
      \ref{functional_estimates_pro} it holds
      \begin{equation*}
            \int_0^t |v_{i x}(s, \, y)| ds \leq \unpo \delta_1
            \quad \forall \, y \, \in \, [0, \, L ]
            \quad i=1, \, 2
      \end{equation*}
      and
      \begin{equation*}
            \int_0^t |w_i (s, \, y)|ds \leq \unpo \delta_1
            \quad \forall \, y \, \in \, [ 0, \, L ]
            \quad i=1, \, 2.
      \end{equation*}
      We also have
      \begin{equation*}
            \int_0^t |w_{ix }(s, \, y)| ds \leq \unpo \delta_1
            \quad \forall \, y \, \in \, [0, \, L ]
            \quad i=1, \, 2.
      \end{equation*}
\end{pro}
In the previous proposition the functions $w_i$ are of course
defined by relation $u_t = w_1 \tilde{r}_1 + w_2 r_2$. Putting
together Proposition \ref{functional_estimates_pro} and
\ref{other_wrt_time_pro} and the decomposition
\eqref{eq_decomposition_w} one gets
\begin{equation*}
      \int_0^t |p_{1 x} - \lambda_1 p_1 |(s, \, y) ds \leq \unpo
      \delta_1,
      \qquad
      \int_0^t |p_{2 x} - \hat{\lambda}_2 p_2 |(s, \, y) ds \leq \unpo
      \delta_1,
      \quad \forall \, y \, \in \, [0, \, L ],
\end{equation*}
and
\begin{equation*}
        \int_0^t |p_{1 x} - \lambda_1 p_1 |(s, \, 0) ds \leq m
      \delta_1,
      \qquad
      \int_0^t |p_{2 x} - \hat{\lambda}_2 p_2 |(s, \, L) ds \leq m
      \delta_1,
\end{equation*}
where the constant $m$ satisfies the hypothesis stated in Section
\ref{BV_estimates}.

The estimates obtained so far will be widely used in next sections
and moreover allow to prove a bound of order $\unpo \delta_1^2$ on
some of the terms that appear on the right hand side of
\eqref{reasons_of_source_term_eq}:
\begin{equation}
\label{estimate_interaction}
\begin{split}
&     \int_0^t \int_0^L \sum_{i, \, j}( |p_i| + |p_{ix }|)(|v_j| +
      |v_{jx }| + |w_j| + |w_{j x}|) (s, \, x) ds dx \\
&     \qquad \leq
      \unpo \delta_1 \int_0^L ( e^{- cx} + e^{c (x -L)} )
      \int_0^t  (|v_j| +
      |v_{jx }| + |w_j| + |w_{j x}|) (s, \, x) ds dx
      \leq \unpo \delta_1^2 \\
\end{split}
\end{equation}
and
\begin{equation}
\label{estimate_boundary}
\begin{split}
&    \int_0^t \int_0^L \sum_{i}|p_{1 x}
     - \lambda_1 p_1|(|p_i| +
     |p_{ix}|)(s, \, x)  dx ds \\
&    \qquad \leq \unpo \delta_1 \int_0^L  e^{- c x} +
     e^{c ( x - L)}  \int_0^t |p_{1 x} - \lambda_1 p_1 |(s, \, x)
     ds dx \leq \unpo \delta_1^2. \\
\end{split}
\end{equation}

\subsection{Interaction functionals}
\label{interaction_functionals}

In this section we introduce three nonlinear functionals and we
use them to bound those terms in the right hand side of
\eqref{reasons_of_source_term} due to interaction between waves of
different families and those due to the fact that the speed
$\sigma_1$ is not constant. The form of the functionals is exactly
the same considered in \cite{BiaBrevv}, with some more
technicalities due to the presence of the boundary.

\subsubsection{Interaction among waves of different families}

We claim that the condition
\begin{equation*}
       \int_0^t \int_0^L |\tilde{s}_1 (s,  \, x) | ds dx \leq
       \unpo \delta_1
       \qquad
       \int_0^t \int_0^L |\tilde{s}_2 (s, \,x) | ds dx \leq
       \unpo \delta_1
\end{equation*}
implies
\begin{equation}
\label{interaction1}
       \int_0^t \int_0^L \sum_{i \neq j}
       \Bigg( |v_i| \Big( |v_j |+  |w_j| \Big) +
      |w_i w_j|  \Bigg) (s, \, x)ds dx \leq
      \unpo \delta_1^2.
\end{equation}
We will prove only that
\begin{equation}
\label{interaction1_case}
      \int_0^t \int_0^L |v_1 v_2 |(s, \, x) ds dx \leq \unpo
      \delta_1^2,
\end{equation}
because the other terms in \eqref{interaction1} can be dealt with
analogously: see for example \cite{AnBia}.

Let $2c$ be the separation speed introduced in
\eqref{eq_separation_speed} and let $P(\xi)$ be defined as
follows:
\begin{equation*}
  P(\xi) : =
  \left\{
  \begin{array}{ll}
       e^{c \xi} / 2 c
        \qquad \; \xi < 0 \\
       1 / 2c
        \qquad \; \; \; \;
        \xi \geq 0
  \end{array}
  \right.
\end{equation*}
One gets
\begin{equation*}
\begin{split}
     \frac{d}{ds}
&     \bigg( \int_0^L \int_0^L
      P(x -y) | v_1(s, \, x) | \, |v_2(s, \, y) | dx dy \bigg)
       \leq \int_0^L |v_2 (s, \, y)|
      \bigg[ P(x-y) \text{sign} v_1
      ( v_{1 x} - \lambda_1 v_1 ) (s, \, x) \bigg]^{ x = L}_{x=0} dy  \\
&    +   \int_0^L |v_1 (s, \, x)|
     \bigg[ P(x -y) \text{sign} v_2
     ( v_{2 x} - \lambda_2 v_2 ) ( s, \, y) \bigg]^{y=L}_{y=0} dx
      - \int_0^L | v_2 (s, \, y) |
     \bigg[ P'(x-y) |v_1 (s, \, x)| \bigg]^{ x=L}_{x=0} dy
     \\
&    +  \int_0^L |v_1 (s, \, x)|
     \bigg[ P'(x-y) |v_2 (s, \, y) | \bigg]^{y=L}_{y=0}
     + \int_0^L |v_1 (s, \, x) | \int_0^L P(x-y ) | \tilde{s}_1  (s, \,
     y)| dy \\
&    +
     \int_0^L \int_0^L
     \bigg( P'(x-y) \Big( \lambda_1 (s, \, x) - \lambda_2 (s, \, y) \Big) +
     2 P ''(x-y) \bigg) |v_1 (s, \, x)| \,| v_2 (s, \, y)| dx dy. \\
\end{split}
\end{equation*}
One has
\begin{equation*}
\begin{split}
P' ( \lambda_1 - \lambda_2 ) + 2 P'' \leq
       2 ( -c P' + P'') = - \delta_{ s =0}, \quad
       0 \leq  P (s) \leq \frac{1}{2c}, \qquad
       0 \leq  P'(s) \leq \frac{1}{2} \\
\end{split}
\end{equation*}
and moreover from the estimates of Proposition \ref{exp_decay} and
\ref{functional_estimates_pro} it follows that
\begin{equation*}
\begin{split}
&     \int_0^t |v_{1 x} - \lambda_1 v_1 |(s, \, L) \int_0^L
      |v_2 (s, \, y)| dy ds \leq \unpo \delta_1^2
      \qquad
      \int_0^t |v_{2 x} - \lambda_2 v_2 |(s, \, 0) \int_0^L
      |v_2 (s, \, x)| dx ds \leq \unpo \delta_1^2 \\
&     \qquad \qquad \qquad \qquad \qquad \qquad \qquad
      \int_0^t \int_0^L |\tilde{s}_1(s, \, y)| \int_0^L
      |v_1(s, \, x)| dx dy ds \leq \unpo \delta_1^2: \\
\end{split}
\end{equation*}
this completes the proof of the estimate
\eqref{interaction1_case}.

With some technical computations, in Appendix
\ref{interaction2_proof} it is proved
\begin{equation}
\label{interaction2}
      \int_0^t \int_0^L \sum_{i \neq j}
       \bigg( |v_i| \Big( |v_{j x}| +  |w_{j x} | \Big) +
           |w_i v_{j x}|  \bigg) (s, \, x)ds dx \leq
       \unpo \delta_1^2,
\end{equation}
which completes the proof of the estimate
\begin{equation*}
      \int_0^t \int_0^L \sum_{i \neq j}
      \bigg( |v_i| \Big( |v_j |+|v_{j x}| + |w_j| + |w_{j x} | \Big) +
           |w_i| (|w_j|+ |v_{j x}| ) \bigg) (s, \, x)ds dx \leq
      \unpo \delta_1^2.
\end{equation*}

\subsubsection{Length and area functionals}

To prove the estimate
\begin{equation*}
      \int_{0}^{t}
      \int_{0}^{L}
      | v_{1 x} w_1 - v_1 w_{1 x} |(s, \, x) dx ds
      \leq
      \mathcal{O}(1) \delta_1^2,
\end{equation*}
we introduce the curve
\begin{equation}
\label{gamma}
      \gamma (x) =
      \left(
      \begin{array}{cc}
      v_1( x ) \\
      w_1 ( x)
      \end{array}
      \right)
\end{equation}
and the related area functional
\begin{equation*}
     \mathcal{A}( \gamma )(s) =
      \frac{1}{2}
      \int \int _{ y \leq x}
      | \gamma_x \wedge
      \gamma_y | dx dy =
      \frac{1}{2}
      \int_{0}^{L}
      \int_{0}^{x}
      | v_1 (s, \, x)w_1(s, \, y) - v_1(s, \, y) w_1(s, \, x) |
      dx dy.
\end{equation*}
The curve $\gamma_x$ satisfies
\begin{equation*}
     \gamma_{xt} + (\lambda_1 \gamma_x)_x = \gamma_{xxx}
\end{equation*}
and moreover one has
\begin{equation*}
\begin{split}
&     \frac{d \mathcal{A}(s)}{ ds} =
      \frac{1}{2}
       \int_0^L \int_y^L \text{sign} \Big( v_1 (s, \, x) w_1
      (s, \, y) - v_1 (s, \, y) w_1 (s, \, x) \Big)
      \bigg( v_1 (s, \, x) w_1(s, \, y) -
      v_1 (s, \, y) w_1 (s, \, x) \bigg)_{xx} \\
&     \quad - \frac{1}{2}
       \int_0^L \int_y^L \text{sign} \Big( v_1 (s, \, x) w_1
      (s, \, y) - v_1 (s, \, y) w_1 (s, \, x) \Big)
      \bigg( \lambda_1(s, \, x) \Big( v_1 (s, \, x) w_1
      (s, \, y) - v_1 (s, \, y) w_1 (s, \, x)   \Big) \bigg)_x \\
&     \quad + \frac{1}{2}
      \int_0^L \int_0^x \text{sign} \Big( v_1 (s, \, x) w_1
      (s, \, y) - v_1 (s, \, y) w_1 (s, \, x) \Big)
      \bigg( v_1 (s, \, x) w_1(s, \, y) -
      v_1 (s, \, y) w_1 (s, \, x) \bigg)_{yy} \\
&     \quad - \frac{1}{2}
      \int_0^L \int_0^x \text{sign} \Big( v_1 (s, \, x) w_1
      (s, \, y) - v_1 (s, \, y) w_1 (s, \, x) \Big)
      \bigg( \lambda_1(s, \, y)  \Big( v_1 (s, \, x) w_1
      (s, \, y) - v_1 (s, \, y) w_1 (s, \, x)   \Big)\bigg)_y \\
\end{split}
\end{equation*}
and hence
\begin{equation*}
\begin{split}
      \frac{ d \mathcal{A} }{ ds} \leq
&     \frac{1}{2} \int_0^L
      \big|
           v_{1y}(s, \, L)w_1(s, \, y) -
           v_1(s, \, y)w_{1y}(s, \, L)
      \big|dy -
      \frac{1}{2} \int_0^L
      \big|
           v_{1 y}(s, \, y) w_1(s, \, y) - w_{1y}(s, \, y)v_1(s, \, y)
      \big|dy \\
&     - \frac{1}{2} \int_0^L
      \lambda_1 (s, \, L)
      \big|
           v_1(s, \, L) w_1(s, \, y)-
            v_1(s, \, y) w_1(s, \, L )
      \big| dy  \\
&     -
      \frac{1}{2} \int_0^L
      \big|
           v_1(s, \, x) w_{1 x}(s, \, x) - w_1(s, \, x) v_{1 x}(s, \, x)
      \big|dx  +
      \frac{1}{2} \int_0^L
      \big|
           v_{1}(s, \, x)w_{1 x}(s, \, 0) -
           v_{1 x } (s, \, 0 )w_1(s, \, x)
      \big|dx \\
\end{split}
\end{equation*}
Since $A(\gamma)(0) \leq \mathcal{O}(1) \delta_1^2$, one obtains,
using the estimates in Propositions \ref{functional_estimates_pro}
and \ref{other_wrt_time_pro},
\begin{equation*}
    \int_0^t \int_0^L
    \big|
          v_1(s, \, x)w_{1x}(s, \, x)-
          v_{1x}(s, \, x)w_1(s, \, x)
    \big|dx
    \leq
    -
    \int_0^t
    \frac{d\mathcal{A}}{ds} ds
    +
    \mathcal{O}(1) \delta_1^2
    \leq
    \mathcal{O}(1) \delta_1^2.
\end{equation*}

The length functional of the curve \eqref{gamma} is defined as
\begin{equation*}
       \mathcal{L}(\gamma)(s) =
       \int_0^L | \gamma_x | dx =
       \int_0^L \sqrt{ v_1^2 + w_1^2} dx,
\end{equation*}
and will be used to prove the estimate
\begin{equation}
\label{length_functional_eq}
      \int_0^t \int _0^L
      v_1^2 \bigg[ \bigg(
                 \frac{w_1}{v_1}
                 \bigg)_x
            \bigg]^2
      \chi
       dx ds
      \leq
      \mathcal{O}(1)
      \delta_1^2,
\end{equation}
where $\chi$ is the characteristic function of the set
\begin{equation*}
        \bigg\{
               x: \, \bigg| \frac{w_1}{v_1}(x)
                            - \lambda_1^{\ast}
                     \bigg|
              \leq
               3 \hat{\delta}        \bigg\}
\end{equation*}
(see Section \ref{par_gradient_decomposition} for the definition
of $\hat{\delta}$).

 We preliminary observe that the following
equalities hold:
\begin{equation*}
\begin{split}
&       |v_1|   \bigg[ \bigg(
                       \frac{w_1}{v_1}
                       \bigg)_x
        \bigg]^2 =
        \frac{w_{1x}^2 v_1^2 +
        v_{1x}^2w_1^2-
        2v_{1x}w_{1x}v_1w_1}
        {|v_1^3|} \leq
        C
        \frac{ | \gamma_{xx}|^2 |\gamma_x|^2-
        \langle \gamma_x, \, \gamma_{xx} \rangle ^2}
       {|\gamma_x|^3}, \\
&      |\lambda_1 \gamma_x|_x =
        \frac{\langle  \lambda_1 \gamma_x, \,
              (\lambda_1 \gamma_x)_x \rangle }
             {| \lambda_1 \gamma_x|}=
        -  \frac{\langle   \gamma_x, \,
                 (\lambda_1 \gamma_x)_x \rangle }
                {| \gamma_x|}, \\
&     | \gamma_{x}|_{xx} =
      \bigg(
           \frac{\langle  \gamma_x, \, \gamma_{xx} \rangle }
                {|\gamma_x|}
      \bigg)_x
      =
      \frac{\langle  \gamma_x, \, \gamma_{xx} \rangle ^2}
           {-|\gamma_x|^3}
      +
      \frac{\langle  \gamma_x, \, \gamma_{xxx} \rangle }
                {|\gamma_x|}
      +
      \frac{|\gamma_{xx}|^2}{|\gamma_x|}. \\
\end{split}
\end{equation*}
From $\gamma_{xt} + (\lambda_1 \gamma_x)_x = \gamma_{xxx}$, one gets integrating by parts
\begin{equation*}
      \begin{split}
      \frac{d \mathcal{L}}{ds} & =
      \int_0^L
              \frac{\langle \gamma_{xxx}, \, \gamma_x \rangle }
                   {|\gamma_x|} -
      \int_0^L
              \frac{\langle (\lambda_1 \gamma_x)_x, \,
                      \gamma_x \rangle }
                   { |\gamma_x|} \\
    & =
      \int_0^L |\gamma_x|_{xx} +
      \int_0^L
              \frac{\langle  \gamma_x, \, \gamma_{xx} \rangle ^2}
                   {|\gamma_x|^3} -
      \int_0^L \frac{|\gamma_{xx}|^2}{|\gamma_x|} +
      \int_0^L | \lambda_1 \gamma_x |_x. \\
      \end{split}
\end{equation*}
Hence,
\begin{equation*}
       \begin{split}
       \frac{1}{C}
       \int_0^t \int _0^L
                     |v_1| \bigg[ \bigg(
                              \frac{w_1}{v_1}
                             \bigg)_x
                          \bigg]^2
      \chi
       dx ds    \leq&~
      \int_0^t \int_0^L
     \frac{ | \gamma_{xx}|^2 |\gamma_x|^2-
            \langle \gamma_x, \, \gamma_{xx} \rangle ^2}
            {|\gamma_x|^3} dx ds
     \\
    \leq&~  - \int_0^t \frac{d \mathcal{L}}{ds} ds +
     \int_0^t
     \bigg[
           | \gamma_x |_x (s, \, x)
     \bigg]^{x=L}_{x=0}+
     \int_0^T
     \bigg[
           | \lambda_1 \gamma_x | (s, \, x)
     \bigg]^{x=L}_{x=0} ds \leq
     \mathcal{O}(1)
     \delta_1.
     \end{split}
\end{equation*}
In the previous estimate we have used the fact that $v_1, \; w_1,
\; v_{1 x}, \; w_{1 x}$ are integrable with respect  to time and
that their integrals are bounded by $\unpo \delta_1$ (Propositions
\ref{functional_estimates_pro} and \ref{other_wrt_time_pro}).
Since $
 \| v_1 \|_{\infty}$ is bounded by $\mathcal{O}(1)
 \delta_1$, the previous estimate complete the proof of
\eqref{length_functional_eq}.

\subsection{Estimate on the error in choosing the speed}
\label{energy_estimates}

The final estimate is the source term due to the cutoff function
$\theta$. Also this computation is similar to the one performed in
\cite{BiaBrevv}, taking into account the fact that here we have a
double boundary. In Appendix \ref{energy_estimates_proof} one can
find the proof of the estimates
\begin{equation}
\label{energy_estimates_eq}
       \int_{0}^{t} \int_{0}^{L}
       \Big( |v_1|+|w_1|+|v_{1 x}| \Big)\Big(|w_1+\sigma_1 v_1|\Big)(s,x)
       dx ds
       \leq
       \mathcal{O}(1) \delta_1 ^2.
\end{equation}
This ends the proof of the estimate
\begin{equation*}
      \int_0^t \int_0^L |\tilde{s}_i (s, \, x)| ds dx \leq
      \unpo \delta_1^2 \quad i=1, \, 2,
\end{equation*}
and hence of Theorem \ref{BV}.

\section{Stability estimates}
\label{stability_estimates} In this section we prove the second
part of Theorem \ref{main_result}, completing the proof. Since the
ideas are essentially the same as in the boundary free case, we
will only sketch the line of the proof, paying more attention to
the choice of the boundary conditions (which is the new element in
this paper). The result of this section is thus:

\begin{teo}
\label{teo_stability}
      There exist constants $L_1$, $L_2$ s.t. the following holds:
      let $u^1, \; u^2$ be two solutions of the parabolic system
      \begin{equation}
      \label{eq_parabolic}
        u_t + A(u)u_x - u_{xx} = 0,
      \end{equation}
      with initial and boundary data $u^1_0, \; u^1_{b 0}, \; u^1_{b L}$
      and $u^2_0, \; u^2_{b 0}, \; u^2_{b L}$ respectively.
      Then
      \begin{equation}
      \label{estimate_stability2}
      \begin{split}
               \| u^1(t ) - u^2(s) \|_{L^1(0, \, L) }
      &        \leq
               L_1 \Big( \|u^1_0 - u^2_0 \|_{L^1(0, \, L)} +
               \| u^1_{b 0} - u^2_{b 0}\|_{L^1(0, \, + \infty)} +
               \| u^1_{b 0} - u^2_{b L}\|_{L^1(0, \, + \infty)}
               \Big) \\
      &        + L_2 \Big( |t - s| + |\sqrt{t} - \sqrt{s}\, | \Big).
      \end{split}
      \end{equation}
\end{teo}

%
\subsection{Stability with respect to initial and boundary data }
We will prove that, in the hypothesis of Theorem
\ref{teo_stability},
\begin{equation}
      \label{eq_stability_bid}
              \| u^1(t ) - u^2(t) \|_{L^1(0, \, L) }
               \leq
               L_1 \Big( \|u^1_0 - u^2_0 \|_{L^1(0, \, L)} +
               \| u^1_{b 0} - u^2_{b 0}\|_{L^1(0, \, + \infty)} +
               \| u^1_{b 0} - u^2_{b L}\|_{L^1(0, \, + \infty)}
               \Big)
      \end{equation}
Let $z(t, \, x)$ be a first order perturbation of a solution $u(t,
\, x)$ of \eqref{eq_parabolic}. By straightforward computations
one gets that $z$ satisfies
\begin{equation}
\label{equation_z}
      z_t + \big( A(u) z \big)_x - z_{xx} =
      \big( DA(u) u_x \big)z -
      \big( DA(u) z \big) u_x.
\end{equation}
To prove Theorem \ref{teo_stability}, it is
enough to prove that any first order perturbation $z(t, \, x)$
satisfies the bound
\begin{equation}
\label{lipII}
      \| z(t) \|_{L^1(0, \, L )} \leq
      L_1 \Big( \| z( t= 0)\|_{L^1(0, \, L)} +
          \| z( x= 0) \|_{L^1(0, \, + \infty)} +
          \| z( x= L) \|_{L^1(0, \, + \infty)} \Big).
\end{equation}
Indeed, provided \eqref{lipII} holds, a homotopy argument which
can be found in \cite{Bre:con,BiaBre:BV} gives then the Lipschitz
estimate \eqref{eq_stability_bid}.

To prove \eqref{lipII} it is convenient to introduce the auxiliary
variable
\begin{equation*}
      \Upsilon = z_x - A(u) z,
\end{equation*}
which satisfies the equation
\begin{equation}
\label{eq_uu}
\begin{split}
      \Upsilon_t + ( A(u) \Upsilon)_x - \Upsilon_{xx} =
&     \bigg[ DA(u) ( u_x \otimes z - z \otimes u_x)
      \bigg]_x -
      A(u) \bigg[ DA(u) \big( u_x \otimes z - z \otimes u_x
                        \big) \bigg] \\
&     +
      DA(u) \big( u_x \otimes \Upsilon \big) -
      DA(u) \big( u_t \otimes z \big). \\
\end{split}
\end{equation}


%
Let $z_0 (x), \; z_{b \, 0}(t)$ and $z_{b \,  L}(t)$  be the
initial and boundary conditions we impose on $z$: since the final
goal is to apply \eqref{lipII} in the homotopy argument, it is not
restrictive to suppose that $z_0 (x), \; z_{b \,  0}(t)$ and $z_{b
\, L}(t)$ satisfy the same regularity hypothesis as $u$. Indeed,
the solution $z$ of \eqref{equation_z} that is used in the
homotopy argument is on the boundaries and at $t =0$ just the
difference of the solutions $u^1$ and $u^2$ of
\eqref{eq_parabolic}.

Hence we will suppose that  $z_0 (x), \; z_{b \, 0}(t)$ and $z_{b
\, L}(t)$ are regular and that $d^k z_0 / d x^k, \; d^k z_{b \,
0}/dt^k$ and  $d^k z_{b \, L}/dt^k, \; k = 1, \dots n$ are
integrable and have a small $L^1$ norm. Moreover, if $\| u_0^1
-u_0^2 \|_{L^1(0, \, L)}$, $\|u_{b \, 0}^1 - u_{b \, 0}^2
\|_{L^1(0, \, + \infty)}$ or $\|u^1_{b \, L} - u^2_{b \, L}
\|_{L^1(0, \, + \infty)}$ are infinite, then
\eqref{eq_stability_bid} holds trivially, and therefore we can
suppose that $z_0 \in L^1(0, \, L)$, $z_{b \, 0}, \; z_{b \, L}
\in L^1(0, \, + \infty)$.

From the hypothesis on $z_0$ it immediately follows that $\uu(t
=0)$ is regular and small in $L^1$ and sup norm. 

As in the proof of the $BV$ bounds on the solution $u$, the
crucial step to show \eqref{lipII} is the introduction of a
suitable decomposition along travelling waves and double boundary
layers: note, moreover, that $u_x$ satisfies equation
\eqref{equation_z}. Hence, it seems promising to decompose $z$
along the same vectors $\tilde{r}_i (u, \, v_i, \, \sigma_i)$ and
$\hat{r}_i(u, \, p_i)$ that appear in the decomposition
\eqref{decomposition} of $u_x$. This choice actually leads to non
integrable source terms. We will therefore allow the vectors
employed in the decomposition of $z$ to depend not only on the
solution $u$, but also on the perturbation $z$ itself:
\begin{equation*}
\left\{
\begin{array}{ll}
      z = z_1 \tr(u, \, v_1, \, \tau_1)
       + z_2 r_2 + q_1 \hr (u, \, p_1) + q_2 r_2 \\
      \Upsilon = \iota_1 \tr( u, \, v_1, \, \tau_1)
       + \iota_2 r_2. \\
\end{array}
\right.
\end{equation*}
In the previous expression the speed of the travelling waves
described by the vector $\tr$ is not $\sigma_1$, but
\begin{equation*}
  \tau_1 = \theta \bigg(
                        \lambda_1^{\ast} - \frac{z_1 }{
                        \upsilon_1}
                  \bigg) -
  \lambda_1^{\ast}.
\end{equation*}
The function $\theta$ is the cutoff
\begin{equation*}
      \theta(s) =
      \left\{
      \begin{array}{lll}
            s \quad \quad \textrm{if} \; |s|\leq \hat{\delta}     \\
            0 \quad \quad \textrm{if} \; |s|\geq 3 \hat{\delta}   \\
            \textrm{smooth connection if}
                     \quad \hat{\delta} \leq s \leq 3 \hat{\delta}
      \end{array}
      \right.
      \qquad \qquad \hat{\delta} \leq \frac{1}{3}.
\end{equation*}
The proof of \eqref{lipII} is from now on very similar to that of
the $BV$ bounds: one inserts the previous decomposition in the
equations \eqref{equation_z} and \eqref{eq_uu} and obtains the
equations:
\begin{equation}
\label{eq_z_i_q_i}
\begin{array}{llllll}
      z_{1 t}+ ( \lambda_1 z_1)_x - z_{1 xx} = 0
      &
      z_{2 t}+ ( \lambda_2 z_2)_x - z_{2 xx} = \underline{s}_1 (t, \, x) \\
      q_{1 t}+ ( \lambda_1 q_1)_x - q_{1 xx} = 0
      &
      q_{2 t}+ ( \hat{\lambda}_2 q_2)_x - q_{2 xx} = 0 \\
      \iota_{1 t}+ ( \lambda_1 \iota_1)_x - \iota_{1 xx} =
      \underline{s}_{\; 3}(t, \, x)
      &
      \iota_{2 t}+ ( \lambda_2 \iota_2)_x - \iota_{2 xx} =
      \underline{s}_{\, 2}(t, \, x) \\
\end{array}
\end{equation}
As in the proof of the $BV$ bounds, to prove \eqref{lipII} it is
sufficient to show that the condition
\begin{equation*}
      \|z(s)\|_{L^1(0, \, L)} \leq C \delta_1 \quad \forall \, s
      \in [0, \, t]
\end{equation*}
implies
\begin{equation*}
       \int_0^t \int_0^L |\underline{s}_{\; i} (s, \, x)| dx ds
            \leq \unpo \delta_1^2 \quad i = 1, \, 2, \, 3
\end{equation*}
and suitable bounds on the boundary terms. Moreover, in the proof
of the previous implication it is not restrictive to assume
\begin{equation*}
      \int_0^t \int_0^L |\underline{s}_{\; i} (s, \, x)| dx ds
            \leq \unpo \delta_1 \quad i = 1, \, 2, \, 3 ,
\end{equation*}
because a posteriori one finds a bound of order $\delta_1^2$.

Actually, one could observe that while the equations for $u_x$ and
$u_t$ have no source term (see Appendix \ref{explicit_source_t}
for details), the equations \eqref{equation_z} and \eqref{eq_uu}
have nontrivial source terms. However, one can show that both the
source terms in \eqref{equation_z} and \eqref{eq_uu} and the other
terms that contribute to $\underline{s}_i$, $i=1, \, 2, \, 3$ can
be bounded by an expression analogous to the one that appears on
the right side of \eqref{reasons_of_source_term_eq}. The
computations that ensure such an estimate are quite similar to
those performed in the proof of Section
\ref{reasons_of_source_term}.

The proof of \eqref{lipII} can therefore be completed with the
same tools described in Paragraph \ref{BV_estimates}, hence we
will skip all the details.

\subsection{Stability with respect to time}

Let $u(t, \, x)$ be a solution of \eqref{eq_parabolic}: from
Proposition \ref{pro_u_xx} and the observations that follow one
gets
\begin{equation*}
       \|u_{xx}(t)\|_{L^1} \leq
       \left\{
       \begin{array}{ll}
            \unpo \delta_1 / \sqrt{t} \quad t \leq 1 \\
            \unpo \delta_1 \quad t > 1. \\
       \end{array}
       \right.
\end{equation*}
Let $t_1 \leq t_2$: the estimate above implies
\begin{equation}
\label{estimate_L1_wrtt}
\begin{split}
      \|u(t_1 ) - u(t_2)\|_{L^1(0, \, L)}
&     \leq
      \int_{t_1}^{t_2} \bigg\|
      \frac{\partial u}{ \partial t}
      (t, \, x) \bigg\|_{L^1} dt      \leq
      \int_{t_1}^{t_2} (\unpo \| u_x (t, \, x)\|_{L^1} +
      \|u_{xx}(t, \, x))\|_{L^1} ) dt \\
&     \leq
      \unpo \int_{t_1}^{t_2} ( \delta_1 + \delta_1 / \sqrt{t}) dt
      \leq \unpo \delta_1 |t_1 - t_2| + \unpo \delta_1 |\sqrt{t_1} -
      \sqrt{t_2} \, | \\
&     \leq L_2 \Big( |t_1 - t_2|+ |\sqrt{t_1} - \sqrt{t_2}\, | \Big). \\
\end{split}
\end{equation}
This completes the proof of Theorem \ref{teo_stability} and hence
of Theorem \ref{main_result}.
%
%

\section{The vanishing viscosity limit}
\label{par_semigroup}

In this section we prove Theorem \ref{T:2}. The proof proceeds in
two steps: first, by using the results of Theorem
\ref{main_result}, we obtain that there exists a subsequence of
solutions $u^\epsilon$ to the problem
\begin{equation*}
 \left\{
      \begin{array}{lllll}
            u_t + A (u) u_x =0 ,
            \quad
            x \in \, ]0, \, l[ \; \;
            t \in \, ]0, + \infty [ \\
            \\
            u(0, x) = \bar{u}_0 (x)\\
            \\
            u(t, 0) = \bar{u}_{b \, 0}(t) \qquad
            u(t, l) = \bar{u}_{b l}(t) \\
     \end{array}
     \right.
\end{equation*}
which converges to a Lipschitz semigroup. Then we use the
machinery of viscosity solutions to complete the proof, showing
the uniqueness of the limit. In particular, we exhibit explicitly
the boundary Riemann solver.


Let $\semie$ the solution of the system \eqref{vvapproximation}:
from Theorem \ref{BV} one gets that the total variation of the
solution of system \eqref{rescaled} is uniformly bounded with
respect to time and hence, by a change of variables, $\semie$
satisfies
\begin{equation*}
       \bv \big\{ \semie \big\}
       , \; \Big| \semie( x) \Big| \leq \unpo \delta_1
       \quad
       \forall \, t>0, \; x \in \, [0, \, l], \;
       \ee > 0
\end{equation*}
and for any $\bar{u}_0 \in \domainu \; \bar{u}_{b \, 0}, \,
\bar{u}_{b\, l} \in \domainub.$ By Helly's theorem, for every
sequence $\ee_n \to 0^+$ and for any $t \ge 0 $ there exists a
subsequence, which we still call $\ee_n$ for simplicity, such that
$p_t^{\; \ee_n}[\bar{u}_0, \, \bar{u}_{b \, 0}, \, \bar{u}_{b \,
l}]$ converges in $L^1(0, \, l)$. The stability with respect to
time and to initial and boundary data ensures that, by a standard
diagonalization procedure, one can find a function
\begin{equation*}
\begin{array}{ccccc}
       p &:&
     [0, \, + \infty [ \,  \times \, \domainu
       \times \domainub \times \domainub
       & \to & \domain \\
& &     (t, \, \bar{u}_0, \, \bar{u}_{b \, 0}, \, \bar{u}_{b \, l})
       & \mapsto &
       \semi \\
\end{array}
\end{equation*}
such that, up to subsequences,
\begin{equation*}
      p_t^{\; \ee_n}(t)[\bar{u}_0, \, \bar{u}_{b \, 0}, \,
      \bar{u}_{b \, l}] \to \semi \quad L^1(0, \, l)\quad
      \forall \, t \ge 0, \;
      \bar{u}_0 \in \domainu, \;  \bar{u}_{b \, 0},
        \, \bar{u}_{b \, l} \in \domainub.
\end{equation*}
Moreover, one can verify that the function
\begin{equation}
\label{eq_semigroup}
\begin{array}{ccccc}
       S &:&
      [0, \, + \infty [ \,  \times \, \domainu
       \times \domainub \times \domainub
       & \to& \domain \times \domainub \times \domainub \\
& &     (t, \, \bar{u}_0, \, \bar{u}_{b \, 0}, \, \bar{u}_{b \, l})
       & \mapsto&
       \bigg( \semi, \, \bar{u}_0( \, \cdot \, + t), \,
       \bar{u}_{b \, 0}( \, \cdot \, + t), \,
       \bar{u}_{b \, l}(\, \cdot \, +t ) \bigg)\\
\end{array}
\end{equation}
satisfies the semigroup properties, together with the Lipschitz estimate
\begin{equation}
\label{E:lip2}
      \begin{split}
             \Bigl\| \semi  - p_s [ \bar{v}_{0}, \, \bar{v}_{b \,  0}, \,
             \bar{v}_{b \, l} ] \Bigr\|_{L^1} \leq
      &      L_1 \bigg(  \|\bar{v}_{ 0} - \bar{u}_{0}\|_{L^1(0, \, l)}+
            \|\bar{v}_{b 0} - \bar{u}_{b 0}\|_{L^1(0, \, + \infty)}
            \\
      &     +  \|\bar{v}_{b l} - \bar{u}_{b l}\|_{L^1(0, \, + \infty)} \bigg)+
             L_2 |t -s|,        \\
      \end{split}
\end{equation}

We now make use of the tool of viscosity solution, which was first
introduced in \cite{Bre:Gli}.

\subsection{The Riemann solver and the boundary Riemann solver}
\label{par_riemann} 
A crucial step in the proof of the uniqueness of the vanishing
viscosity limit is the local description of the vanishing
viscosity solution in case of piecewise constant data, which
however has an interest in its own. The aim of this section is to
characterize the limit as $\ee_n \to 0^+$ of the solution of
\begin{equation}
\label{eq_briemann}
\left\{
\begin{array}{lll}
     u_t + A(u) u_x = \ee_n u_{xx} \\
     u(0, \, x) =
      \left\{
      \begin{array}{ll}
      u^+ \quad x > 0 \\
      u^- \quad x < 0 \\
      \end{array}
      \right. \\
     u(t, \, 0) \equiv  u_{b \, 0}
      \qquad
      u(t, \, l) \equiv u_{b \, l} \\
\end{array}
\right.
\end{equation}
where $u^+, \; u^-,  \; u_{b \, 0}$ and $u_{b \, l}$ are
constants. In the following, we will write "solution to the
Riemann problem" meaning "vanishing viscosity solution to the
Riemann problem".

In \cite{AnBia,BiaBrevv} it is shown that the solution of
\eqref{eq_briemann} is defined locally: to solve
\eqref{eq_briemann} it is therefore sufficient to characterize the
vanishing viscous solutions in the following three cases:
\begin{enumerate}
\item
     the Cauchy problem with datum
     \begin{equation*}
     u_0(x) =
     \left\{
     \begin{array}{ll}
           u^- \quad x < 0 \\
           u^+ \quad x > 0 \\
     \end{array}
     \right.
     \end{equation*}
\item the boundary problem at $x=0$
     \begin{equation*}
     \left\{
     \begin{array}{ll}
           u(0, \, x) \equiv u_{0} \\
           u(t, \, 0) \equiv u_{b \, 0} \\
     \end{array}
     \right.
     \end{equation*}
\item the boundary problem at $x=l$
     \begin{equation*}
     \left\{
     \begin{array}{ll}
           u(0, \, x) \equiv u_{0} \\
           u(t, \, l) \equiv u_{b \, l} \\
     \end{array}
     \right.
     \end{equation*}
\end{enumerate}
The second and the third case are clearly analogous, and therefore
 in Section
\ref{par_boundary_solver} we will deal only with the second one.
In the following section, instead, we will recall for completeness
the essential steps of the construction of the solution in case 1:
we refer to \cite{BiaBrevv} for an exhaustive account.

In any of
the three cases the crucial step is the definition of two families
of admissible states, as it will be
clearer in the following.

\subsection{The non conservative Riemann solver}
\label{par_non_cons}

Since in this case the construction of the
first and the second curve of admissible states is the same, we
will describe only the construction of the first curve $T^1u_r$ of
the
 states that can be connected
by waves of the first family to a right state $u_r$.
For a general reference, see \cite{BiaBrevv}.

Consider the family $\Upsilon \subset \mathcal{C}^0( \, [0, \, s
]; \mathbb{R}^n \times \mathbb{R} \times \mathbb{R})$ of curves
\begin{equation*}
      \tau \mapsto (u(\tau), \, v_1 (\tau), \, \sigma_1(\tau)),
      \quad \tau \, \in \, [0, \, s]
\end{equation*}
with
\begin{equation*}
      | u (\tau) - u^{\ast}| \leq \ee,
       \quad
       |v_1| \leq \ee,
       \quad
       |\lambda_1^{\ast}- \sigma_1 (\tau) |
       \leq \ee.
\end{equation*}
The function $f_1(\tau)$ related to the curve $\gamma \in \uu$ is
defined as
\begin{equation*}
      f_1 (\tau) = \int_0^{\tau} \lambda_1 \big( u(\varsigma) \big)
      d \varsigma.
\end{equation*}
Let $\tr$ be the generalized eigenvector of the travelling waves
of first family (see Section \ref{par_travelling_waves} for the
proper definition of $\tr$). By the contraction map principle, one
can show that if $s$ is small enough then for any $\tau \in [0, \,
s]$ there is a solution $(\hat{u}, \, \hat{v}_1, \,
\hat{\sigma}_1)$ of the following system:
\begin{equation*}
\left\{
\begin{array}{lll}
     \hat{u}( \tau) = {\displaystyle u_r + \int_0^{\tau} \tr
      \big( \hat{u}(\varsigma), \, \hat{v}_1 (\varsigma), \,
      \hat{\sigma}_1 (\varsigma) d \varsigma } \\
     \hat{v}_1 (\tau) = \mathrm{conc}_{[0, \, s]} f_1 (\tau) -
      f_1 (\tau) \phantom{{\displaystyle \int}} \\
     \hat{\sigma}_1 (\tau)=
     {\displaystyle \frac{ d \mathrm{conc}_{[0, \, s]} f_1}{d \tau}}. \\
\end{array}
\right.
\end{equation*}
We indicate with $\mathrm{conc}_{[0, \, s]} f_1$ the
concave
envelope of $f_1$ in the interval $[0, \, s]$.

The curve of admissible states passing through $u_r$ is defined as
$T^1_{s} u_r = \hat{u}(s)$. Indeed, let
\begin{equation*}
       \tilde{u}( x/t) =
       \left\{
       \begin{array}{lll}
            T^1_s u_r
             \quad x / t < \sigma_1 (s)\\
            \hat{u}(\tau)
             \quad \sigma_1 (\tau) = x / t \\
            u_r
             \quad x /t > \sigma_1 (0):
       \end{array}
       \right.
\end{equation*}
one can show that any sequence of vanishing viscosity solution of
the Riemann problem with data $(u_r, \, T^1_s u_r )$ converges to
$\tilde{u}$. Moreover, the curve $T^1_s u_r$ is Lipchitz
continuous.

\subsection{The boundary Riemann solver}
\label{par_boundary_solver}

In this paragraph we will construct
the vanishing viscosity solution in case 2. We will proceed as
follows: we will construct two curves of admissible states $Z^1$
and $Z^2$ and given a right state $u_0$ and a left state $u_{b \,
0}$, we will show that there is a couple $(s_1, \, s_2)$ such that
\begin{equation*}
      Z^1_{s_1}  \circ Z^2_{s_1} u_0 = u_{b \, 0}.
\end{equation*}
The waves of the second family are entering the domain: it is
therefore quite reasonable to suppose that they are not influenced
by the presence of the boundary and therefore the second
admissible curve will be the one defined in the previous
paragraph, $Z^2_s u_0= T^2_s u_0$. Let $\bar{u} = Z^2_{s_2} u_0$ be the value reached
throughout the waves of the first family.

\begin{figure}
\caption{the graphic and the orbit of a boundary layer of the
first family connecting $u_{b \, 0}$ to $\bar{u}$: when $\ee \to
0^+$ the graphic is pressed against the axis $x=0$}
\label{fig_boundary_layers}
\begin{center}
\psfrag{a}{$\bar{u}$} \psfrag{b}{$u_{b \, 0}$}
\includegraphics[scale=0.4]{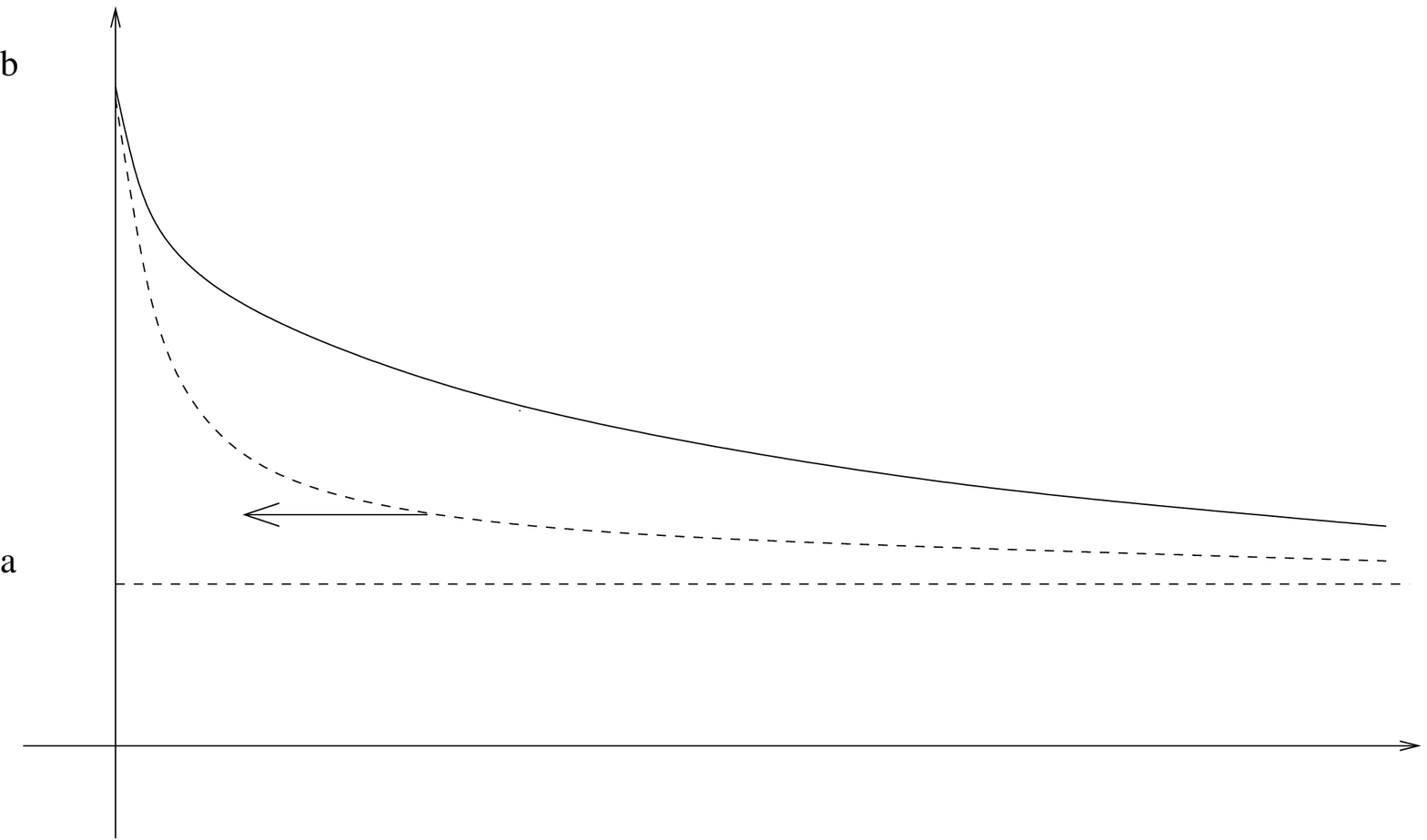}
\hfill
\psfrag{W}{$p_2$} \psfrag{V}{$p_1$} \psfrag{Z}{$u$}
\psfrag{v}{$u_{b \, 0}$} \psfrag{u}{$\bar{u}$}
\psfrag{o}{$(u^{\ast}, \, 0, \, 0)$}
\includegraphics[scale=0.4]{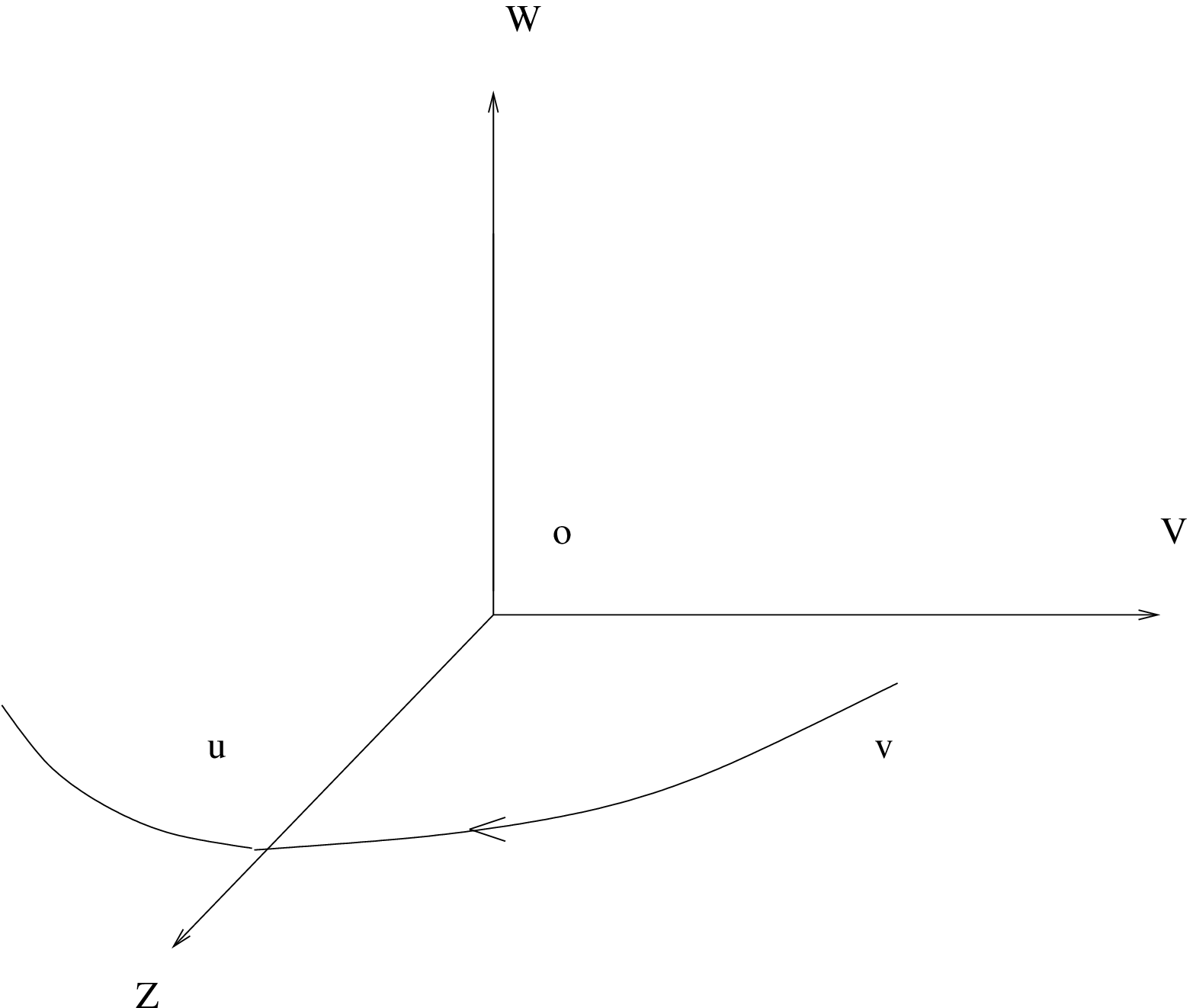}
\end{center}
\end{figure}

The waves of the first family are leaving the domain and are
therefore affected by the boundary datum. To understand their
behavior, it is convenient to focus the attention on the boundary
layers of the first family, i.e. on the solution of
\begin{equation}
\label{eq_stat}
\begin{split}
       u_{xx}= A(u)u_x
\end{split}
\end{equation}
that are exponentially decreasing to an equilibrium as $x \to
+\infty$. One can now go back to the problem
\begin{equation*}
       A(u^{\ee}) u^{\ee}_x = \ee u^{\ee}_{xx}
\end{equation*}
and let $\ee \to 0^+$. Since $u^{\ee}( x) = u(x/\ee)$, we get
\begin{equation}
\label{eq_boundary_layers}
      \lim_{\ee \to 0^+} u^{\ee}( 0^+) =
      \lim_{x \to + \infty} u( x).
\end{equation}
Such a behavior is illustrated in figure
\ref{fig_boundary_layers}.

The value $\lim_{\ee \to 0^+} u^{\ee}(0^+)$ is the state reached
throughout the waves of the second family: we called it $\bar{u}$.
It also represents the trace of the hyperbolic limit on the
boundary $x=0$. From \eqref{eq_boundary_layers} it follows that
the states which can be connected to $\bar{u}$ by boundary layers
are the initial points of orbits that decrease exponentially to
$\bar{u}$, i.e. that lay on the stable manifold throughout
$\bar{u}$.

The stable manifold at the equilibrium point $(\bar{u}, \, 0)$ of
the system
\begin{equation}
\label{eq_blayer}
\left\{
\begin{array}{ll}
     u_x = p \\
     p_x = A(u) p \\
\end{array}
\right.
\end{equation}
is parameterized by the projection $p_1$ of $p$ on the stable
space. Passages analogous to those in Section \ref{par_double_bp}
ensure that the stable manifold is characterized by the relation
$p= p_1 \breve{r}_1 (p_1)$ for a suitable vector function
\begin{equation*}
       \breve{r}_1=
       \left(
       \begin{array}{cc}
            1 \\
            f(p_1)
       \end{array}
       \right).
\end{equation*}
One imposes $u_1 (+ \infty) = \langle l_1,  \, \bar{u} \rangle $
and from the second equation gets
\begin{equation*}
      u_1 (0) =  \langle l_1, \, \bar{u} \rangle  -
      p_{1}(0) \exp \bigg( \int_0^{+ \infty}
      \lambda_1 \Big( u_1 \big( p_1 (0), \, x \big)\Big)d x \bigg).
\end{equation*}
Since $\lambda_1 \leq -c < 0$, the previous map is invertible and
one can express $p_1(0)$ as a function of $u_1(0).$ The inverse
map is clearly regular.

We parameterize the stable manifold by $ s_1 : = u_1 - <l_1 , \,
\bar{u}>$ and obtain (for some suitable regular function $z$) the
map
\begin{equation}
\label{eq_stable_man}
      Z^1_{s_1} \bar{u} =
      \left(
      \begin{array}{cc}
             \langle l_1 , \, \bar{u} \rangle  + s_1 \\
            z(s_1)
      \end{array}
      \right),
\end{equation}
defined on a small enough interval $[0, \, s]$.
%

The vanishing viscosity solution of
\begin{equation}
\label{eq_riemann}
\left\{
\begin{array}{lll}
     u_t + A(u) u_x =0 \\
     \\
     u(t, \, 0) \equiv u_0
      \qquad
      u(0, \, x) \equiv u_{b \, 0}
\end{array}
\right.
\end{equation}
can be constructed patching together the curve described so far.
Let
\begin{equation*}
      u_{b \, 0} =  Z_{s_1}^1 \circ T_{s_2}^2 u_0:
\end{equation*}
thanks to a version of the implicit function theorem valid for
Lipschitz maps (see \cite{Cl}), one can reconstruct from $u_0$ and
$u_{b \, 0}$ the couple $(s_1, \, s_2)$. The vanishing viscosity
solution of \eqref{eq_riemann} is then given by
\begin{equation*}
       u(t, \, x) =
       \left\{
       \begin{array}{lll}
            T^2_{s_2} u_0
             & x / t < \sigma_2
              (s_2)\\
            \hat{u}(\tau)
             & \sigma_2 (\tau) = x / t \\
            u_0
             & x /t > \sigma_2 (0).
       \end{array}
       \right.
\end{equation*}
One gets in particular that the trace of the solution at $x =0$ is
not necessarily the boundary value $u_{b \, 0}$, but it is the
intermediate state $T^2_{s_2} u_0$.

\begin{rem}
\label{rem_comparison} In the case of systems in conservation
form, with only linearly degenerate or genuinely non linear
fields, a boundary Riemman solver was introduced in
\cite{DubLeFl}. In that paper, it was introduced the following
admissibility condition on the trace $u(t, \, 0^+) = \bar{u}$ of
the solution of \eqref{eq_riemann}: the solution in the sense of
Lax \cite{Lax} of the Riemann problem
\begin{equation*}
\left\{
\begin{array}{ll}
       u_t + f(u)_x =0 \\
       u(0, \, x)=
       \left\{
       \begin{array}{lll}
              \bar{u}_b \qquad x < 0 \\
              \\
              \bar{u} \; \qquad x > 0 \\
       \end{array}
       \right.
\end{array}
\right.
\end{equation*}
is composed only of waves with non positive speed. Such a
condition is in general different from \eqref{eq_stable_man} and
therefore the two boundary Riemann solvers do not coincide.

On the other side, in \cite{Good, SabTou:mixte, Ama, AmaCol} it
was considered a quite general boundary condition: more precisely,
let $N$ be the dimension of the system
\begin{equation*}
      u_t + f(u)_x =0
\end{equation*}
and let $p$ be the number of positive eigenvalues of $Df(u)$,
which is supposed to be constant. Let ${b: \mathbb{R}^N \to
\mathbb{R}^p}$ be a regular enough function such that $Db(u)$ is
injective on the space generated by the $p$ eigenvectors of
$Df(u)$ associated to positive eigenvalues; then, given $g: [0, \,
+ \infty[ \to \mathbb{R}^p$, the boundary condition considered in
\cite{Good, SabTou:mixte, Ama, AmaCol} is $g(t) = b \big(u(t, \,
0^+) \big)$. Such a definition, which in the original papers was
introduced in the case of conservative systems with only linearly
degenerate or genuinely non linear fields, is compatible with the
boundary Riemann solver defined by the vanishing viscosity limit.
Indeed, in our case $N=2$, $p=1$: let $b(u)$ be equal to the
coordinate of $u$ along the curve of admissible states $T^2 u_0$,
i.e. let $b(u) = s_2$ if $u= Z^1_{s_1} \circ T^2_{s_2} u_0$.
Moreover, let $g$ be the coordinate of $\bar{u}_b$ along the same
curve: with this choice, the condition
\begin{equation*}
      u(t, \, 0^+) =
      T^2_{s_2} u_0 \qquad
      \bar{u}_b = Z^1_{s_1} \circ T^2_{s_2} u_0
\end{equation*}
is equivalent to $g(t) = b\big( u(t, \, 0^+) \big)$.
\end{rem}

\subsection{Viscosity solutions}
\label{par_viscosity_solutions}

Before giving the definition of viscosity solution we have to
introduce some preliminary notation; moreover, in the following we
will use the spaces $\domainu$, $\domainub$, $\domain$ that have
been defined in the introduction (equation \eqref{eq_domain} and
previous lines).

Let $u(t, \, x)$ be a function such that, for any $t$, $u(t) \in
\mathcal{D}_0$: given a point $(\tau, \, \xi) \in ]0, \, l[ \times
[0, \, + \infty[$, let $A^{\mathfrak{b}}=A\big(u(\tau, \, \xi)
\big)$ and let $U^{\mathfrak{b}}_{(u, \, \tau, \, \xi)}$ be the
solution of the linear Cauchy problem
\begin{equation*}
       w_t + A^{\mathfrak{b}}w_x = 0
       \qquad
       w(0, \, x) = u(\tau, \, x).
\end{equation*}
 Viceversa, let
$U^{\sharp}_{(u, \, \tau, \, \xi)}$ be the solution (defined in
Section \ref{par_non_cons}) of the Riemann problem
\begin{equation*}
\begin{split}
&      u_t + A(u) u_x = 0 \\
&      u(0, \, x) =
       \left\{
       \begin{array}{ll}
            u(\tau, \, \xi^-)
             \qquad
             x < 0 \\
            u(\tau, \, \xi^+)
             \qquad
             x >0 \\
       \end{array}
       \right. \\
\end{split}
\end{equation*}
The previous limits are well defined, since $u(\tau ) \in BV(0, \,
l)$. Given a function $\bar{u}_{b \, 0} \in \domainub$, the
definition of $U^{\sharp}_{(u, \, \tau, \, \xi)}$ can be extended
naturally to the case $\xi=0$: it is enough to define
$U^{\sharp}_{(u, \, \bar{u}_{b 0}, \, \tau)}$ as the solution
(described in Section \ref{par_boundary_solver}) of the boundary
Riemann problem
\begin{equation*}
\left\{
\begin{array}{ll}
      u_t + A(u) u_x = 0 \\
      u(0, \, x) \equiv  u(\tau, \, 0^+)
       \qquad
       u(t, \, 0) \equiv \bar{u}_{b \, 0}( \tau^+). \\
\end{array}
\right.
\end{equation*}
Given a function $\bar{u}_{b \, l} \in \domainub$, the definition
of $U^{\sharp}_{(u, \, \bar{u}_{b l}, \, \tau)}$ is clearly
analogous.
\begin{say}
\label{def_vs}
      Let $u(t, \, x)$ such that for any $t$, $u(t) \in
      \mathcal{D}_0$ and such that the function
      $t \mapsto u(t, \, \cdot \, )$ is continuous
      in $L^1_{loc}$ and let
      $\bar{u}_{b \, 0}, \; \bar{u}_{b \, l} \in \domainub$ and
      $\bar{u}_0 \in \domainu$.

      Then $u$ is a viscosity solution of the system
      \begin{equation}
      \label{eq_the_problem_II}
      \left\{
      \begin{array}{lllll}
            u_t + A (u) u_x =0 ,
            \quad
            x \in \, ]0, \, l[, \;
            t \in \, ]0, + \infty [ \\
            \\
            u(0, x) = \bar{u}_0 (x)\\
            \\
            u(t, 0) = \bar{u}_{b \, 0}(t) \qquad
            u(t, l) = \bar{u}_{b \, l}(t) \\
      \end{array}
      \right.
      \end{equation}
      if and only if the followings hold: \\
      (i) $u(0) = \bar{u}_0$ \\
      (ii) for every $\beta> 0$ and for every point
      $(\tau, \, \xi)$ with $\xi \neq 0, \; l$
       \begin{equation*}
             \lim_{h \to 0^+} \frac{1}{h}
             \int_{\max \{0, \, \xi - \beta h \}}^{\min\{l, \, \xi + \beta h\}}
             \big| u( \tau+ h, \, x) -
             U_{(u, \, \tau, \, \xi )}^{\sharp}(h, \, x - \xi)
             \big| dx = 0
       \end{equation*}
       (iii) for every $\beta > 0$ and for every $\tau > 0$
        \begin{equation*}
             \lim_{h \to 0^+} \frac{1}{h}
             \int_0^{\min\{l, \, \beta h\}}
             \big| u( \tau+ h, \, x) -
             U_{(u, \, \bar{u}_{b 0}, \,  \tau )}^{\sharp}(h, \, x )
             \big| dx = 0;
        \end{equation*}
        and
        \begin{equation*}
             \lim_{h \to 0^+} \frac{1}{h}
             \int_{\max \{0, \, l - \beta h \}}^{l}
             \big| u( \tau+ h, \, x) -
             U_{(u, \, \bar{u}_{b l}, \,  \tau )}^{\sharp}(h, \, x )
             \big| dx = 0
        \end{equation*}
       (iv) there exist constants $C$ and $\beta'$
       such that for every point $(\tau, \, \xi)$
       with $\xi \neq 0, \, l$ and for every $\rho>0$ small enough
       \begin{equation*}
             \limsup_{h \to 0^+} \frac{1}{h}
             \int_{\max \{0, \, \xi- \rho  +
             \beta' h \}}^{\min \{l, \, \xi + \rho - \beta' h\}}
             \big| u( \tau+ h, \, x) -
             U_{(u, \, \tau, \, \xi)}^{\mathfrak{b}}(h, \, x - \xi)
             \big| dx \leq C \bigg( \bv \big( u(\tau), \,
             ] \xi - \rho, \, \xi + \rho[ \big) \bigg)^2.
       \end{equation*}
\end{say}

The previous definition may appear a bit complex: note, however,
that, since $\rho$ and $h$ can be arbitrarily small, it is not
restrictive to suppose
\begin{equation*}
\begin{array}{ll}
 \max \{0, \, \xi - \beta h \}=
      \xi - \beta h
& 
      \min\{l, \, \xi + \beta h\}=
      \xi + \beta h \\
     \max \{0, \, \xi - \rho + \beta' h \} =
      \xi - \rho + \beta' h
& 
      \min \{l, \, \xi + \rho - \beta' h \} =
      \xi + \rho - \beta' h \\
     \max \{ 0, \, l - h \beta \} = l - h \beta
& 
      \min \{ l, h \beta \} = h \beta
\end{array}
\end{equation*}

The definition of viscosity solution ensures, roughly speaking,
that a function is well approximated by the solution of a suitable
linear problem and of a suitable Riemann problem.

The following proposition ensures that viscosity solutions
coincide indeed with vanishing viscosity limits. The proof is very
simile to that of the analogous property stated in \cite{BiaBrevv}
(Lemma 15.2, page 308) and will be therefore omitted.
\begin{pro}
\label{pro_viscosity}
       Let $\bar{u}_0 \in \domainu $ and $\bar{u}_{b \, 0}, \;
       \bar{u}_{b \, l} \in \domainub$. Let
       $p_t (\bar{u}_0, \, \bar{u}_{b \, 0}, \, \bar{u}_{b
      \,l})$ be a vanishing viscosity solution of the system
      \eqref{eq_the_problem_II}:
      then $p_t (\bar{u}_0, \,
      \bar{u}_{b \, 0}, \, \bar{u}_{b \,l})$ is a viscosity
      solution of the same system. \\
      Viceversa, if $u(t, \, x)$ is a viscosity solution
      of the problem \eqref{eq_the_problem_II}
      then
      \begin{equation*}
             u(t) = p_t (\bar{u}_0, \,
             \bar{u}_{b \, 0}, \, \bar{u}_{b \,l})
             \quad \forall \, t \ge 0.
      \end{equation*}
\end{pro}

From the previous result it immediately follows the uniqueness of
the semigroup: indeed, let by contradiction $p^1_t (\bar{u}_0,
\,\bar{u}_{b \, 0}, \, \bar{u}_{b \,l})$ and $p^2_t (\bar{u}_0, \,
\bar{u}_{b \, 0}, \, \bar{u}_{b \,l})$ be two different vanishing
viscosity solutions. The function $p^1_t (\bar{u}_0, \,\bar{u}_{b
\, 0}, \, \bar{u}_{b \,l})$ is hence a viscosity solution of
problem \eqref{the_problem} by the first part of Proposition
\ref{pro_viscosity}. Then $p^1_t (\bar{u}_0, \,\bar{u}_{b \, 0},
\, \bar{u}_{b \,l})= p^2_t (\bar{u}_0, \,\bar{u}_{b \, 0}, \,
\bar{u}_{b \,l})$ for any $t \ge 0$ by the second part of the
proposition.

\appendix

\section{Appendix}
\subsection{Appendix to Section \ref{parabolic_estimates}}
\subsubsection{Proof of Proposition \ref{estimate_kernels_pro}}
\label{proof_pro_kernels} In the following, for simplicity we will
suppose $\lambda_i^{\ast} = \lambda^{\ast}_2 > 0$, since the case
$\lambda_i^{\ast} =
\lambda^{\ast}_1 < 0$ is analogous.\\
We denote by
\begin{equation*}
       \Gamma ^{\lambda_2^{\ast}   }  (t, x, y) = ( 1 - e^{- x y / t})
        G(t, \, x - y - \lambda_2^{\ast} t)
\end{equation*}
the solution of the equation
\begin{equation}
\label{eq_lin}
      z_t + \lambda_2^{\ast} z_x -z_{xx}=0
\end{equation}
in the first quadrant with zero boundary datum and Cauchy datum $
  \delta_y
$. The following estimates have been proved in \cite{BiaBrevv}:
\begin{equation}
\label{estimate_from_biabre}
\begin{split}
&    \|\Gamma^{\lambda_2^{\ast}}(t, \, y)\|_{L^1(0, \, + \infty)}
     \leq
      \unpo
      \qquad
      \bigg| \int_0^{+ \infty}
             \bigg|
             \int_y^{+ \infty} \Gamma_x^{\lambda_2^{\ast}}
             (t, \, x, \, \xi) d \xi
      \bigg| dx \bigg| \leq \unpo \qquad \forall \, t \in \mathbb{R}^+ \\
&      \|\Gamma_x^{\lambda_2^{\ast}}(t, \, y)\|_{L^1(0, \, +
\infty)} \leq
      \frac{\unpo}{\sqrt{t}}
      \qquad
      \bigg| \int_0^{+ \infty }
      \bigg|
             \int_y^{+ \infty} \Gamma_{xx}^{\lambda_2^{\ast}}
             (t, \,  x, \,  \xi) d \xi
      \bigg| dx  \bigg|
      \leq \frac{\unpo}{\sqrt{t}} \quad \forall \, t \in (0,1) \\
\end{split}
\end{equation}
Let $G(t, \, x) = \exp ( -x^2 / 4t ) / 2 \sqrt{ \pi t}$: we will
use the notation
\begin{equation*}
       G(t, \, x -  \lambda_2^{\ast} t) = G^{\lambda_2^{\ast}}(t, \, x ).
\end{equation*}
The estimate on the $L^1$ norm of $\Delta^{\lambda_2^{\ast}}$ in
Proposition \ref{estimate_kernels_pro} can be obtained via the
maximum principle applied to equation \eqref{eq_lin}: indeed,
\begin{equation*}
      0 \leq \Delta^{\lambda_2^{\ast}}(t, \, x, \, y) \leq G^{\lambda_2^{\ast}}(t, \, x-y) \quad \forall \, t \ge 0
      \quad x, \, y \in ]0, \, L[,
\end{equation*}
and therefore $\| \Delta^{\lambda_2^{\ast}} (t, \, y) \|_{L^1}
\leq 1$.

To prove the estimate on the $L^1$ norm of
$\Delta_x^{\lambda_2^{\ast}}$ it is convenient to write
$\Delta^{\lambda_2^\ast}$ as
\begin{equation*}
\begin{split}
         \Delta^{ \lambda_2^{\ast}} (t, \, x, \, y)  =
&        \Gamma ^{\lambda_2^{\ast}   }  (t, \,  x, \, y)
        + \phi^{\, \lambda_2^{\ast}} (t, \, x, \, y)
         \sum_{m \, \neq 0} \Bigl[
         G (t, x + 2mL -y) -
         G (t, x+ 2mL +y ) \Bigr], \\
\end{split}
\end{equation*}
with
\[
\phi^{\, \lambda_2^{\ast}} (t, \, x, \, y) = \exp \bigg(
      \frac{\lambda_2^{\ast}}{2}(x-y) -
      \frac{ ( \lambda_2^{\ast})^2 }{4}t \bigg).
\]

Since $\lambda_2^{\ast}>0$, for $m > 0$ it holds
\begin{equation*}
\begin{split}
&     \bigg| \frac{\partial }{\partial x}
      \bigg( \phi^{\, \lambda_2^{\ast}}
       (t, \, x, \, y) G(t, \, x-y + 2mL)
      \bigg)
      \bigg|  \\
&     \qquad
      =    \phi^{\, \lambda_2^{\ast}} (t, \, x, \, y)  \bigg|
          \frac{\lambda_2^{\ast}}{2}  G(t, \, x-y +
          2mL)
      +   G_x(t, \, x-y + 2mL)
      \bigg| \leq
      \big|
          G_x^{\lambda_2^{\ast}}(t, \, x-y +2mL)
      \big|,
      \phantom{\bigg(} \\
\end{split}
\end{equation*}
and similarly
\begin{equation*}
     \bigg| \frac{\partial}{\partial x}
     \bigg( \exp \bigg(
     \frac{\lambda_2^{\ast}}{2}(x-y) -
      \frac{\lambda^{\ast\,
              2}}{4}t \bigg) G(t, \, x+y + 2mL)
              \bigg)
      \bigg|  \leq
      \big|
          G_x^{\lambda_2^{\ast}}(t, \, x + y +2mL)
      \big|. \phantom{\bigg(}
\end{equation*}

The terms of the series corresponding to $m < 0$ can be estimated
as follows: let $n : = - m$ then
\begin{equation*}
\begin{split}
&     \bigg| \frac{\partial }{\partial x}
      \bigg( \phi(t, \, x,  \,  y) G (t, \, x - y -2nL) \bigg) \bigg|
      \\
& \quad     = \phi^{\, \lambda_2^{\ast}} (t, \, x, \,  y)
      \bigg|
      - G_x (t, \, 2nL - x
      +y) - \frac{\lambda_2^{\ast}}{2} G (t, \, 2nL - x
      +y)+ \lambda_2^{\ast} G (t, \, 2nL - x
      +y)
      \bigg| \\
&   \quad   \leq
    \Big| G^{\lambda_2^{\ast}}_x (t, \, 2nL -x + y )
    \Big|
      + \lambda_2^{\ast} \Big| G^{\lambda_2^{\ast}} (t, \, 2nL -x + y )
      \Big| \phantom{\bigg(} \\
\end{split}
\end{equation*}
and similarly
\begin{equation*}
      \bigg| \frac{\partial }{\partial x}
      \bigg( \phi(t, \, x - y) G_x (t, \, x + y -2nL) \bigg) \bigg|
      \leq \Big| G^{\lambda_2^{\ast}}_x (t, \, 2nL -x - y ) \Big|
      + \lambda_2^{\ast} \Big| G^{\lambda_2^{\ast}} (t, \, 2nL -x - y )
      \Big|.
\end{equation*}
Since $\|G_x^{\lambda_2^{\ast}}(t)\|_{L^1} \leq \unpo / \sqrt{t}$,
one obtains
\begin{equation*}
\begin{split}
&     \|\Delta_x^{\lambda_2^{\ast}} (t, \, y)\|_{L^1(0, \, L)}
      \leq
      \int_0^L | \Gamma_x^{\lambda_2^{\ast}}(t, \, x, \, y) | dx +
      \int_{2L}^{+ \infty} \bigg( \big|
      G_z^{\lambda_2^{\ast}}(t, \, z + y) \big| +
      \big| G_z^{\lambda_2^{\ast}}(t, \, z-y) \big| \bigg) dz \\
&     + \int_{L}^{+ \infty}
      \bigg( \big| G_z^{\lambda_2^{\ast}}(t, \, z + y) \big|
      +  \big| G_z^{\lambda_2^{\ast}}(t, \, z-y) \big| \bigg) dz +
      \lambda_2^{\ast} \int_L^{+ \infty}
      \bigg( \big| G_z^{\lambda_2^{\ast}}(t, \, z + y) \big|  +
      \big| G_z^{\lambda_2^{\ast}}(t, \, z -y ) \big|  \bigg) dz
      \leq \frac{\unpo}{ \sqrt{t}} .\\
\end{split}
\end{equation*}

In the following estimates, we will suppose $y < L/2$: by symmetry
this is not restrictive. Observe that, for $y < L/2$
\begin{equation}
\label{estimate_piccina}
       x + y - 2L  < - L / 2 < 0 \quad \forall \, x \in [0, \, L].
\end{equation}
This assumption corresponds to the fact that
the most singular part in $\Delta^{\lambda^\ast}$ is
 collected in $\Gamma^{\lambda^\ast}$, i.e. it is given by $G(t,x+y) - G(t,x-y)$.
 If $y > L/2$, then the most singular part
 would be given by $G(x-y) - G(x+y-2L)$.

One has
\begin{equation*}
\begin{split}
      \tilde{\Delta}^{\lambda_2^{\ast}}
&      (t,\, x, \, y) =
     \int_y^L \Gamma_x^{\lambda_2^{\ast}}(t, \, x, \, \xi) d \xi +
      \int_y^L \phi_x(t, \, x , \,  \xi) \sum_{m \neq 0} \Bigl[ G(t, \, x- \xi
      + 2mL) - G(t, \, x + \xi + 2mL ) \Bigr] d \xi \\
&     +
      \int_y^L \phi(t, \, x , \,  \xi) \sum_{m \neq 0} \Bigl[ G_x(t, \, x- \xi
      + 2mL) - G_x(t, \, x + \xi + 2mL ) \Bigr] d \xi \\
    =&~ \int_y^L \Gamma_x^{\lambda_2^{\ast}}(t, \, x, \, \xi) d \xi
      - \int_y^L \bigg\{ \phi_x(t, \, x , \,  \xi)
      \sum_{m \neq 0} G(t, \, x +\xi
      + 2mL) - \phi(t, \, x , \,  \xi) \sum_{m \neq 0} G_x(t, \, x+  \xi
      + 2mL) \bigg\} d \xi \\
&     + \int_y^L \bigg\{ \phi_x(t, \, x , \,  \xi)
      \sum_{m \neq 0} G(t, \, x - \xi
      + 2mL) +
      \phi(t, \, x, \,  \xi) \sum_{m \neq 0} G_x(t, \, x- \xi
      + 2mL) \bigg\} d \xi \\
     =&
      \int_y^L \Gamma_x^{\lambda_2^{\ast}}(t, \, x, \, \xi) d \xi +
      \sum_{m \neq 0}
      \phi(t, \, x, \, y) G(t, \, x + y + 2mL ) -
      \sum_{m \neq 0}
      \phi(t, \, x, \, L) G(t, \, x + L + 2mL) \\
&     -
      \int_y^L \sum_{m \neq 0} \lambda_2^{\ast}
      \phi(t, \, x , \,  \xi) G(t, \, x + \xi + 2mL) d \xi
     + \sum_{m \neq 0}
      \phi(t, \, x , \, y) G(t, \, x - y + 2mL ) \\
      &-
      \sum_{m \neq 0}
      \phi(t, \, x, \, L) G(t, \, x -L + 2mL). \\
\end{split}
\end{equation*}
The integrability of the first term follows from
\eqref{estimate_from_biabre}, the other terms are clearly
integrable because of the quadratic exponential decay of the heat
kernel $G$: hence $\|\tilde{\Delta}^{\lambda_2^{\ast}}(t,
y)\|_{L^1} \leq \unpo$.

The function $\tilde{\Delta}_x^{\lambda_2^{\ast}}$ can be written
as follows:
\begin{equation*}
\begin{split}
      \tilde{\Delta}_x^{\lambda_2^{\ast}}(t, \, x, \, y) =
      &~
      \int_y^L \Gamma_{xx}^{\lambda_2^{\ast}}(t, \, x, \, \xi) d \xi +
      \sum_{m \neq 0}
      \frac{\lambda_2^{\ast}}{2}
      \phi(t, \, x, \, y) G(t, \, x + y + 2mL ) \\
&     +
      \sum_{m \neq 0}
      \phi(t, \, x, \, y) G_x(t, \, x + y + 2mL ) -
      \sum_{m \neq 0}
      \frac{\lambda_2^{\ast}}{2}
      \phi(t, \, x, \, L) G(t, \, x + L + 2mL) \\
&      -
      \sum_{m \neq 0}
      \phi(t, \, x, \, L) G_x(t, \, x + L + 2mL) -
      \int_y^L \sum_{m \neq 0}
      \frac{(\lambda_2^{\ast })^2}{2}
      \phi(t, \, x , \,  \xi) G(t, \, x + \xi + 2mL) d \xi \\
      &-
      \int_y^L \sum_{m \neq 0}
      \lambda_2^{\ast }
      \phi(t, \, x , \,  \xi) G_x(t, \, x + \xi + 2mL) d \xi + \sum_{m \neq 0}
      \frac{\lambda_2^{\ast}}{2}
      \phi(t, \, x , \, y) G(t, \, x, - y + 2mL )\\
      & + \sum_{m \neq 0}
      \phi(t, \, x , \, y) G_x(t, \, x, - y + 2mL ) - \sum_{m \neq 0}
      \frac{\lambda_2^{\ast}}{2}
      \phi(t, \, x, \, L) G(t, \, x -L + 2mL) \\
      & -
      \sum_{m \neq 0}
      \phi(t, \, x, \, L) G_x(t, \, x -L + 2mL), \\
\end{split}
\end{equation*}
and hence with computations similar to those performed before
one gets
\begin{equation*}
      \|\tilde{\Delta}^{\lambda_2^{\ast}}_x(t, \, y)\|_{L^1} \leq
      \frac{\unpo}{\sqrt{t}}
      \quad \forall \, t \leq 1 \; \; y \in ]0, \, L[.
\end{equation*}

If one derives the explicit formula \eqref{J_0} for $J^{\lambda \,
L }$ and then integrate by parts gets
\begin{equation*}
      \int_0^L | J_x^{\lambda_2^{\ast} \, L}(t, \, x) | dx =
      \int_0^L \bigg| \lambda_2^{\ast} C e^{\lambda_2^{\ast} x} dx -
      \lambda_2^{\ast} C \int_0^L  \tilde{\Delta}^{\lambda_2^{\ast}} (t, \, x, \, y)
      e^{\lambda_2^{\ast} y} dy \bigg| \, dx \leq \unpo,
\end{equation*}
where we have used the estimate
$\|\tilde{\Delta}^{\lambda_2^{\ast}}\|_{L^1} \leq \unpo$. By
symmetry it follows $\|J_x^{\lambda_2^{\ast} \, 0}\|_{L^1} \leq
\unpo$.

Deriving $J_x^{\lambda_2^{\ast} \, L}$ one obtains
\begin{equation*}
      \|J_{xx}^{\lambda_2^{\ast} \, L}(t)\|_{L^1} \leq
      \int_0^L | C ( \lambda_2^{\ast })^2 e^{\lambda_2^{\ast} x}| dx  +
      C \lambda_2^{\ast} \int_0^L \bigg|
      \int_0^L \tilde{\Delta}_x^{\lambda_2^{\ast}}(t,
      \, x, \, y) e^{\lambda_2^{\ast} y} dy  \bigg| \, dx
      \leq \frac{\unpo}{ \sqrt{t}},
\end{equation*}
thanks to the estimate on
$\|\tilde{\Delta}_x^{\lambda_2^{\ast}}\|$. By symmetry one gets
$\|J_{xx}^{\lambda_2^{\ast} \, 0}(t)\|_{L^1} \leq \unpo /
\sqrt{t}$: this concludes the proof of Proposition
\ref{estimate_kernels_pro}.

\subsection{Appendix to Section \ref{gradient_decomposition}}

\subsubsection{Explicit source terms}
\label{explicit_source_t}

We want to find the equations satisfied by the quantities $ v_1 $,
$ v_2 $, $p_1 $, $p_2$, $w_1$, $w_2$: we will use the
decomposition
\begin{equation*}
      \left\{
       \begin{array}{ll}
       u_x = v_1 \tilde{r}_1 ({u,\, v_1, \, \sigma_1}) +
                   v_2 r_2 +
                   p_1 \hat{r}_1 (u, p_1 ) +
                   p_2 r_2 \\
       &   \\
        u_t = w_1 \tilde{r}_1 ({u,\, v_1, \, \sigma_1}) +
                   w_2 r_2
       \end{array}
       \right.
       \quad   \sigma_1 = \lambda_1(u^{\ast})
                              - \theta
                              \bigg(
                                    \frac{ w_1}{v_1} + \lambda_1(u^{\ast})
                              \bigg).
\end{equation*}
and insert it in the parabolic equation
\begin{equation*}
       u_t + A(u)u_x -u_{xx}=0.
\end{equation*}

A derivation w.r.t. $x$ gives
\begin{equation*}
\begin{array}{ll}
      \tilde{r}_{1 x} = \mathrm{D} \tilde{r}_1
                        (  v_1 \tilde{r}_1 + v_2 r_2 +
                           p_1 \hat{r}_1 +p_2 r_2      ) +
                           v_{1 x} \tilde{r}_{1 v}+
                           \sigma_{1 x} \tilde{r}_{1 \sigma} \\
      \hat{r}_{1 x}=  \mathrm{D} \hat{r}_1
                        (  v_1 \tilde{r}_1 + v_2 r_2 +
                           p_1 \hat{r}_1 +p_2 r_2      ) +
                           p_{1 x} \hat{r}_{1 p}.
      \end{array}
\end{equation*}
Recalling that
\begin{equation*}
\begin{split}
&
      \hat{\lambda}_2 :=
      \lambda_2 - p_1
      \langle  \hat{\ell}_2, \,
        \mathrm{D} \hat{r}_1 r_2  \rangle  \\
&      A(u) \tilde{r}_1 =  v_1 \mathrm{D} \tilde{r}_1 \tilde{r}_1
       + \lambda_1 \tilde{r}_1 +
       v_1 ( \lambda_1 -\sigma_1 ) \tilde{r}_{1 v}\\
&      A(u) \hat{r}_1 =  p_1 \mathrm{D} \hat{r}_1 \hat{r}_1 +
                          \lambda_1 \hat{r}_1 +
                           p_1 \lambda_1 \hat{r}_{1 p} \\
\end{split}
\end{equation*}
one gets
\begin{equation}
\label{u_t}
\begin{split}
      u_t   & =  \, u_{xx} - A(u) u_x  \\
            & = v_{1 x} \tilde{r}_1 +
                v_1 \tilde{r}_{1 x}+
                p_{1 x} \hat{r}_1 +
                p_1 \hat{r}_{1 x} +
                v_{2 x} r_2 +
                p_{ 2 x} r_2 
        - v_1 A(u) \tilde{r}_1 -
                p_1 A(u) \hat{r}_1 -
                \lambda_2 v_2 r_2 -
                \lambda_2 p_2 r_2 \\
            & = ( v_{1 x} - \lambda_1 v_1 )( \tilde{r}_1 + v_1 \tilde{r}_{1 v})+
                v_1^2 \sigma_1 \tilde{r}_{1 v}   
        + v_1 v_2 \mathrm{D}\tilde{r}_1 r_2 +
                v_1 p_1 \mathrm{D}\tilde{r}_1 \hat{r}_1 +
                v_1 p_2 \mathrm{D}\tilde{r}_1 r_2 +
                v_1 \sigma_{ 1 x} \tilde{r}_{1 \sigma} \\
            & \quad + ( p_{1 x} - \lambda_1 p_1 )( \hat{r}_1 + p_1 \hat{r}_{1 p}) +
                v_1 p_1 \mathrm{D}\hat{r}_1 \tilde{r}_1 +
                v_2 p_1 \mathrm{D}\hat{r}_1 r_2 
        + ( v_{2 x} - \lambda_2 v_2 ) r_2 +
                ( p_{2 x} - \hat{ \lambda}_2 p_2) r_2. \\
\end{split}
\end{equation}
We multiply the previous expressions by $\ell_1$ and by
$\tilde{\ell}_2$, the vectors of the dual basis of $( \tilde{r}_1,
\, r_2 )$: we obtain
\begin{equation}
\label{decomposition_of_u_t}
     \begin{array}{ll}
     w_1 = v_{1 x} - \lambda_1 v_1 +
           p_{1x} - \lambda_1 p_1 \\
     w_2 =  v_{2 x} - \lambda_2 v_2 +
                  p_{2x} - \hat{\lambda}_2 p_2 +
                  e(t,\, x),
    \end{array}
\end{equation}
where the error term $e(t, \, x)$ satisfies the estimate
\eqref{reasons_of_source_term_eq} in Paragraph 8.2.2.

Deriving \eqref{u_t}, one obtains
\begin{equation}
\label{u_t_x}
\begin{split}
      u_{ t x}= &\Big( v_{1 xx} - ( \lambda_1 v_1 )_x \Big) \tilde{r}_1 +
                 \Big( v_1 \big( v_{ 1 xx} - ( \lambda_1 v_1 )_x  \big)
                        + 2 v_{1 x} \big( v_{1 x} - \lambda_1 v_1 \big)+
                        ( v_1^2 \sigma_1 )_x
                 \Big) \tilde{r}_{1 v}   \\
           &  +  \Big(  v_1 \big( v_{ 1 x} - \lambda_1 v_1 \big)
                 \Big)  \mathrm{D}\tilde{r}_1 \tilde{r}_1 +
                 \Big(  v_2 \big( v_{1 x} - \lambda_1 v_1 \big) +
                 \big( v_1 v_2 \big)_x
                 \Big)  \mathrm{D} \tilde{r}_1 r_2     \\
           &  +  \Big(   p_1 \big( v_{1 x} - \lambda_1 v_1 \big) + ( v_1 p_1 )_x
                 \Big)  \mathrm{D}\tilde{r}_1 \hat{r}_1 +
                 \Big(   p_2 \big( v_{1 x} - \lambda_1 v_1 \big) + ( v_1 p_2 )_x
                 \Big)  \mathrm{D}\tilde{r}_1 r_2    \\
           &  +  \Big(  \sigma_{1 x} ( v_{1 x} - \lambda_1 v_1 ) +
                         ( v_1 \sigma_{1 x})_x
                 \Big) \tilde{r}_{1 \sigma} +
                 \Big(  v_1 \big( v_{ 1 x} -  \lambda_1 v_1  \big)
                        +  v_1^2 \sigma_1
                 \Big) ( \tilde{r}_{1 v} )_x  \\
           &  +  v_1 v_2 (  \mathrm{D}\tilde{r}_1 r_2 )_x +
                 v_1 p_1 (  \mathrm{D}\tilde{r}_1 \hat{r}_1 )_x +
                 v_1 p_2 (   \mathrm{D}\tilde{r}_1 r_2 )_x +
                 v_1 \sigma_{1 x} ( \tilde{r}_{ 1 \sigma} )_x
                 \phantom{\Big)}\\
           &  +  \Big(  p_{1 xx} - ( \lambda_1 p_1 )_x
                 \Big) \hat{r}_1 +
                 \Big( p_1 \big(  p_{1 xx} - ( \lambda_1 p_1 )_x \big) +
                        2 p_{1 x} \big( p_{1 x}- \lambda_1 p_1 \big)
                 \Big) \hat{r}_{1 p} \\
           &  +  \Big( v_1 \big( p_{1 x} - \lambda_1 p_1 \big) +
                        ( v_1 p_1 )_x
                 \Big) \mathrm{D} \hat{r}_1 \tilde{r}_1 +
                 \Big( v_2 \big( p_{1 x} - \lambda_1 p_1 \big) +
                        ( v_2 p_1 )_x
                 \Big) \mathrm{D} \hat{r}_1 r_2
                 \phantom{\Big(}\\
           &  +  \Big( p_1 \big( p_{1 x} - \lambda_1 p_1 \big)
                 \Big) \mathrm{D} \hat{r}_1 \hat{r}_1 +
                 \Big(  p_2 \big( p_{1 x} - \lambda_1 p_1 \big)
                 \Big)
                 \mathrm{D} \hat{r}_1 r_2+
                 \Big( p_1 \big( p_{1 x} - \lambda_1 p_1 \big)
                 \Big) ( \hat{r}_{1 p} )_x
                 \phantom{\Big(} \\
           &  +  v_1 p_1 ( \mathrm{D} \hat{r}_1 \tilde{r}_1 )_x +
                 v_2 p_1 (  \mathrm{D} \hat{r}_1 r_2 )+
                   \Big( v_{2 xx} - ( \lambda_2 v_2 )_x \Big) r_2 +
                 \Big( p_{2 xx} - ( \lambda_2 p_2 )_x \Big) r_2 .
                 \phantom{\Big(} \\
\end{split}
\end{equation}
On the other hand,
\begin{equation}
\label{u_x_t}
      u_{x t} =  v_{1 t} \tilde{r}_1 + v_1 \tilde{r}_{1 t} +
                 v_{2 t} r_2 +
                 p_{1 t} \hat{r}_1 + p_1 \hat{r}_{1 t}+
                 p_{2 t} r_2,
\end{equation}
where
\begin{equation}
\label{derivative_of_vectors_wrt_t}
      \begin{array}{ll}
      \tilde{r}_{1 t} = \mathrm{D} \tilde{r}_1
                        (  w_1 \tilde{r}_1 + w_2 r_2  ) +
                        v_{1 t} \tilde{r}_{ 1 v} +
                        \sigma_{ 1 t} \tilde{r}_{1 \sigma} \\
      \hat{r}_{1 t}   = \mathrm{D} \hat{r}_1
                        (  w_1 \tilde{r}_1 + w_2 r_2  ) +
                        p_{1 t} \hat{r}_{ 1 p}.
     \end{array}
\end{equation}
We equal \eqref{u_t_x} and \eqref{u_x_t} and we use
\eqref{derivative_of_vectors_wrt_t}, obtaining
\begin{equation*}
      \begin{split}
      0   & =           \, u_{t x} - u_{x t} \\
          & =         \Big( v_{1 xx} - (\lambda_1 v_1)_x -v_{1 t} \Big)
                              \tilde{r}_1
      +   \Big( v_1 \big( v_{1 xx} - (\lambda_1 v_1)_x -v_{1 t}\big)+
                            2 v_{1 x} \big( v_{1 x} - \lambda_1 v_1 \big)+
                           ( v_1^2 \sigma_1 )_x
                      \Big) \tilde{r}_{1 v} \\
          & \quad +   \Big( v_1 ( v_{1 x} - \lambda_1 v_1 ) - v_1 w_1
                      \Big) \mathrm{D}\tilde{r}_1 \tilde{r}_1 +
                      \Big( v_2 ( v_{1 x} - \lambda_1 v_1 ) +
                          ( v_1 v_2)_x  - v_1 w_2
                      \Big) \mathrm{D}\tilde{r}_1 r_2 \\
          & \quad +   \Big( p_1 \big( v_{1 x} - \lambda_1 v_1 \big) +
                            ( v_1 p_1 )_x
                      \Big) \mathrm{D}\tilde{r}_1 \hat{r}_1 +
                      \Big( p_2 \big( v_{1 x} - \lambda_1 v_1 \big) +
                            ( v_1 p_2 )_x
                      \Big) \mathrm{D}\tilde{r}_1 r_2  \\
          & \quad
      +   \Big( \sigma_{1 x} \big( v_{1 x} - \lambda_1 v_1 \big) +
                            ( v_1 \sigma_{1 x} )_x -
                            \sigma_{1 t} v_1
                      \Big) \tilde{r}_{1 \sigma}
     +    \Big( v_{1 x} v_1 + v_1^2 ( \sigma_1 - \lambda_1 )
                      \Big) ( \tilde{r}_{ 1 v})_x +
                      v_1 v_2 ( \mathrm{D} \tilde{r}_1 r_2 )_x \\
         &  \quad +   v_1 p_1 ( \mathrm{D} \tilde{r}_1 \hat{r}_1 )_x +
                      v_1 p_2 ( \mathrm{D} \tilde{r}_1 r_2  )_x +
                      v_1 \sigma_{1 x} ( \tilde{r}_{ 1 \sigma} )_x +
                      \Big( p_{1 xx} - (\lambda_1 p_1 )_x -p_{1 t} \Big) \hat{r}_1
                      \phantom{\Big(}\\
         & \quad +    \Big( p_1  \big( p_{1 xx} -
                      (\lambda_1 p_1 )_x -p_{1 t} \big) +
                            2 p_{1 x} \big( p_{1 x} - \lambda_1 p_1 \big)
                      \Big) \hat{r}_{1 p} 
     +    \Big( v_1 \big( p_{1 x} - \lambda_1 p_1 \big) +
                            ( v_1 p_1 )_x - w_1 p_1
                      \Big) \mathrm{D} \hat{r}_1 \tilde{r}_1 \\
         & \quad +    \Big(  v_2 ( p_{1 x} - \lambda_1 p_1 ) +
                            ( v_2 p_1 )_x - w_2 p_1
                      \Big) \mathrm{D} \hat{r}_1 r_2
     +   p_1  \Big( p_{1 x} - \lambda_1 p_1 \Big)
                             \mathrm{D} \hat{r}_1 \hat{r}_1 +
                      \Big( p_2 \big( p_{1 x} - \lambda_1 p_1 \big) \Big)
                      \mathrm{D} \hat{r}_1 r_2 \\
         & \quad +    \Big(  p_1 ( p_{1 x} - \lambda_1 p_1 ) \Big)
                      ( \hat{r}_{1 p} )_x +
                      v_1 p_1 ( \mathrm{D} \hat{r}_1 \tilde{r}_1 )_x +
                      v_2 p_1 ( \mathrm{D} \hat{r}_1 r_2 )_x
         +    \Big( v_{2 xx} -  ( \lambda_2 v_2 )_x - v_{2 t} \Big) r_2 \\
     & \quad +
                      \Big( p_{2 xx} -  ( \hat{\lambda}_2 p_2 )_x - p_{2 t} \Big) r_2
                      \phantom{\Big(}\\
         & =          \Big( v_{1 xx} - (\lambda_1 v_1)_x -v_{1 t} \Big)
                                    \tilde{r}_1 +
                      \Big( p_{1 xx} -  (\lambda_1 p_1 )_x -p_{1 t} \Big)
                                    \hat{r}_1
                       \phantom{\Big(}             \\
         & \quad +    \Big( v_{2 xx} -  ( \lambda_2 v_2 )_x -
                      v_{2 t} \Big) r_2 +
                      \Big( p_{2 xx} -  ( \hat{\lambda}_2 p_2 )_x
                      - p_{2 t} \Big) r_2 +
                       s_1 ( t,  x). \\
      \end{split}
\end{equation*}
We can therefore impose the conditions
\begin{equation*}
      \begin{array}{llll}
      v_{1 t} + ( \lambda_1 v_1 )_x - v_{1 xx} = 0 \\
      p_{1 t} + ( \lambda_1 p_1 )_x - p_{1 xx} = 0 \\
      v_{2 t} + ( \lambda_2 v_2 )_x - v_{2 xx} =
      \, \langle  \tilde{\ell}_2 ( t, x) , s_1 ( t, x) \rangle  \, =
                         \tilde{s}_1 ( t, x)\\
      p_{2 t} + ( \hat{\lambda}_2 p_2 )_x - p_{2 xx} = 0,
      \end{array}
\end{equation*}
where $ ( \ell_1, \tilde{\ell}_2) $ is the dual basis of $ (
\tilde{r}_1, r_2 )$.

To compute the equations satisfied by $w_1,
\; w_2$ we will use
\begin{equation*}
      \begin{split}
             u_{tt} =
      &~      u_{xx \, t} - \big( A(u) u_x \big)_t \\
      =&~ u_{t \, xx } -
             \big( A(u) u_t \big)_x + DA(u)
             \big( u_x \otimes u_t
                   - u_t \otimes u_x \big),\\
      \end{split}
\end{equation*}
which follows from
\begin{equation*}
      \begin{split}
      &      \big( A(u) u_x \big)_t =
             DA (u) (u_t \otimes u_x) +
             A ( u) u_{xt} \\
      &      \big( A(u) u_t \big)_x =
             DA(u) ( u_x \otimes u_t) +
             A(u) u_{tx}. \\
      \end{split}
\end{equation*}
Straightforward computations ensures that
\begin{equation*}
      \begin{split}
             u_{x  t} - A(u) u_t =
      &      ( w_{ 1x} - \lambda_1 w_1) \tilde{r}_1 +
             w_1 ( v_{1x} - \lambda_1 v_1 ) \tilde{r}_{1 v} \\
      &      \quad +
             w_1 v_1 \sigma_1 \tilde{r}_{1 v} +
             w_1 v_2 \mathrm{D} \tilde{r}_1 r_2 +
             w_1 p_1 \mathrm{D} \tilde{r}_1 \hat{r}_1 +
             w_1 p_2 \mathrm{D} \tilde{r}_1 r_2 \\
      &      \quad +
             w_1 \sigma_{1 x} \tilde{r}_{1 \sigma} +
             ( w_{ 2 x} -  \lambda_2 w_2 ) r_2  \\
      \end{split}
\end{equation*}
and
\begin{equation}
\label{eq_derivation}
      \begin{split}
             DA(u)
      &       \bigg( u_x \otimes u_t
                  - u_t \otimes u_x \bigg)   =
             v_1 w_2 \tilde{r}_1 \otimes r_2 +
             v_2 w_1 r_2 \otimes \tilde{r}_1 +
             p_1 w_1 \hat{r}_1 \otimes \tilde{r}_1 +
             p_1 w_2 \hat{r}_1 \otimes r_2 +
             p_2 w_1 r_2 \otimes \tilde{r}_1 \\
      &      \quad -
             w_1 v_2 \tilde{r}_1 \otimes r_2 -
             w_1 p_1 \tilde{r}_1 \otimes \hat{r}_1 -
             w_1 p_2 \tilde{r}_1 \otimes r_2 -
             w_2 v_1 r_2 \otimes \tilde{r}_1 -
             w_2 p_1 r_2 \otimes \hat{r}_1. \\
      \end{split}
\end{equation}
Hence
\begin{equation*}
      \begin{split}
            u_{tt}
      &     = w_{ 1 t} \tilde{r}_1 +
            w_{ 2 t} r_2 \phantom{\Big(} \\
      &     =
            - w _1^2  \mathrm{D} \tilde{r}_1 \tilde{r}_1  -
            w_1 w_2 \mathrm{D} \tilde{r}_1 r_2 -
            w_1 v_{1 t}  \tilde{r}_{1 v} -
            w_1 \sigma_{1 t} \tilde{r}_{ 1 \sigma} +
            \Big( w_{1 xx} - ( \lambda_1 w_1 )_x  \Big) \tilde{r}_1 +
            \Big( v_1 \big( w_{ 1x} - \lambda_1 w_1 \big) \Big)
            \mathrm{D} \tilde{r}_1 \tilde{r}_1 \\
      &     \quad +
            \Big( v_2 ( w_{1 x} - \lambda_1 w_1 ) \Big)  \mathrm{D} \tilde{r}_1 r_2 +
            \Big( p_1 ( w_{ 1x } - \lambda_1 w_1 ) \Big)
            \mathrm{D} \tilde{r}_1 \hat{r}_1 +
            \Big( p_2 ( w_{1x} - \lambda_1 w_1 ) \Big) \mathrm{D} \tilde{r}_1 r_2 +
            \Big( v_{1 x} ( w_{1 x} - \lambda_1 w_1 )\Big) \tilde{r}_{1v} \\
        &   \quad +
            \Big( \sigma_{1 x} ( w_{1 x} - \lambda_1 w_1 ) \Big)
            \tilde{r}_{ 1 \sigma} +
            \Big( w_{1 x} ( v_{1 x} - \lambda_1 v_1 ) \Big)
             \tilde{r}_{1 v} +
            \Big( w_1 ( v_{1 x} - \lambda_1 v_1 ) \Big)
            ( \tilde{r}_{1 v} )_x +
            w_1 v_1 \sigma_1 ( \tilde{r}_{1 v} )_x \\
      &     \quad +
            ( w_1 v_1 \sigma_1 )_x \tilde{r}_{ 1v } +
            ( w_1 v_2 )_x  \mathrm{D} \tilde{r}_1 r_2 +
            w_1 v_2 (  \mathrm{D} \tilde{r}_1 r_2 )_x +
            ( w_1 p_1 )_x  \mathrm{D} \tilde{r}_1 \hat{r}_1 +
            w_1 p_1 (  \mathrm{D} \tilde{r}_1 \hat{r}_1 )_x
            \phantom{\Big(} \\
      &     \quad +
            ( w_1 p_2 )_x  \mathrm{D} \tilde{r}_1 r_2 +
            w_1 p_2 (  \mathrm{D} \tilde{r}_1 r_2 ) _x +
            ( w_1 \sigma_{ 1x} )_x \tilde{r}_{1 \sigma} +
            w_1 \sigma_{1 x} ( \tilde{r}_{1 \sigma} )_x +
            \Big( w_{2 xx} - ( \lambda_2 w_2 )_x \Big) r_2 \\
      &     \quad +
            ( p_1 w_2 )_x  \mathrm{D} \hat{r}_1 r_2 +
            p_1 w_2 (  \mathrm{D} \hat{r}_1 r_2 )_x
            \phantom{\Big(} \\
      &     \quad + DA(u) \Big(
             v_1 w_2 \tilde{r}_1 \otimes r_2 +
             v_2 w_1 r_2 \otimes \tilde{r}_1 +
             p_1 w_1 \hat{r}_1 \otimes \tilde{r}_1 +
             p_1 w_2 \hat{r}_1 \otimes r_2 +
             p_2 w_1 r_2 \otimes \tilde{r}_1 \\
      &      \quad -
             w_1 v_2 \tilde{r}_1 \otimes r_2 -
             w_1 p_1 \tilde{r}_1 \otimes \hat{r}_1 -
             w_1 p_2 \tilde{r}_1 \otimes r_2 -
             w_2 v_1 r_2 \otimes \tilde{r}_1 -
             w_2 p_1 r_2 \otimes \hat{r}_1 \Big)  \\
      &      =   \Big( w_{1 xx} - ( \lambda_1 w_1 )_x  \Big) \tilde{r}_1 +
              \Big( w_{2 xx} - ( \lambda_2 w_2 )_x \Big) r_2 +
             s_2 (t, \, x).\\
      \end{split}
\end{equation*}
One can check that, since $A$ is triangular,
\begin{equation*}
      \langle \ell_1, \, DA(u)
             \big( u_x \otimes u_t
                   - u_t \otimes u_x \big) \rangle  =0
\end{equation*}
and therefore the equations satisfied by $w_i \; i=1, \, 2$ are
\begin{equation}
\label{w}
      \begin{split}
       &     w_{ 1t} + ( \lambda_1 w_1 )_x - w_{1 xx} = 0 \\
       &     w_{ 2t} + ( \lambda_2 w_2 )_x - w_{ 2 xx} =
             \langle  \tilde{\ell}_2(t, \, x), \, s_2 (t, \, x) \rangle  =
             \tilde{s}_2 (t, \, x). \\
      \end{split}
\end{equation}

\subsubsection{Proof of Proposition \ref{reasons_of_source_term}}
\label{source_pro} Equation \eqref{decomposition} and
\eqref{eq_theta} ensure that, since,
\begin{equation*}
   \sigma_1 = \lambda_1^{\ast} -
                 \theta \bigg( \frac{w_1}{v_1} +
                              \lambda_1^{\ast}
                        \bigg),
\end{equation*}
then
\begin{equation*}
      \begin{split}
      &     \sigma_{ 1 x} = - \theta'\bigg( \frac{w_1}{v_1} +
                              \lambda_i^{\ast}
                        \bigg) \,
            \Bigg( \frac{w_1}{v_1} \Bigg)_x =
            -
            \bigg( \frac{w_{ 1x } v_1 - v_{1x} w_1}{ v_1^2}
            \bigg) \theta',  \\
      &     \qquad \qquad
            | v_1^2 \sigma_{1x} |= \mathcal{O}(1)
            | w_{1 x} v_1 - v_{1 x} w_1 |, \phantom{\Bigg(}\\
      &     \; \qquad \qquad \sigma_{1 x} \neq 0 \iff
            \bigg| \frac{w_1}{v_1} - \lambda_1^{\ast}
            \bigg|  \leq 3 \hat{\delta}. \phantom{\Bigg(}\\
      \end{split}
\end{equation*}
Most of the terms in $\tilde{s}_i(t, \, x) \; i=1, \, 2$ and $e(t,
\, x)$ can be directly reduced to those in Proposition
\ref{reasons_of_source_term}. The terms which requires some
technicalities are:
\begin{enumerate}
\item
\begin{equation*}
      |p_{1 x } - \lambda_1 p_1 | \,  | \langle \tilde{\ell}_2(u, \, v_1, \, \sigma_1),
      \; \hr(u, \, p_1) \rangle|
      \leq \unpo (|p_1| + |v_1 |)  |p_{1 x } - \lambda_1 p_1 |.
\end{equation*}
Indeed,
\begin{equation*}
\begin{split}
&     | \langle \tilde{\ell}_2(u, \, v_1, \, \sigma_1),
      \; \hr(u, \, p_1) \rangle | \leq
      |\langle  \tilde{\ell}_2, \; \hr - r_1^{\ast}  \rangle | +
      |\langle \tilde{\ell}_2 - \ell^{\ast}_2, \, r_1^{\ast}  \rangle |
      \leq \unpo (|p_1| + |v_1|).
\end{split}
\end{equation*}
We have denoted by $r_1^{\ast}$ the first eigenvector of the
matrix $A(u^{\ast})$ and by $(\ell_1, \, \ell^{\ast}_2)$ the dual
base
of $(r_1^{\ast}, \, r_2).$ \\
\item
\begin{equation*}
        \begin{split}
        \quad &  \bigg| 2 v_{1 x} ( v_{1 x} - \lambda_1 v_1 )+
                  ( v_1^2 \sigma_1 )_x  \bigg|=
        \bigg| 2 v_{1 x}
          \bigg( w_1 -
                (p_{1x}-\lambda_1 p_1)
          \bigg)+
          2 v_1 v_{1x} \sigma_1  +
          v_1^2 \sigma_{1 x}\bigg|  \\
        & \leq  \bigg| 2 v_{1 x} ( w_1 + v_1 \sigma_1 )\bigg|  +
             \bigg| 2 v_{1 x} ( \lambda_1 p_1 - p_{1 x})\bigg|  +
            \mathcal{O}(1) \bigg| v_{1 x} w_1 - v_1 w_{1 x} \bigg|.  \\
        \end{split}
\end{equation*}
\item
$  \Big| \sigma_{1 x} ( v_{1 x} - \lambda_1 v_1 ) +
             ( v_1 \sigma_{1 x} )_x -
              \sigma_{1 t} v_1
      \Big|   | \tilde{r}_{1 \sigma} |
$: some computations ensures that
$$
      \Big( \frac{w_1}{v_1} \Big)_x
       (v_{1x} - \lambda_1 v_1) +
       v_{1 x}  \Big( \frac{w_1}{v_1} \Big)_x +
       v_1  \Big( \frac{w_1}{v_1} \Big)_{xx} -
       \Big( \frac{w_1}{v_1} \Big)_t v_1
       = 0.
$$
Hence, since $ | \tilde{ r}_{1 \sigma} | = \mathcal{O}(1) | v_1|$,
one gets
$$
      \Big| \sigma_{1 x} ( v_{1 x} - \lambda_1 v_1 ) +
             ( v_1 \sigma_{1 x} )_x -
              \sigma_{1 t} v_1
      \Big| \,  | \tilde{r}_{1 \sigma} | \,
      \leq   \mathcal{O}(1)
      \chi_{ \{ |\lambda_1^{\ast} - w_1 / v_1| \leq 3 \hat{\delta} \} }  v_1^2
      \Big| \Big( \frac{w_1}{v_1} \Big)_x \Big|^2.
$$
\item 
$
   | - w_1 \sigma_{1 t} + \sigma_{ 1x} ( w_{1 x} - \lambda_1 w_1) +
   ( w_1 \sigma_{1 x} )_x | \, | \tilde{r}_{1 \sigma}|:
$ \\
since
$$
   -w_1 \theta' \bigg( \frac{w_1}{ v_1} \bigg)_t +
   w_{1x}  \theta' \bigg( \frac{w_1}{ v_1} \bigg)_{x} -
   \lambda_1 w_1 \theta' \bigg( \frac{w_1}{ v_1} \bigg)_x +
   w_{1x} \theta' \bigg( \frac{w_1}{ v_1} \bigg)_x +
   w_1 \theta' \bigg( \frac{w_1}{ v_1} \bigg)_{xx} =0,
$$
one is left to the estimate
$$
  \bigg| \,  \theta'' \bigg( \frac{w_1}{v_1} \bigg)_x \bigg|^2
  | w_1 v_1| \leq
  \mathcal{O}(1)
  v_1^2  \chi_{ \{ |\lambda_1^{\ast} - w_1 / v_1| \leq 3 \hat{\delta} \} }
  \bigg| \bigg( \frac{w_1}{v_1} \bigg)_x \bigg|^2
 .
$$
\item
$
 |  v_1 \sigma_{1x} ( \tilde{r}_{1 \sigma})_x|
$ : first of all, we observe that that $ \theta'(s) \neq 0 $
implies  $ |w_1| \leq \mathcal{O}(1) |v_1|$
 and therefore
$$
     | v_1 \sigma_{1x}|=
     |v_1 \theta' |
     \Big|
          \Big( \frac{w_1}{v_1} \Big)_x
     \Big| =
     \Big|
         \frac{w_{1x}v_1 -
               v_{1x}w_1}
              {v_1^2}
     \Big| |v_1 \theta'|
     \leq
     \mathcal{O}(1)( |w_{1x}|+
                     |v_{1x}|).
$$
We develop
\begin{equation*}
     | ( \tilde{r}_{1 \sigma})_x| =
     |  ( \tilde{r}_{1x})_{ \sigma}| \leq
     \mathcal{O}(1)
     \bigg(
          | v_1| + |v_2 |+ |p_1| +| p_2 |+ |v_{1x}|
     \bigg) +
     \mathcal{O}(1)|v_1 \sigma_{1x}|.
\end{equation*}
Since
\begin{equation*}
     \begin{array}{ccccc}
     \theta' \neq 0 &
     \Rightarrow &
     |w_1| = |v_{1x} - \lambda_1 v_1+
              p_{1x} - \lambda_1 p_1|
     \leq \mathcal{O}(1) |v_1| &
     \Rightarrow &
     |v_{1x}| \leq
     \mathcal{O}(1) |v_1|+
     |p_{1x} - \lambda_1 p_1|,
    \end{array}
\end{equation*}
one has
\begin{equation*}
     \begin{split}
     | v_{1x} \sigma_{1x} v_1| =&~
     \Big|
         \frac{w_{1x}v_1 -
               v_{1x}w_1}
              {v_1^2}
     \Big| |\theta' v_1 v_{1x}| \\
  \leq&~  \mathcal{O}(1) |w_{1x}v_1 - w_1 v_{1x}|+
     \mathcal{O}(1) \bigg(|w_{1x}| +
     \mathcal{O}(1) |v_{1x}|\bigg)
     |p_{1x}-\lambda_1 p_1|.\\
    \end{split}
\end{equation*}
Using the previous estimates, we get
\begin{equation*}
      \begin{split}
      &   |  v_1 \sigma_{1x} ( \tilde{r}_{1 \sigma})_x|
         \leq
          \mathcal{O}(1) |w_{1x}v_1 -
                          v_{1x} w_1|+
          \mathcal{O}(1)
          \big(
                |v_1| + |w_{1x}|
          \big)
          \,
          \big(
               |v_2|+ |p_1|+|p_2|
          \big)           \\
      &   \quad +
           \mathcal{O}(1) |w_{1x}v_1 - w_1 v_{1x}|+
             \mathcal{O}(1) \bigg(|w_{1x}| +
          \mathcal{O}(1) |v_{1x}|\bigg)
         |p_{1x}-\lambda_1 p_1 |+
          \mathcal{O}(1) v_1^2
           \chi_{ \{ |\lambda_1^{\ast}  - w_1 / v_1 | \leq 3 \hat{\delta} \} }
          \Big( \frac{w_1}{v_1} \Big)_x^2. \\
      \end{split}
\end{equation*}
\item
\begin{equation*}
\begin{split}
     |v_1 (w_{1 x} - \lambda_1 w_1 ) - w_1^2| =&~
      | v_1 w_{1 x} - v_{1 x}w_1 + v_{1x} w_1
       - \lambda_1 v_1 w_1 -w_1^2| \\
     \leq&~  | v_1 w_{1 x} - v_{1 x}w_1 | +
       |w_1 (v_{1 x}  - \lambda_1 v_1 - w_1)| \\
       \leq&~
        | v_1 w_{1 x} - v_{1 x}w_1 | +
       |w_1 (p_{1 x}  - \lambda_1 p_1)|. \\
\end{split}
\end{equation*}
\item
\begin{equation*}
      \begin{split}
      | w_{1x} ( v_{1 x} - \lambda_1 v_1 &) + ( w_1 v_1 \sigma_1)_x +
             (  w_{1x } - \lambda_1 w_1 )v_{1x} | \\
      &       = | 2 w_{1 x} (v_{1 x} - \lambda_1 v_1 - w_1 )
             - w_{1 x} v_{1 x } + \lambda_1 w_{1 x} v_1 +
             2 w_{1 x} w_1  \\
      &      \quad  + w_{1 x} v_1 \sigma_1 + w_1 v_{1 x} \sigma_1 + w_1 v_1 \sigma_{1 x} +
             w_{1 x} v_{1 x} - \lambda_1 w_1 v_{1 x} | \\
       &     = | 2 w_{1 x} (p_{1 x} - \lambda_1 p_1 )+
             2 w_{1 x} (w_1 + \sigma_1  v_1 ) -
             \sigma_1 w_{1 x} v_1 \\
      &      \quad + \lambda_1 w_{1 x} v_1 +
             \sigma_1 w_1 v_{1x } + (w_{1 x} v_1 - w_1 v_{1 x})
             \theta' (w_1 / v_1 ) - \lambda_1 w_1 v_{1 x} |
             \\
      &      \leq 2 |  w_{1 x} (p_{1 x} - \lambda_1 p_1 ) | +
             2 |  w_{1 x} (w_1 + \sigma_1  v_1 ) |  \\
      &      \quad + |(\lambda_1 - \sigma_1 ) (w_{1 x} v_1 - w_1 v_{1 x})|
             + (w_{1 x} v_1 - w_1 v_{1 x})
             \theta' (w_1 / v_1 ) | \\
      \end{split}
\end{equation*}
\item
$$
  | w_1 ( w_1 + \sigma_1 v_1 - p_{1 x} + \lambda_1 p_1 )| =
   | w_1 ( w_1 + \sigma_1 v_1 ) +
  w_1 ( p_{1 x} - \lambda_1 p_1 )|
$$
\item
$$
   | w_1 \sigma_{1 x} ( \tilde{r}_{1 \sigma})_x | \leq
    | v_1 \sigma_{1 x} ( \tilde{r}_{1 \sigma})_x |,
$$
and therefore one comes back to case (5).
\end{enumerate}
This completes the proof of the estimate
\eqref{reasons_of_source_term}.

\subsection{Appendix to Paragraph \ref{BV_estimates}}

\subsubsection{Proof of the estimate \eqref{exp_decay_px}}
\label{exp_decay_px_proof}

It is convenient to introduce a representation formula for $p_i,
\; i=1, \, 2$. To this end, two new convolution kernels are
needed: let $I^{\lambda_i^{\ast} \, 0}(t, \, s, \, x)$ be the
solution of the equation
\begin{equation*}
      I^{\lambda_i^{\ast} \, 0}_t + \lambda_i^{\ast}I^{\lambda_i^{\ast} \, 0}_x
      - I^{\lambda_i^{\ast} \, }_{xx}=0,
\end{equation*}
with boundary and initial data
\begin{equation*}
      I^{\lambda_i^{\ast} \, 0}(0, \,s, \,  x) \equiv 0,
      \quad
      I^{\lambda_i^{\ast} \, 0}(t, \, s, \, 0) = \delta_{t =s},
      \quad
      I^{\lambda_i^{\ast} \, 0}(t, \, s, \, L) \equiv 0.
\end{equation*}
Without specifying the explicit expression of $I^{\lambda_i^{\ast}
\, }$, we observe that
\begin{equation*}
      \int_0^{+ \infty} I^{\lambda_i^{\ast} \,0 }(t, \, s, \, x) ds =
      J^{\lambda_i^{\ast} \, 0}(t, \, x)
\end{equation*}
(see equation \eqref{J_0} for the definition of
$J^{\lambda^{\ast}_i \, 0}$). Analogously, let $
I^{\lambda_i^{\ast} \, L}(t, \, x)$ be the solution of the
equation
\begin{equation*}
      I^{\lambda_i^{\ast} \, L}_t + \lambda_i^{\ast}I^{\lambda_i^{\ast} \, L}_x
      - I^{\lambda_i^{\ast} \, L}_{xx}=0,
\end{equation*}
with boundary and initial data:
\begin{equation*}
      I^{\lambda_i^{\ast} \, 0}(0, \, s, \, x) \equiv 0,
      \quad
      I^{\lambda_i^{\ast} \, 0}(t, \, s, \,  0) \equiv 0,
      \quad
      I^{\lambda_i^{\ast} \, 0}(t, \, s, \,  L) = \delta_{t =s}.
\end{equation*}
By construction, it satisfies
\begin{equation*}
       \int_0^{+ \infty} I^{\lambda_i^{\ast} \, L }(t, \, s, \, x) ds =
      J^{\lambda_i^{\ast} \, L}(t, \, x)
\end{equation*}
(see equation \eqref{eq_J_L} for the definition of
$J^{\lambda_i^{\ast} \, L}$). If $t \leq 1$ the function $p_1$
admits the following representation formula:
\begin{equation*}
\begin{split}
     p_1(t, \, x) =&~ \int_0^{+ \infty}
      I^{\lambda_1^{\ast} \, 0}(t, \, s, \,  x)p_1(s, \, 0) ds +
      \int_0^{+ \infty}
      I^{\lambda_1^{\ast} \, L}(t, \, s, \,  x)p_1(s, \, L) ds +
      \\
&     +
      \int_0^t \int_0^L \Delta^{\lambda_1^{\ast}}( t-s, \, x, \,
      y) \Big( (\lambda_1^{\ast} - \lambda_1) p_{1 y} -
      \lambda_{1 y} p_1 \Big)
      (s, \, y)dy ds ,  \\
\end{split}
\end{equation*}
and hence
\begin{equation*}
\begin{split}
     p_{1 x}(t, \, x)=&~ \int_0^{+ \infty}
      I_x^{\lambda_1^{\ast} \, 0}(t, \, s, \,  x)p_1(s, \, 0) ds +
      \int_0^{+ \infty}
      I_x^{\lambda_1^{\ast} \, L}(t, \, s, \,  x)p_1(s, \, L) ds \\
&     + \int_0^t \int_0^L
      \Delta_x^{\lambda_1^{\ast}}( t-s, \, x, \, y)
      \Big( (\lambda_1^{\ast} - \lambda_1) p_{1 y} - \lambda_{1
      y}p_1
      \Big)(s, \, y) dy ds.\\
\end{split}
\end{equation*}
From the expression of $\Delta^{\lambda_1^{\ast}}$, given by
formula \eqref{Delta_product}, it follows that
\begin{equation*}
      \bigg\|
      \Delta_x^{\lambda_1^{\ast}}(t, \, \,  \cdot\, , \, y)
      \exp(c(\, \cdot \, -y) / 2 )
      \bigg\|_{L^1} \leq \frac{\unpo}{\sqrt{t}},
\end{equation*}
and from the previous observations
\begin{equation*}
\begin{split}
&      \bigg| \exp(c x / 2) \int_0^{\infty}
       I_x^{\lambda_i^{\ast} \,0 }(t, \, s, \, x) ds
       \bigg| =
       \big| \exp (c x/ 2) J_x^{\lambda_i^{\ast} \, 0}(t, \, x)
       \big| \leq \unpo \\
&      \qquad \qquad \bigg| \int_0^{+ \infty}
       I_x^{\lambda_i^{\ast} \, L}(t, \, s, \, x) ds
       \bigg| =
       \big| J_x^{\lambda_i^{\ast} \, L}
       (t, \, x) \big|.\\
\end{split}
\end{equation*}
Hence
\begin{equation*}
\begin{split}
     \big| \exp& ( c x/ 2) p_1(t, \, x) \big| =
      \bigg| \exp (c x/2 ) \int_0^{ + \infty}
      I_x^{\lambda_1^{\ast} \, 0}(t, \, s, \,  x)p_1(s, \, 0) ds
      \bigg| +
      \bigg| \exp (cx / 2)\int_0^{+ \infty}
      I_x^{\lambda_1^{\ast} \, L}(t, \, s, \,  x)
      p_1(s, \, L) ds \bigg| \\
&     + \bigg| \exp (cx /2)\int_0^t \int_0^L
      \Delta_x^{\lambda_1^{\ast}}( t-s, \, x, \,
      y) \Big( (\lambda_1^{\ast} - \lambda_1) p_{1 y}
      - \lambda_{ 1 y} p_1  \Big)
      (s, \, y) dy ds \bigg| \\
      \leq&~
      \unpo |p(x= 0) |_{\infty} +
      \unpo |p(x= L) |_{\infty}  \phantom{\bigg|}\\
&     +  \unpo \delta_1
      \bigg| \int_0^t \bigg( \sup_{y}p_{1 y}( s, \, y)
      \exp (c y / 2) \bigg)
      \int_0^L \Delta_x^{\lambda_1^{\ast}}(t-s, \, x, \, y)
      \exp\big(c (x -y)/2  \big) ds dy \bigg| + \unpo \delta_1^2\\
\end{split}
\end{equation*}
and therefore
\begin{equation*}
       | \sup_{x}p_{1 x}( t, \, x)
       \exp (c x/ 2) | \leq \unpo \delta_1
       \quad \forall \, t \leq 1.
\end{equation*}
The estimate
\begin{equation*}
       \sup_{x} \big| p_{2 x}( t, \, x)
       \exp \big( c( L- x) /2 \big) \big|
       \leq \unpo \delta_1
       \quad \forall \, t  \leq 1
\end{equation*}
follows by symmetry.

If $ t>1$ the following representation formula holds:
\begin{equation*}
\begin{split}
      p_{1 x}(t, \, x)=
&     \int_0^L p_1 (t-1,\ , y) \Delta^{\lambda_1^{\ast}}(1, \, x,
      \, y) dy +
      \int_{t - 1}^{+ \infty}
      I_x^{\lambda_1^{\ast} \, 0}(1, s, \,  x)p_1(s, \, 0) ds
      \\
&     + \int_{t-1}^{+ \infty}
      I_x^{\lambda_1^{\ast} \, L}(1, \, s, \,  x) p_1(s, \, L) ds \cr
      & +
      \int_0^1 \int_0^L \Delta_x^{\lambda_1^{\ast}}( 1 - s, \, x, \,
      y) \Big( (\lambda_1^{\ast} - \lambda_1) p_{1 y} - \lambda_{1 y}
      p_1  \Big)(t- 1+s, \, y) dy ds.\\
\end{split}
\end{equation*}
It follows that
\begin{equation*}
      |\sup_x p_{1x} (t, \, x) \exp ( cx/2 )| \leq
      \unpo \delta_1 \quad \forall \, t> 1,
\end{equation*}
and hence by symmetry
\begin{equation*}
      \big| \sup_x p_{2 x} (t, \, x) \exp \big( c(L-x )/2
      \big) \big| \leq
      \unpo \delta_1 \quad \forall \, t> 1.
\end{equation*}
This concludes the proof of \eqref{exp_decay_px}.

\subsubsection{Proof of Proposition \ref{other_wrt_time_pro}}
\label{other_wrt_time_par}

We will perform the computations only for
$v_2, \; w_2$ and $w_{2 x}$, since those for $v_1, \; w_1$ and
$w_{1 x}$ follow by symmetry.

{\bf Three new convolution kernels:} the solution of equation
\begin{equation}
\label{eq_v}
      Q_t + \lambda_2^{\ast}Q_x - Q_{xx} =0
\end{equation}
with boundary conditions
\begin{equation*}
      Q(0, \, x) = \delta_y,
      \qquad
      Q(t, \, 0) = 0,
      \qquad
      Q_x(t, \, L)=0
\end{equation*}
is
\begin{equation}
\label{eq_def_theta}
      Q(t, \, x) = \Theta^{\driftd}(t, \, x, \, y) : =
      \int_0^x \phi(t, \, z, \, y) \bigg( \sum_m
      G_z (z + 2mL - y) + G_z (z + 2mL + y)  \bigg) dz
\end{equation}
As in Section \ref{parabolic_estimates}, we use the notation
\begin{equation*}
   \phi^{\, \lambda_2^{\ast}}(t, x, y) =
   \exp
      \bigg(
           \frac{
                 \, \, \lambda_2^{\ast}}{2}
           \, (x -y) -
           \frac{
                \, \, (  \lambda_2^{\ast})^2 }{4}t\,
      \bigg).
\end{equation*}
and $G(t, \, x)= \exp(-x^2/ 4t)/ 2 \sqrt{\pi t} $.
Note that, by construction,
\begin{equation}
\label{eq_theta_x}
      \Theta^{\lambda_2^{\ast}}_x (t, \, 0, \, y) \equiv 0
      \quad \forall \, t \ge 0, \quad y \in ]0, \, L[ \, .
\end{equation}
Moreover, an argument similar to that used in Section
\ref{par_maximum} ensures that a maximum principle holds for
equation \eqref{eq_v}, in other words if
\begin{equation*}
      Q(0, \, x) \leq 0, \qquad Q(t, \, 0) \leq 0, \qquad Q_x(t, \, L) \leq 0,
\end{equation*}
then $Q(t, \, x) \leq 0 \; \forall \, t, \; x$.

The solution of \eqref{eq_v} with boundary conditions
\begin{equation*}
       Q(0, \, x) = 0,
       \qquad
       Q(t, \, 0)= 1,
       \qquad
       Q_x(t, \, L) =0
\end{equation*}
is
\begin{equation}
\label{eq_conv_ker_b}
      B^{\driftd}(t, \, x) = 1 - \int_0^L \Theta^{\driftd}(t, \,
      x, \, y) dy.
\end{equation}
In the following, we will need another convolution kernel,
$\tilde{\Theta}^{\lambda_2^{\ast}}(t, \, x, \, y)$, such that
\begin{equation}
\label{eq_tilde_theta}
       \tilde{\Theta}_y^{\driftd}(t, \, x, \, y) =
       - \Theta_x^{\driftd}(t, \, x, \, y ).
\end{equation}
We arbitrarily impose $\tilde{\Theta}^{\driftd}(t, \, x, \, L)
\equiv 0 \; \forall \, t, \; x$ and define
\begin{equation*}
      \tilde{\Theta}^{\driftd}(t, \, x, \, y) : =
      \int_y^L \Theta_x^{\driftd}(t, \, x, \, \xi) d \xi.
\end{equation*}
Recalling \eqref{eq_theta_x}, we observe that
$\tilde{\Theta}^{\driftd}(t, \, x, \, y)$ is the derivative with
respect to $x$ of a function $z$ such that
\begin{equation}
\label{eq_zeta}
\begin{split}
&      z_x (t, \, 0, \, y) \equiv 0
       \quad
       z_x (t, \, L, \, y) \equiv 0
       \quad
       z(0,  \, x, \, y) =
       \left\{
       \begin{array}{ll}
              0 \quad 0 < x \leq y \\
              1 \quad y \leq x < L \\
       \end{array}
       \right. \\
&      \qquad
       \qquad
       \qquad
       \qquad
       z_t + \lambda_2^{\ast} z_x - z_{xx}=0. \\
\end{split}
\end{equation}
It follows that $\tilde{\Theta}^{\driftd}(t, \, x, \, y)$
satisfies
\begin{equation*}
      \tilde{\Theta}^{\driftd}(t, \, 0, \, y)
      \equiv 0
      \quad
      \tilde{\Theta}^{\driftd}(t, \, L, \, y)
      \equiv 0
      \quad
      \tilde{\Theta}^{\driftd}(0, \, x, \, y) =
      \delta_y
\end{equation*}
and hence actually
\begin{equation}
\label{eq_tilde_theta_delta}
      \tilde{\Theta}^{\driftd}(t, \, x, \, y)
      \equiv \Delta^{\driftd}(t, \, x, \, y),
\end{equation}
where $\Delta^{\driftd}$ is the convolution kernel defined by
\eqref{Delta_product}. In the following, however, for sake of
clearness we will write $\tilde{\Theta}^{\driftd}(t, \, x, \, y)$
when we want to underline that the relation \eqref{eq_tilde_theta}
holds. From the identity \eqref{eq_tilde_theta_delta} and the
estimates \eqref{estimate_kernels} it follows
\begin{equation}
\label{eq_theta_x_integral}
       \| \tilde{\Theta}^{\driftd}(t, \, x, \, y)\|_{L^1}
       \leq \unpo \quad
       \| \tilde{\Theta}_x^{\driftd}(t, \, x, \, y) \|_{L^1}
       \leq \frac{\unpo}{\sqrt{t}}
       \quad \forall \, t \leq 1.
\end{equation}
Moreover, let $z$ be as in \eqref{eq_zeta} and let $B^{\driftd}$
be defined by \eqref{eq_conv_ker_b}, then $z(t, \, x, \, 0) +
B^{\driftd}(t, \, x) \equiv 1$ and hence
\begin{equation}
\label{eq_theta_b}
      \tilde{\Theta}^{\driftd}(t, \,x, \, 0) +
      B_x^{\driftd}(t, \, x) = 0.
\end{equation}
Such an identity, together with \eqref{eq_theta_x_integral},
implies
\begin{equation}
       \|B_x^{\driftd}(t, \, x)\|_{L^1} \leq \unpo
       \quad
       \| B_{xx}^{\driftd}(t, \, x) \|_{L^1} \leq
       \frac{\unpo}{\sqrt{t}}
       \quad t \leq 1.
\end{equation}
Since the kernels introduced so far will be used to prove the
integrability of $v_{2x}$ with respect to time, one has to prove
that they are integrable on small time intervals.
\begin{itemize}
\item
\begin{equation}
\label{estimate_theta_x_wrtt}
         \int_0^1 |\tilde{\Theta}_x^{\lambda_2^{\ast}}
                 (t, \, x, \, y)| dt
         = \int_0^1 |\Delta_x^{\lambda_2^{\ast}}
                 (t, \, x, \, y)| dt
         \leq \mathcal{O}(1)
         \quad \forall x \in \, [0, \, L],
         \quad \forall \,  y \in \, ]0, \, L[
\end{equation}
\begin{proof}
One can check that
\begin{equation}
\label{eq_basic_wrtt}
       \int_0^1 | G_x^{\lambda_2^{\ast}}
       (t, \, x - y) | dt
       \leq \unpo
       \qquad
       \int_0^1 | G^{\lambda_2^{\ast}}
       (t, \, x - y) | dt
       \leq \unpo
       \quad \forall \, x, \, y \in \mathbb{R}.
\end{equation}
Since
\begin{equation*}
\begin{split}
     \Delta_x^{\driftd}(t, \, x, \, y) =
&    \bigg(
           \phi(t, \, x, \, y) \sum_{m \ge 0}
           G(t, \, x -y + 2 mL )
     \bigg)_x    -
     \bigg(
           \phi(t, \, x, \, y) \sum_{m \ge 0}
           G(t, \, x + y + 2 mL )
     \bigg)_x \\
&    +
     \bigg(
           \phi(t, \, x, \, y) \sum_{n > 0}
           G(t, \, x -y - 2 nL )
     \bigg)_x     -
     \bigg(
           \phi(t, \, x, \, y) \sum_{n > 0}
           G(t, \, x + y - 2 nL )
     \bigg)_x, \\
\end{split}
\end{equation*}
one gets
\begin{equation*}
\begin{split}
      | \Delta_x^{\driftd}(t, \, x, \, y) |
&     \leq
      \sum_{m \ge 0}
      |G_x^{\lambda_2^{\ast}}(t, \, x - y + 2mL )|
      +
      \sum_{m \ge 0}
      |G_x^{\lambda_2^{\ast}} (t, \, x+ y + 2mL)| \\
&     \quad  +
      \sum_{n > 0} |G_x^{\driftd}(t, \, 2nL + y -x)| +
      \driftd \sum_{n > 0} |G^{\driftd}(t, \, 2nL + y -x)| \\
&     \quad  +
      \sum_{n > 0} |G_x^{\driftd}(t, \, 2nL - y -x)| +
      \driftd \sum_{n > 0} |G^{\driftd}(t, \, 2nL - y -x)|. \\
\end{split}
\end{equation*}
Since
\begin{equation*}
     |G_x^{\lambda_2^{\ast}}(t, \, z + 2mL ) | \leq
      e^{-mL} |G_x^{\lambda_2^{\ast}}(t, \, z  ) |
     \qquad
      |G^{\lambda_2^{\ast}}(t, \, z + 2mL ) | \leq
      e^{-mL} |G^{\lambda_2^{\ast}}(t, \, z  ) |
\end{equation*}
if $m \ge 0$, $t \leq 1$ and $z$ is large enough, from the
previous estimates and from \eqref{eq_basic_wrtt} one deduces
\eqref{estimate_theta_x_wrtt}.
\end{proof}
\item From equation \eqref{eq_theta_b} and the previous estimate
it follows
\begin{equation}
\label{eq_b_xx_wrtt}
      \int_0^1 |B_{xx}^{\driftd}(t, \, x) | dt  \leq
      \unpo \quad \forall \, x \in [0, \, L].
\end{equation}
\end{itemize}
{\bf A representation formula for $\boldsymbol{v_2}:$} it is
convenient to introduce the auxiliary function
\begin{equation*}
  V_2 (t, \, x) = \int_0^x
               v_2 (t, \, \xi) d \xi,
\end{equation*}
which satisfies the equation
\begin{equation*}
       V_{2 t} + \lambda_2 V_{2 x} - V_{2 xx} =
       \tilde{S}_1 (t, \, x),
\end{equation*}
where
\begin{equation*}
  \tilde{S}_1 (t, \, x) =
  \int_{ 0}^{x} \tilde{s}_1 (t, \, \xi)
  d \xi .
\end{equation*}
The boundary and initial conditions of $V_2 (t, \, x)$ are
\begin{equation*}
      V_2 (0, \, x) = \int_{0}^{x}
                         v_2 (0, \, \xi) d \xi,
      \qquad
      V_2 (t, \, 0) = \int_0^t (v_{2 x} - \lambda_2 v_2) (s, \, 0) ds,
      \qquad
      V_{2x} (t, \, L) = 0.
\end{equation*}
The convolution kernels \eqref{eq_def_theta} and
\eqref{eq_conv_ker_b} provide the representation formula
\begin{equation}
\label{eq_v_big}
\begin{split}
             V_{2} (t, \, x ) =
      &      \int_0^L {\Theta}^{\lambda_2^{\ast}} (t, \, x, \, y)
             V_2 (0, \, y) dy +
             \int_0^t B (t -s, \, x)
             ( v_{2x} - \lambda_2 v_2 )
             (s,\, 0) ds \\
      &      +
              \int_0^t \int_0^L
              \Theta^{\lambda_2^{\ast}}( t-s, \, x, \, y)
              \Big(
                  \big( \lambda_2^{\ast} - \lambda_2 \big)
                      v_2 \Big)
                    (s, \, y) dy ds \\
      &      +
              \int_0^t \int_0^L
              \Theta^{\lambda_2^{\ast}}( t-s, \, x, \, y)
              \tilde{S}_1
                    (s, \, y) dy ds.  \\
     \end{split}
\end{equation}
Since
\begin{equation*}
      \tilde{\Theta}^{\driftd}(t, \, x, \, 0) +
      B_x^{\driftd}(t, \, x) \equiv 0,
      \qquad
      \tilde{S}_1 (t, \, 0) \equiv 0,
\end{equation*}
from \eqref{eq_v_big} one gets
\begin{equation*}
      \begin{split}
            V_{2x} (t, \, x) =
      &     v_2 (t, \, x) =
            \int_0^L \tilde{\Theta}^{\lambda_2^{\ast}}
            (t, \, x, \, y) v_{2 } (0, \, y) dy +
            \int_0^t B^{\lambda_2^{\ast}}_x (t -s, \, x)
            \Big( v_{2 x} - \driftd v_2 \Big)(s, \, 0) ds \\
      &      +
            \int_0^t \int_0^L \tilde{\Theta}^{\lambda_2^{\ast}}
            (t-s, \, x, \, y) \tilde{s}_1 (s, \, y) dy ds      +
            \int_0^t \int_0^L \tilde{\Theta}^{\lambda_2^{\ast}}
            ( t-s, \, x, \, y)
            \Big(
                  \big( \lambda_2^{\ast} - \lambda_2 \big)
                 v_2
            \Big)_y
            (s, \, y) dy ds \\
      \end{split}
\end{equation*}
and
\begin{equation*}
      \begin{split}
             v_{2x} (t, \, x) =
      &     \int_0^L \tilde{\Theta}_x^{\lambda_2^{\ast}}
            (t, \, x, \, y) v_{2 } (0, \, y) dy +
            \int_0^t B^{\lambda_2^{\ast}}_{xx}
            (t -s, \, x)
            \Big( v_{2 x} - \driftd v_2 \Big)(s, \, 0) ds \\
      &      +
            \int_0^t \int_0^L \tilde{\Theta}_x^{\lambda_2^{\ast}}
            (t-s, \, x, \, y) \tilde{s}_1 (s, \, y) dy ds \\
      &     +
            \int_0^t \int_0^L \tilde{\Theta}_x^{\lambda_2^{\ast}}
            ( t-s, \, x, \, y)
            \Big(
                 \big(
                 \lambda_2^{\ast} - \lambda_2 \big) v_{2y} -
                 \lambda_{2y} v_2
            \Big)
            (s, \, y) dy ds . \\
      \end{split}
\end{equation*}
From the estimate \eqref{eq_theta_x_integral},
\eqref{estimate_theta_x_wrtt} and \eqref{eq_b_xx_wrtt} on the
convolution kernels it follows
\begin{equation*}
\begin{split}
     \int_0^1 | v_{2 x}  (t, \, x)| dt \leq&~
     \| v_2 (0) \|_{L^1} \sup_{x, \, y}
     \int_0^1 |\tilde{\Theta}^{\driftd}(t, \, x, \, y)| dt \\
& ~    +
      \unpo
      \bigg( \int_0^1 \bigg\{ (v_{2 x} - \lambda_2 v_2 )
      (  s, \, 0)  + (\lambda_2^{\ast} - \lambda_2 )
      v_2 (s, \,0 ) ds
      \bigg\}ds \bigg) \\
&~     + \bigg( \int_0^1 |\tilde{s}_1
      (s)|_{\infty} ds \bigg)
      \, \bigg( \int_0^1    \frac{\unpo}{\sqrt{t}} \, dt \bigg)+
      \bigg( \int_0^1 \frac{\unpo}{\sqrt{s}} \, ds \bigg) \,
      \bigg( \delta_1  \sup_{y } \int_0^1 |v_{2 y}|(s, \, y) ds  + \delta_1^2
      \bigg) \\
      \leq&~
      \mathcal{O}(1) \delta_1,
\end{split}
\end{equation*}
for all $x \in [0, L]$.
If $t > 1$ we can use for $v_{2x} $ the expression
\begin{equation}
\label{equation_v2x}
      \begin{split}
             v_{2x} (t, \, x) =
      &     \int_0^L \tilde{\Theta}_x^{\lambda_2^{\ast}}
            (1, \, x, \, y) v_{2} (t-1, \, y) dy +
            \int_0^1 B^{\lambda_2^{\ast}}_{xx}
            (1-s, \, x)
            \big( v_{2 x} - \driftd v_2 \big)(t-1+s, \, 0) ds \\
      &      +
            \int_0^1 \int_0^L  \tilde{\Theta}_x^{\lambda_2^{\ast}}
            (1-s, \, x, \, y) \tilde{s}_1 (t-1+ s, \, y) dy ds \\
      &     +
            \int_0^1 \int_0^L \tilde{\Theta}_x^{\lambda_2^{\ast}}
            ( 1-s, \, x, \, y)
            \Big(
                 \big(
                 \lambda_2^{\ast} - \lambda_2 \big) v_{2y} -
                 \lambda_{2y} v_2
            \Big)
            (t-1+s, \, y) dy ds.\\
      \end{split}
\end{equation}
Computations analogous to the previous ones lead to
\begin{equation*}
   \int_1^T | v_{2x} (s, \, x)| ds \leq \mathcal{O}(1) \delta_1.
\end{equation*}
Hence
\begin{equation*}
      \int_0^T | v_{2x} (s, \, x)| ds \leq \mathcal{O}(1) \delta_1
      \quad \forall \, T \, > 0,  \; \; x \in [0, \, L] .
\end{equation*}
{\bf The integrability of $\boldsymbol{w_2}$ with respect to
time}: it holds
\begin{equation}
\label{estimate_w_wrt}
      \int_0^t |w_{2 }(s, \, y)| ds \leq \unpo \delta_1
      \quad \forall \, t > 0, \quad
      \forall \, y \, \in \, [0, \, L].
\end{equation}
\begin{proof}
We preliminary observe that
\begin{equation*}
      w_2 (0, \, x) = \langle \tilde{\ell}_2, \, u_t (0, \, x)
      \rangle,
      \qquad
      w_2 (t, \, 0) = \langle \tilde{\ell}_2, \,  u'_{b \, 0}(t)
      \rangle,
      \qquad
      w_2 (t, \, L)= \langle \tilde{\ell}_2, \, u'_{b \, L}(t)
      \rangle,
\end{equation*}
where $\tilde{\ell}_2$ satisfies $\langle \tilde{\ell}_2, \, r_2
\rangle \equiv 1$ and $\langle \tilde{\ell}_2, \, \tilde{r}_1
\rangle \equiv 0$. Hence
\begin{equation*}
      \| w_2 ( t =0) \|_{L^1(0, \, L)} \leq \unpo \delta_1,
      \qquad
      \| w_2 (x=0) \|_{L^1(0, \, + \infty)} \leq  \delta_1,
      \qquad
       \| w_2 ( x =L) \|_{L^1(0, \, + \infty)} \leq
       \delta_1.
\end{equation*}
Let $2c$ be the separation speed defined by
\eqref{eq_separation_speed}, let $K$ be a compact neighborhood of
the value $u^{\ast}$ defined by \eqref{E:noBD} and let $C >0$
satisfy
\begin{equation*}
      0 < c \leq \lambda_2(u) \leq C
      \quad \forall \; u \in K.
\end{equation*}
If $y \in ]0, \, L[$, the estimate \eqref{estimate_w_wrt} can be
obtained applying Lemma \ref{functional_pro} to the functional
\begin{equation}
\label{eq_functional_w2_wrttx}
      P_y( x) =
      \left\{
      \begin{array}{lll}
             a \big( 1 - e^{- C x} \big) \quad
             \quad x \leq y \\
             \\
             b \big(  e^{- c x } - e^{-cL} \big)
             \quad x > y,
      \end{array}
      \right.
\end{equation}
where $a$ and $b$ satisfy
\begin{equation}
\label{eq_conditions}
       \left\{
       \begin{array}{lll}
              a \big( 1 - e^{- Cy} \big) =
              b \big(  e^{- cy } - e^{-cL} \big) \\
              \\
              a C e^{-C y} + bc e^{-c y}=1. \\
       \end{array}
       \right.
\end{equation}
By straightforward computations, from \eqref{eq_conditions} one
gets that
the functional $P_y$ satisfies
\begin{equation*}
\begin{split}
&     P_y (0) = P_y (L) =0, \qquad
      0 \leq P_y (x) \leq P_y(y) \leq \unpo,
      \quad
      P'_y (0) \leq \unpo,
      \quad
       - P_y'(L)  \leq \unpo,
      \quad \forall \, L > > 1 \\
&     \qquad \qquad \qquad \qquad  \qquad
      \qquad  P''_y (x) + \lambda_2 P'_y (x) \leq -
      \delta_{x =y}. \\
\end{split}
\end{equation*}
Since $w_2$ satisfies
\begin{equation*}
      w_{2 t} + (\lambda_2 w_2)_x - w_{2 xx} =
      \tilde{s}_2(t, \,x),
\end{equation*}
Lemma \ref{functional_pro} ensures that
\begin{equation*}
\begin{split}
        \int_0^t |w_{2 }(s, \, y)| ds
&       \leq
        \unpo \int_0^L |w_2 (0, \, x)| dx +
        \unpo \int_0^t \int_0^L
        |\tilde{s}_2 (s, \, x)| dx ds \\
&       \quad +
        \unpo \int_0^t |w_2 (s, \, 0)| ds +
        \unpo \int_0^t |w_2 (s, \, L)| ds \\
&       \leq
        \unpo \delta_1
        \quad \forall \, y \in ]0, \, L[. \phantom{\int} \\
\end{split}
\end{equation*}
\end{proof}

{\bf Integrability of $\boldsymbol{w_{2x}}$ with respect to time}:
it holds
\begin{equation}
\label{estimate_wx_wrtt}
      \int_0^t |w_{2 x}(s, \, x)| ds \leq \unpo \delta_1
      \quad \forall \, t>0.
\end{equation}
\begin{proof}
From the representation
\begin{equation*}
      \begin{split}
             w_{2x} (t, \, x) =
      &      \int_0^L \Delta_x^{\lambda_2^{\ast}}
             (t, \, x, \, y) w_2 (0, \, y) dy  +
              \int_0^t \int_0^L \Delta_x^{\lambda_2^{\ast}} (t -s, \, x, \, y)
              \tilde{s}_2 (s, \, y) dy ds       \\
       &       +
             \int_0^{t}  \int_0^L \Delta_x^{\lambda_2^{\ast}}
             (t-s, \, x, \, y)
             \Big( ( {\lambda_2^{\ast}} - \lambda_2 \big) w_{2y} -
                   \lambda_{2 y} w_2
             \Big) (s, \, y) ds dy +
             w_2 (0, \, L) J_x^{\driftd \, L} (t, \, x)
             \\
      &        +
             w_2 (0, \, 0) J_x^{\driftd \, 0}(t, \, x) +
             \int_0^{t}
             J_x^{\lambda_2^{\ast} \, 0}(t - s, \, x)
             w'_{2} (s, \, 0) ds        +
             \int_0^{t}
             J_x^{\lambda_2^{\ast} \, L} (t- s, \, x)
             w'_{2 } (s, \, L) ds   \\
      \end{split}
\end{equation*}
it follows
$$
  \int_0^1 | w_{2x}| (t, \, x) dx \leq \mathcal{O}(1) \delta_1.
$$
If $t \geq 1$ one can write
\begin{equation*}
\begin{split}
      w_{2x} (t, \, x) =
      &      \int_0^L \Delta_x^{\lambda_2^{\ast}}
             (1, \, x, \, y) w_2 (t-1, \, y) dy  +
              \int_0^1 \int_0^L
              \Delta_x^{\lambda_2^{\ast}} (1 -s, \, x, \, y)
              \tilde{s}_2 (t-1+ s, \, y) dy ds       \\
       &       +
             \int_0^{1}  \int_0^L \Delta_x^{\lambda_2^{\ast}}
             (1-s, \, x, \, y)
             \Big( ( {\lambda_2^{\ast}} - \lambda_2 \big) w_{2y} -
                   \lambda_{2 y} w_2
             \Big) (t-1+s, \, y) ds dy \\
        &     +
             w_2 (t - 1, \, L) J_x^{\driftd \, L} (1, \, x)
             +
             w_2 (t-1, \, 0) J_x^{\driftd \, 0}
             (1, \, x) \phantom{\int}\\
        &    +  \int_0^{1}
             J_x^{\lambda_2^{\ast} \, 0}(1 - s, \, x)
             w'_{2} (t-1 +s, \, 0) ds        +
             \int_0^{1}
             J_x^{\lambda_2^{\ast} \, L} (1- s, \, x)
             w'_{2 } (t-1+ s, \, L) ds   \\
      \end{split}
\end{equation*}
and obtains
$$
  \int_1^T |w_{2x}| (t, \, x) dt \leq \mathcal{O}(1) \delta_1.
$$
This concludes the proof of \eqref{estimate_wx_wrtt}.
\end{proof}

\subsubsection{Proof of the estimate \eqref{interaction2}}
\label{interaction2_proof}

We need three preliminary results: \\
$\bullet$ For any $t \leq 1$, the following holds:
\begin{equation}
      \qquad \label{estimate_tilde_theta_x}
      | \tilde{\Theta}_x^{\driftd}(t, \, x,\, y)|
      \leq
      a (t, \, x-y) + b (t, \, x)
      \quad \|a(t)\|_{L^1(-L, \, L)}, \; \;
      \|b(t) \|_{L^1(-L, \, L)}
      \leq \frac{\mathcal{O}(1)}{ \sqrt{t}  }.  \quad
\end{equation}
{\sl Proof of \eqref{estimate_tilde_theta_x}} In the following,
$\alpha(t, \, x - y)$ and $\beta(t, \, x)$ will denote functions
that satisfy
\begin{equation*}
      \| \alpha(t)\|_{L^1(-L, \, L)}, \quad
      \|\beta(t)\|_{L^1(-L, \, L)} \leq \frac{\unpo}{ \sqrt{t}}.
\end{equation*}
By the identities \eqref{Delta_product} and
\eqref{eq_tilde_theta_delta},
\begin{equation*}
      \tilde{\Theta}_x^{\driftd}(t, \, x,\, y) =
      \Delta_x^{\driftd}(t, \, x,\, y) =
    \bigg( \phi^{\, \lambda_2^{\ast}}(t, \, x, \, y)
    \sum_{m \,  = \,- \infty}^{m \,  = \, + \infty}
      G (t, x + 2mL -y) -
      G (t, x+ 2mL +y ) \bigg)_x.
\end{equation*}
One has
\begin{equation*}
\begin{split}
      \bigg|  \bigg(
      \phi^{\, \lambda_2^{\ast}}(t, \, x, \, y)
      \sum_{m \,  = \,- \infty}^{m \,  = \, +
      \infty}
      G (t, x + 2mL -y) \bigg)_x \bigg| \leq
&      \sum_{m \ge 0} \Big| \, G_x^{\driftd}
      (t, \, x-y + 2mL) \Big|     +
      \driftd \sum_{n > 0} G^{\driftd}(t, \, 2nL - x+ y) \\
&      +
      \sum_{n >0}
       G_x^{\driftd}(t, \, 2nL - x+ y)
       \leq \alpha(x -y) ,
      \phantom{\bigg(}  \\
\end{split}
\end{equation*}
where we have set $n: = -m $. To complete the proof of
\eqref{estimate_tilde_theta_x}, it is convenient to observe that
\begin{equation*}
       G^{\driftd}_x (t, \, x + y) \leq G^{\driftd}(t, \, x )
       \quad
       \forall \, x \ge \big( \driftd t + \sqrt{2t}\big),
       \quad \forall \,
       y \ge 0
\end{equation*}
and that
\begin{equation*}
\begin{split}
&      |G^{\driftd}_x (t, \, x + y)| \leq
        G^{\driftd}_x (t, \, x) +
        G_x(t, \, \sqrt{2t})
        \, \chi_{ \displaystyle{
        \{ \, 0 \leq  \, x \leq \sqrt{2t} + \driftd t \}}  }
         \leq \beta(x) \\
&       |G^{\driftd}_x (t, \, 2L - x - y)| \leq
        G^{\driftd}_x (t, \, L - x) +
        G_x (t, \, \sqrt{2t})
        \, \chi_{ \displaystyle{
        \{ \,  L - \sqrt{2t} -  \driftd t \leq \, x \, \leq L \}}  }
         \leq \beta(x), \quad \forall \, x, \; y \in [0, \, L]
\end{split}
\end{equation*}
where $\chi_E$ denotes the characteristic function of the set $E$.
Hence
\begin{equation*}
\begin{split}
      \bigg|  \bigg(
&     \phi^{\, \lambda_2^{\ast}}(t, \, x, \, y)
      \sum_{m \,  = \,- \infty}^{m \,  = \, +
      \infty}
      G (t, x + 2mL + y) \bigg)_x \bigg|
      \leq
      \sum_{m >  0} G_x^{\driftd}
      (t, \, x+ y + 2mL)     +
      G_x^{\driftd}(t, \, x + y) \\
&     \quad +
      \driftd \sum_{n > 0} G^{\driftd}(t, \, 2nL - x - y)  +
      \sum_{n >0}
       G_x^{\driftd}(t, \, 2nL - x - y) \\
&     \leq \sum_{m >0 } G_x^{\driftd}
      (t, \, x+  2mL) + \beta (x ) +
      \driftd \sum_{n > 0} G^{\driftd}(t, \, L - x)   +
      \sum_{n > 1} G_x^{\driftd}(t, \, (2n- 1) L - x )  +
      G_x^{\driftd}(t, \, 2L - x- y) \\
&      \leq \beta(x), \phantom{\sum} \\
\end{split}
\end{equation*}
which concludes the proof of \eqref{estimate_tilde_theta_x}.
$\qquad \qquad \qquad \qquad \qquad \qquad \qquad \qquad \qquad
\qquad \qquad \qquad \qquad \Box$

\noindent $\bullet$ If $t \leq 1$ then
\begin{equation}
\label{estimate_vx_wrtx}
      \int_0^L |v_{ 2 x }(t, \, x)| dx \leq \frac{\unpo
      \delta_1}{\sqrt{t}}.
\end{equation}
\begin{proof} Let $t \leq 1$. From the equality
\begin{equation}
\label{eq_uxx}
       u_{xx}= v_1 \Big( \bigd \tr u_x  + v_{1 x} \tilde{r}_{1 v}
       + \sigma_{1 x} \tilde{r}_{1 \sigma} \Big) + v_{1 x} \tr +
       p_1 \Big( \bigd \hr u_x + p_{1x } \hat{r}_{1 p}
       \Big) + p_{1x } \hr + v_{2 x} r_2 +
       p_{2 x} r_2,
\end{equation}
and from the bounds $\|p_{1 x}(t) \|_{L^1} \leq \unpo \delta_1$
and $\|u_{xx}(t)\| \leq \unpo \delta_1 / \sqrt{t}$, it follows
that
\begin{equation*}
      \|v_{1x}(t)\| = \|\langle \ell_1 , \, u_{xx}(t)
      \rangle  - p_{1 x}(t)\|_{L^1} \leq
      \frac{\unpo \delta_1}{\sqrt{t}},
\end{equation*}
where $\ell_1 = (1, \, 0)$. Hence
\begin{equation*}
       \|w_1(t)\|_{L^1} \leq \unpo \|v_1(t)\|_{L^1}+
       \|v_{1 x}(t)\|_{L^1} + \unpo \| p_1(t)\|_{L^1} +
       \|p_{1 x}(t)\|_{L^1} \leq \frac{\unpo \delta_1}{\sqrt{t}}.
\end{equation*}
From the estimates
\begin{equation*}
\begin{split}
&     \|w'_1 (x= 0)\|_{L^1 (0, \, + \infty)} =
      \|\langle \ell_1, \, u''_{b \, 0} \rangle  \|_{L^1 (0, \, + \infty)}
      \leq \delta_1 \\
&     \|w'_1 (x= L)\|_{L^1 (0, \, + \infty)} =
      \|\langle \ell_1, \, u''_{b \, L} \rangle  \|_{L^1 (0, \, + \infty)}
      \leq \delta_1 \\
&     \|w_1 (t =0)\|_{L^1(0, \, L)} =
      \|\langle \ell_1, \, u''_0 - A(u_0) u'_0  \rangle \|_{L^1(0, \, L)}
      \leq \unpo \delta_1, \\
\end{split}
\end{equation*}
and from the representation formula
\begin{equation}
\label{eq_representation}
\begin{split}
      w_{1 x}(t, \, x) =
&     \int_0^L \Delta^{\lambda_1^{\ast}}_x
      (t,\, x, \, y) w_1 (0, \, y)dy +
      J_x^{\lambda_1^{\ast} \, 0}(t, \, x) w_1(0, \, 0)
      + J_x^{\lambda_1^{\ast} \, L}(t, \, x) w_1(0, \, L) \\
&     + \int_0^t J_x^{\lambda_1^{\ast} \, 0}
     (t - s, \, x) w_1' (s,\, 0)
      ds +
      \int_0^t J_x^{\lambda_1^{\ast} \, L}(t - s, \, x) w_1' (s, \, L)
      ds \\
&     + \int_0^t \int_0^L
      \Delta^{\lambda_1^{\ast}}_x(t -s, \, x, \, y)
      \Big( ( \lambda_1^{\ast}- \lambda_1 ) w_{1 y}
      - \lambda_{1 y} w_1 \Big)
      (s, \, y) ds dy, \\
\end{split}
\end{equation}
it follows that
\begin{equation*}
       \|w_{1 x}(t)\|_{L^1} \leq \frac{\unpo \delta_1}{\sqrt{t}}.
\end{equation*}
Hence
\begin{equation*}
      \|\sigma_{1 x } (t) v_1 (t)\|_{L^1} =
      \bigg\|
      \theta' \bigg( w_{1 x }(t) -
      \frac{w_1}{v_1 } v_{1 x}(t) \bigg)
      \bigg\|_{L^1} \leq \frac{\unpo \delta_1}{\sqrt{t}}.
\end{equation*}
and therefore from \eqref{eq_uxx} one gets
\eqref{estimate_vx_wrtx}.
\end{proof}

\noindent $\bullet$ If $t \ge 1$ then
\begin{equation}
\label{estimate_v2_wrt_x}
      \int_0^L |v_{2 x}(t, \, x))| dx \leq \unpo \delta_1
\end{equation}
\begin{proof} One can repeat the same computations performed to prove
\eqref{estimate_vx_wrtx}, using, instead of
\eqref{eq_representation}, the following representation formula
(which holds if $t \ge 1$):
\begin{equation*}
\begin{split}
      w_{1 x}(t, \, x) =
&     \int_0^L \Delta^{\lambda_1^{\ast}}_x
      (1,\, x, \, y) w_1 (t -1, \, y)dy +
      J_x^{\lambda_1^{\ast} \, 0}(1, \, x) w_1(t-1, \, 0)
      + J_x^{\lambda_1^{\ast} \, L}(1, \, x) w_1(t-1, \, L) \\
&     + \int_0^1 J_x^{\lambda_1^{\ast} \, 0}
     (1 - s, \, x) w_1' (t-1+s,\, 0)
      ds +
      \int_0^1 J_x^{\lambda_1^{\ast} \, L}(1 - s, \, x) w_1' (t-1+s, \, L)
      ds \\
&     + \int_0^1 \int_0^L
      \Delta^{\lambda_1^{\ast}}_x(1 -s, \, x, \, y)
      \Big( ( \lambda_1^{\ast}- \lambda_1 ) w_{1 y}-
      \lambda_{1 y} w_1 \Big)
      (t-1+s, \, y) ds dy. \\
\end{split}
\end{equation*}
\end{proof}

Let
\begin{equation*}
      \mathcal{I}(T):= \sup_{
      \begin{array}{cc}
             \scriptstyle{\tau \in (- T, T)} \\
             \scriptstyle{x \in (-L, \, L)}
      \end{array}}
      \integralt \integralx
      |v_1 (t, \, x)| \, |v_{2 x}(t - \tau, \, x - \xi )| dt dx.
\end{equation*}
It holds:
\begin{equation*}
       \int_0^T
       \int_0^L |v_1(t, \, x)| \,
       |v_{2 x}(t, \, x)| dx dt \leq \mathcal{I}(T).
\end{equation*}
Moreover, thanks to the estimates \eqref{estimate_vx_wrtx} and
\eqref{estimate_v2_wrt_x},
\begin{equation*}
\begin{split}
&     \int_{ \max \{ 0, \, \tau \}}^{ \max \{  2, \, 2 + \tau \} }
      \!\! \integralx
      |v_1(t, \, x)| \, |v_{2 x}(t - \tau, \, x - \xi)|
      \leq
      \unpo\| v_1 \|_{L^{\infty}} \delta_1 \int_0^2
      \bigg\{
      1 + \frac{1}{\sqrt{t}}
      \bigg\} dt \leq
       \unpo \delta_1^2 . \\
\end{split}
\end{equation*}
Hence we are left to estimate the term
\begin{equation*}
      \int_{ \max \{2, \, 2 + \tau \}}^{  \min \{ T, \, T + \tau \}}
      \integralx |v_1 (t, \, x)| \, |v_{2 x}( t - \tau, \, x - \xi)|
      dx dt
\end{equation*}
in the case $T \ge 2$: to do this, we will exploit the
representation formula \eqref{equation_v2x} and the estimate
\eqref{estimate_tilde_theta_x}.

One has
\begin{equation*}
\begin{split}
&     \integraltt \integralx v_1 (t, \, x )
      \int_0^L \tilde{\Theta}_x^{\lambda_2^{\ast}}(1, \, x- \xi, \, y)
      v_2(t -1 - \tau, \, y ) \\
&     \leq \integraltt \integralx v_1 (t, \, x )
      \int_0^L a(1, \, x- \xi- y)
      v_2(t -1 - \tau, \, y )  \\
&     \quad +
      \integraltt \integralx v_1 (t, \, x )
      \int_0^L b(1, \, x - \xi )
      v_2(t -1 - \tau, \, y) \\
&     \leq \int_{-L}^L a(1, \, z)
      \integraltt
      \int_{ \max \{ 0, \, z + \xi \}}^{ \min \{ L, \, L + \xi +z \}}
      v_1 (t, \, x) v_2( t-1 - \tau, \, x - z - \xi) d\xi \\
&     \quad +
      \integralx b(1, \, x - \xi)
      \Bigg( \integraltt v_1(t, \, x)
      \bigg( \int_0^L v_2 (t - 1 - \tau, \,  y ) dy \bigg) dt
      \Bigg) dx \leq \unpo \delta_1^2,
\end{split}
\end{equation*}
and
\begin{equation*}
\begin{split}
&      \integraltt \integralx
       v_1(t, \, x)
       \int_0^1 \int_0^L
       \tilde{\Theta}_x^{\lambda_2^{\ast}}(1-s, \, x - \xi, \,  y)
       \Big( (\lambda_2^{\ast} - \lambda_2) v_{2 y} \Big)
       ( t-\tau -1 +s, \, y) dy ds dx dt \\
&      \leq \delta_1 \integraltt \integralx
       v_1(t, \, x)
       \int_0^1 \int_0^L
       a(1-s, \, x - \xi - y)
       v_{2 y}( t- \tau - 1 +s, \, y) dy ds dx dt \\
&      \quad + \delta_1
       \integraltt \integralx
       v_1(t, \, x)
       \int_0^1 \int_0^L
       b(1-s, \, x - \xi)
       v_{2 y}( t- \tau- 1 +s, \, y) dy ds dx dt \\
&      \leq \delta_1 \int_0^1 \int_{-L}^L
       a(1-s, \, z)
       \Bigg( \int_{\max\{0, \, \xi + z \}}^{ \min \{ L, \, L + z+ \xi \}}
       \integraltt v_1 (t, \, x) v_{2 x} (t - \tau-1 +s, \, x - \xi
       -z)dx dt
       \Bigg) dz ds \\
&      \quad + \delta_1 \int_0^1  \integralx
       b(1-s, \, x - \xi)
       \Bigg( \integraltt v_1(t, \, x)
       \bigg( \int_0^1 v_{2 y} (t - \tau -1+s, \, y) dy \bigg) dt
        \Bigg) dx ds \\
&      \leq \unpo \delta_1 \mathcal{I}(T) + \unpo \delta_1^3.
\end{split}
\end{equation*}
With analogous computations one can estimate the other terms that
comes from the representation formula \eqref{equation_v2x} and
hence prove that $\mathcal{I}(T) \leq \unpo \delta_1^2$.
\subsubsection{Proof of the estimate \eqref{energy_estimates_eq}}
\label{energy_estimates_proof} Since in the following we will
often refer to equations \eqref{decomposition} and
\eqref{eq_theta}, we recall them:
\begin{equation*}
      \sigma_1 = \lambda^{\ast}_1 -
                 \theta \bigg( \frac{w_1}{v_1} +
                              \lambda_1^{\ast}
                        \bigg),
\end{equation*}
where the cut-off $\theta$ is given by
\begin{equation*}
      \theta(s) =
      \left\{
      \begin{array}{lll}
             s \quad \quad \textrm{if} \; |s|\leq \hat{\delta}     \\
             0 \quad \quad \textrm{if} \; |s|\geq 3 \hat{\delta} \\
             \textrm{smooth connection if}
                     \quad \hat{\delta} \leq s \leq 3 \hat{\delta}
      \end{array}
      \right.
      \delta_1 < < \hat{\delta} \leq \frac{1}{3}.
\end{equation*}
It follows that $ | w_1 + \sigma_1 v_1| \neq 0 $
 only when the function $\theta$
is not the identity, i.e. when $| w_1 + \lambda^{\ast}_1 v_1
|>\hat{\delta} |v_1| $. Since
\begin{equation*}
       w_1=v_{1 x} - \lambda_1 v_1 + p_{1 x} - \lambda_1 p_1 ,
\end{equation*}
the condition $ | w_1 + \sigma_1 v_1| \neq 0 $ implies
\begin{equation*}
      |v_{1 x} +( \lambda_1^{\ast} - \lambda_1 ) v_1 +
       p_{1 x} - \lambda_1 p_1 | > \hat{ \delta} |v_1|.
\end{equation*}
There are therefore two possible cases:
\begin{enumerate}
\item
\begin{equation*}
     |v_{1 x} +( \lambda_1^{\ast} - \lambda_1 ) v_1 |
         \ge \frac{1}{2} \hat{\delta} |v_1|,
\end{equation*}
and therefore, since $| \lambda_1^{\ast} - \lambda_1 | \leq
  \mathcal{O}(1) \delta_1$ and $ \delta_1 < < \hat{\delta}$,
\begin{equation*}
  |v_{1 x}| \ge  \frac{\hat{\delta}}{3} |v_1|.
\end{equation*}
\item
\begin{equation*}
      |v_{1 x}| < \frac{\hat{\delta}}{3} |v_1| \quad
                \Longrightarrow            \quad
      | p_{1 x} - \lambda_1 p_1 | > \frac{\hat{\delta}}{2} |v_1|.
\end{equation*}
\end{enumerate}
If case 1 holds, then
\begin{equation*}
\begin{split}
       |w_1 + \sigma_1 v_1|  = &~
          |v_{1 x} +( \sigma_1  - \lambda_1 ) v_1 +
           p_{1 x} - \lambda_1 p_1 | \cr
       \leq&~
       | v_{1 x} | + \delta_1 |v_1| +
           |p_{1 x} - \lambda_1 p_1| \leq
           \mathcal{O}(1)
        |v_{1 x} | + |p_{1 x} - \lambda_1 p_1| \\
\end{split}
\end{equation*}
and therefore
\begin{equation*}
\begin{split}
&              \Big( |v_1|+|w_1|+|v_{1 x}| + |w_{1 x}| \Big)
             \Big( |w_1+ \sigma_1 v_1| \Big)  \leq
             \mathcal{O}(1)
              \Big(
              |v_{1 x}| + |p_1| +|p_{1 x }| + |w_{1 x}| \Big)
              \Big( \mathcal{O}(1)
                |v_{1 x}|+ |p_{1 x}- \lambda_1 p_1| \Big) \\
&                \leq \unpo
                 \Big( |v_{1 x}| +
                 |p_1| +|p_{1 x }| + |w_{1 x}| \Big)
                 |p_{1 x}- \lambda_1 p_1| +
                 \unpo |v_{1 x}|^2 +
                 \unpo  |v_{1 x}| \Big( |p_1 | + |p_{1 x}| \Big) +
                 \unpo |w_{1 x}|^2. \\
\end{split}
\end{equation*}

Since
\begin{equation*}
      |p_1|, \, |p_{1 x}| \leq \unpo \delta_1 \exp (- c x/2),
\end{equation*}
it follows that, if case 1 holds, then one is left to prove
\begin{equation}
\label{real_energy_estimate}
        \int_{0}^{T} \int_{0}^{L}
        \chi_{ \;\{| (w_1 / v_1 ) + \lambda_1^{\ast}| \ge
        \hat{\delta} \}}
               \Big (|v_{1 x}|^2 + |w_{1 x}|^2 \Big) (t,x)
                     dx dt
  \leq
  \mathcal{O}(1) \delta_1 ^2.
\end{equation}
On the other hand, if case 2 holds then
\begin{equation*}
        |v_{1 x} +( \sigma_1  - \lambda_1 ) v_1 +
           p_{1 x} - \lambda_1 p_1 | \leq
            \frac{4}{3} \hat{\delta} |v_1 | + |p_{1 x} - \lambda_1 p_1 | \leq
         \mathcal{O}(1) |p_{1 x} - \lambda_1 p_1|,
\end{equation*}
and therefore
\begin{equation*}
\begin{split}
&       \int_0^T \int_0^L
       \Big( |v_1|+|w_1|+|v_{1 x}| + |w_{1 x}| \Big)
        \Big(|w_1+\sigma_1 v_1| \Big)(s, \, x) ds dx
         \\
&      \quad \leq   \mathcal{O}(1)
       \int_0^T \int_0^L \Big( |v_1|+|w_1|+|v_{1 x}| +
       |w_{1 x}| \Big)
       | p_{1 x} - \lambda_1 p_1|(s, \, x) ds dx \leq \unpo \delta_1^2,
\end{split}
\end{equation*}
thanks to the exponential decay of $
 |p_1|
$ and $
 |p_{1 x}|
$.\\

To prove \eqref{real_energy_estimate} it is convenient to
introduce a new cutoff function:   
\begin{equation*}
      \psi(s) =
      \left\{
      \begin{array}{lll}
            0 \quad \quad \textrm{if} \; |s|\leq 3/5\,  \hat{\delta} \\
            1 \quad \quad \textrm{if} \; |s|\geq 4/5 \, \hat{\delta} \\
            \textrm{smooth connection if} \quad 3/5 \,
            \hat{\delta} \leq |s| \leq 4/5 \, \hat{\delta}.
      \end{array}
      \right.
\end{equation*}
Moreover, in the following we will only prove that
\begin{equation}
\label{real_energy_estimate2}
        \int_{0}^{T} \int_{0}^{L}
        \chi_{ \;\{| (w_1 / v_1 ) + \lambda_1^{\ast}| \ge
        \hat{\delta} \}}
               |v_{1 x}|^2  (t,x)
                     dx dt
  \leq
  \mathcal{O}(1) \delta_1 ^2,
\end{equation}
because the estimate
\begin{equation*}
        \int_{0}^{T} \int_{0}^{L}
        \chi_{ \;\{| (w_1 / v_1 ) + \lambda_1^{\ast}| \ge
        \hat{\delta} \}}
              |w_{1 x}|^2  (t,x)
                     dx dt
  \leq
  \mathcal{O}(1) \delta_1 ^2.
\end{equation*}
can be obtained with similar techniques.

As we have already observed,
 it is sufficient
to show
\begin{equation*}
      \int_{0}^{T} \int_{0}^{L}
               |v_{1 x}|^2 \psi \Big( \argument \Big) (t, x)
                     dx dt
  \leq
  \mathcal{O}(1) \delta_1 ^2.
\end{equation*}
Multiplying the equation
\begin{equation*}
       v_{1 t} + ( \lambda_1 v_1 )_x - v_{1 xx} = 0
\end{equation*}
by $\psi v_1$, we get
\begin{equation}
\label{developement}
\begin{split}
      0 & =  \integral \Bigg( \, \frac{d}{dt} \big( \energy \psi \big) -
          \energy ( \psi_t + \lambda_1 \psi_x -\psi_{xx} )+
          \psi |v_{1 x}|^2
         +  \energy \lambda_{1 x} \psi -
          v_1^2 \psi_{xx} \Bigg)dx dt \\
        & + \int_{0}^{T} \Bigg[ \psi v_1 ( \lambda_1 v_1 - v_{1 x})
        \bigg]^{x=L}_{x=0} dt +
          \int_{0}^{T} \Bigg[ \energy ( \psi_x - \lambda_1 \psi )
          \bigg]^{x=L}_{x=0} dt. \\
\end{split}
\end{equation}
Indeed,
\begin{equation*}
      \frac{d}{dt} \big( \energy \psi \big) = v_1 v_{1 t} \psi + \energy \psi_t
\end{equation*}
and
\begin{equation*}
\begin{split}
       &     \integral \big( \lambda_1 v_1 - v_{1 x} \big)_x \psi v_1 dx dt  =
          \integral ( v_{1 x} - \lambda_1 v_1 ) ( \psi v_1 )_x dx dt +
          \int_{0}^{T} \Bigg[ \psi v_1 ( \lambda_1 v_1 - v_{1 x} )\bigg]^{x=L}_{x=0} dt \\
       &  \quad = \integral \psi_x \bigg( \energy \bigg)_x +
          \psi v_{1 x}^2 -\lambda_1 \psi_x v_1 ^2 -
          \lambda_1 \psi \bigg( \energy \bigg)_x dx dt +
         \int_{0}^{T} \Bigg[ \psi v_1 (  \lambda_1 v_1 - v_{1 x} )\bigg]^{x=L}_{x=0} dt \\
       &  \quad = \integral \psi v_{1 x}^2 + \bigg( \energy \bigg) ( \lambda_{1 x} \psi
          - \lambda_1 \psi_x + \psi_{xx} -2 \psi_{xx} ) dx dt
         + \int_{0}^{T} \Bigg[ \psi v_1 ( \lambda_1 v_1 - v_{1 x} )\bigg]^{x=L}_{x=0}
         dt \\
       &  \quad \quad +  \int_{0}^{T}
       \Bigg[ \energy ( \psi_x - \lambda_1 \psi  )\bigg]^{x=L}_{x=0} dt .\\
\end{split}
\end{equation*}
One can develop the term $
 \psi_t + \lambda_1 \psi_x - \psi_{xx}
$ and, since
\begin{equation}
\label{developement2}
\begin{split}
&
                    \psi_t = \psi' \bigg( \frac{w_{1 t}v_1
                    - w_1 v_{1 t}}{v_1^2}   \bigg),   \qquad   \qquad     
                    \psi_x = \psi' \bigg( \frac{w_{1 x}v_1 -
                    w_1 v_{1 x}}{v_1^2} \bigg), \\
                  & \psi_{xx} = \psi'' \bigg( \frac{w_1}{v_1} \bigg)_x^2 +
                    \psi' \bigg( \frac{w_{1 xx} v_1 - v_{1 xx} w_1}{v_1^2}-
                    2 \frac{ v_{1 x}
                    ( w_{1 x}v_1 -
                    w_1 v_{1 x}) }
                    {v_1^3} \bigg),\\
\end{split}
\end{equation}
one obtains
\begin{equation*}
\begin{split}
 v_1^2(\psi_t + \lambda_1 \psi_x - \psi_{xx})=&~
      \, \psi' v_1 ( w_{1 t} + ( \lambda_1 w_1 )_x - w_{1 xx} ) -
      \psi'  w_1 ( v_{1 t} + ( \lambda_1 v_1 )_x - v_{1 xx} )\\
      & ~ -
       \psi'' v_1^2 \bigg( \frac{w_1}{v_1} \bigg)_x^2 +
           2 \psi' v_{1 x} v_1 \bigg( \frac{w_1}{v_1} \bigg)_x. \\
\end{split}
\end{equation*}
Thus, inserting the last formula into \eqref{developement}, we
obtain
\begin{equation*}
\begin{split}
      \integral \psi |v_{1 x}|^2  =&~
                - \frac{1}{2} \int_0 ^L \bigg[v_1^2 dx
                \bigg]^{t = T}_{t=0}
        + \int_{0}^{T} \bigg[ \psi v_1 (  v_{1 x} - \lambda_1 v_1)
        \bigg]^{x=L}_{x=0} dt +
        \int_{0}^{T} \bigg[ \energy ( \psi_x - \lambda_1 \psi )
        \bigg]^{x=L}_{x=0} dt   \\
     &~ - \frac{1}{2}
         \integral \psi'' v_1^2 \bigg( \frac{w_1}{v_1} \bigg)_x^2 +
                               \psi' v_{1 x} v_1 \bigg( \frac{w_1}{v_1} \bigg)_x +
                               v_1^2 \psi_{xx} -
                               \energy \lambda_{1 x} \psi. \\                                               \\
\end{split}
\end{equation*}
The boundary terms are bounded by $\unpo \delta_1^2$ since $\| v_1
\|_{L^{\infty}} \leq \unpo \delta_1$ and thanks to the estimates
of Proposition \ref{functional_estimates_pro}. Since by
\eqref{length_functional_eq}
\begin{equation*}
       \int_0^T \int_0^L
       \chi_{\{ | \lambda_1^{\ast} + w_1/ v_1 | \leq 3 \hat{\delta}
       \}} v_1^2 \bigg( \frac{w_1}{v_1} \bigg)_x^2 dx ds \leq
       \unpo \delta_1^2,
\end{equation*}
we are left to estimate the following terms:
\begin{equation*}
\begin{split}
\bullet \qquad     \int_0^T & \int_0^L
      \bigg|
      \psi' v_{1 x} v_1 \bigg( \frac{w_1}{v_1} \bigg)_x
      \bigg| ds dx \leq
      \int_0^T \int_0^L
      \bigg|
      \psi' v_{1 x} \bigg( w_{1 x} - \frac{w_1}{v_1} v_{1x}
      \bigg)
      \bigg| ds dx \\
&     \leq
      \unpo \int_0^T \int_0^L
      \bigg|
      \psi '\bigg( |v_1| + |p_{1 x} - \lambda_1 p_1| \bigg)
      \bigg( w_{1 x} - \frac{w_1}{v_1} v_{1x}
      \bigg)
      \bigg| ds dx \\
&     \leq
      \unpo \int_0^T \int_0^L
      \big| v_1 w_{1 x} - v_{1 x} w_1 | ds dx +
      \unpo \int_0^T \int_0^L
      |p_{1 x} - \lambda_1 p_1| \Big(|w_{1 x}| + \unpo | v_{1 x}| \Big)
\end{split}
\end{equation*}
Indeed, if $\psi' \neq 0$ then $|\lambda_1^{\ast} - w_1 / v_1|
\leq \hat{\delta}$ and hence
\begin{equation*}
       |v_{1 x }| \leq \unpo |v_1| + |p_{1 x} - \lambda_1 p_1|.
\end{equation*}
\begin{equation*}
\begin{split}
\bullet \qquad \qquad     \int_0^T \int_0^L
      \psi ' \Big( w_{1 xx} v_1 - w_1 v_{1 xx} \Big) ds dx =&~
      \unpo \int_0^T \int_0^L \Big( w_{1 x} v_1 - w_1 v_{1 x}
      \Big)_x ds dx \\
     \leq&~ \unpo
       \int_0^T \bigg[ w_{1 x} v_1 -
       w_1 v_{1 x} \bigg]^{x =L }_{x=0} \leq
       \unpo \delta_1^2 \\
\end{split}
\end{equation*}
\begin{equation*}
\begin{split}
\bullet \qquad \qquad            \Big| \integral \energy
             \lambda_{1 x} \psi \Big| =&~
             \Big| \integral \energy ( \lambda_{1 } -
             \lambda^{\ast}_1)_x \psi \Big| \\
            \leq&~  \Big| \integral ( \lambda_{1 } -
             \lambda^{\ast}_1) \Big( \energy \psi \Big)_x \Big| +
            \Big| \int_0 ^T \bigg[ ( \lambda_{1 } - \lambda^{\ast}_1)
              \energy \psi  \bigg]^{x = L}_{x = 0}   \qquad \qquad \\
          \leq&~
             \mathcal{O}(1) \delta_1  \Big| \int_0 ^T
             \bigg[ \energy \psi \bigg]^{x = L}_{x = 0} +
             \mathcal{O}(1) \delta_1^2
            \leq \mathcal{O}(1) \delta_1^2. \\
\end{split}
\end{equation*}
\vspace{2cm}

\noindent {\bf Acknowledgments.} The author expresses her
gratitude to Stefano Bianchini for having proposed the problem and
for many helpful suggestions. She also wishes to thank Alberto
Bressan for useful remarks.
\bibliography{biblio}

\end{document}